\documentclass[11pt,twoside]{amsart}

\usepackage[dvipsnames]{xcolor}
\usepackage{hyperref}
\hypersetup{colorlinks=true, linkcolor=Blue, citecolor=Maroon}
\usepackage{amsthm, amsmath, amscd, amssymb,centernot,txfonts}
\usepackage{tikz-cd}
\usepackage{lipsum}
\usepackage{multicol}
\usepackage{amsfonts}
\usepackage{mathrsfs}
\usepackage[all]{xy}
\usepackage[left=2.5cm,top=2.5cm,bottom=3cm,right=2.5cm]{geometry}
\usepackage{dirtytalk}
\usepackage{mathtools}
\usetikzlibrary{shapes,arrows}
\usepackage{etoolbox}
\usepackage{dynkin-diagrams}
\usepackage{enumitem}
\usepackage{chngpage}
\usepackage{adjustbox}
\usepackage[normalem]{ulem}
\usepackage{float}
\usepackage{cancel}

\usepackage[utf8]{inputenc}
\usepackage{fourier} 
\usepackage{array}
\usepackage{makecell}
\usepackage{comment}

\usepackage{hyperref}

\newcommand{\longsquiggly}{\xymatrix{{}\ar@{~>}[r]&{}}}
\newcommand\aimplies{\stackrel{\mathclap{\normalfont\mbox{Thm. ~\ref{Def_1}}}}{\implies}}
\newcommand\bimplies{\stackrel{\mathclap{\normalfont\mbox{Cor. ~\ref{Defphi}}}}{\implies}}
\newcommand\cimplies{\stackrel{\mathclap{\normalfont\mbox{Cor. ~\ref{defcanonicalmorphism }}}}{\implies}}

\theoremstyle{plain}
\numberwithin{equation}{section}
\newtheorem{theorem}{Theorem}[section]
\newtheorem{proposition}[theorem]{Proposition}
\newtheorem{lemma}[theorem]{Lemma}
\newtheorem{corollary}[theorem]{Corollary}

\newtheorem{set-up}[theorem]{Set-up}

\theoremstyle{definition}
\newtheorem{remark}[theorem]{Remark}

\newtheorem{question}[theorem]{Question}
\newtheorem{definition}[theorem]{Definition}

\newcommand*{\QEDA}{\hfill\ensuremath{\blacksquare}}
\newcommand*{\QEDB}{\hfill\ensuremath{\square}}

\tikzstyle{decision} = [diamond, draw, , 
    text width=4.5em, text badly centered, node distance=3cm, inner sep=0pt]
\tikzstyle{block} = [rectangle, draw, , 
    text width=10em, text centered, rounded corners, minimum height=2em]
\tikzstyle{block1} = [rectangle, draw, , 
    text width=5em, text centered, rounded corners, minimum height=2em]    
\tikzstyle{line} = [draw, -latex']
\tikzstyle{cloud} = [draw, ellipse,, node distance=3cm,
    minimum height=2em]
\begin{document}

\title[Deformations and moduli of irregular canonical covers with $K^2=4p_g-8$]{Deformations and moduli of irregular canonical covers with $K^2=4p_g-8$}

\author[P. Bangere]{Purnaprajna Bangere$^1$}
%\address{Department of Mathematics, University of Kansas, Lawrence, USA}
%\email{purna@ku.edu}

\author[F.J. Gallego]{Francisco Javier Gallego$^2$}
%\address{Departamento de \'Algebra, Geometr\'ia y Topolog\'ia and Instituto de Matem\'atica Interdisciplinar,
%Universidad Complutense de Madrid, Spain}
%\email{gallego@mat.ucm.es}

\author[J. Mukherjee]{Jayan Mukherjee$^3$}
%\address{The Institute for Computational and Experimental Research in Mathematics, Providence, USA} 
%\email{jayan\textunderscore mukherjee@brown.edu}

\author[D. Raychaudhury]{Debaditya Raychaudhury$^4$}
%\address{The Fields Institute for Research in Mathematical Sciences, Toronto, Canada}
%\email{draychau@fields.utoronto.ca}

\keywords{\textcolor{black}{Deformations of morphisms, moduli of surfaces of general type, canonical covers or surfaces of minimal degree. 
\\ 
\noindent\rule{2.1cm}{0.3mm}
\\
\\
\noindent 1. Department of Mathematics, University of Kansas, Lawrence, USA.
\\
\noindent purna@ku.edu
\\
\\ 
\noindent 2. Departamento de \'Algebra, Geometr\'ia y Topolog\'ia and Instituto de Matem\'atica Interdisciplinar,
Universidad Complutense de Madrid, Spain.
\\
\noindent gallego@mat.ucm.es \qquad (corresponding author)
\\
\\
\noindent 3. Department of Mathematics, University of California, Riverside, USA
%The Institute for Computational and Experimental Research in Mathematics, Providence, USA.
\\
\noindent jayanm@ucr.edu
\\
\\
\noindent 4. Department of Mathematics, University of Toronto Toronto, Canada.
\\
\noindent debaditya.raychaudhury@utoronto.ca
}}

\maketitle

\begin{center}
    \textit{Dedicated to N. Mohan Kumar on the occasion of his 70th birthday}
\end{center}

\begin{abstract}

In this article, we study the  moduli of irregular surfaces of general type with at worst canonical singularities satisfying $K^2 = 4p_g-8$, for any even integer $p_g\geq 4$.
These surfaces also have unbounded irregularity $q$. We carry out our study by investigating the deformations of the canonical morphism $\varphi:X\to \mathbb{P}^N$, where $\varphi$ is a quadruple Galois cover of a smooth surface of minimal degree. These canonical covers are classified in \cite{GP} into four distinct families, one of which is the easy case of a product of curves.
The main objective of this article is to study the deformations of the other three, non trivial, 
unbounded families. We show that any deformation of $\varphi$ factors through a double cover of a ruled surface and, hence, is never birational. More interestingly, we prove that, with two exceptions, a general deformation of $\varphi$  is two--to--one onto its image, whose normalization is a ruled surface of appropriate irregularity. 
We also show that, with the exception of one family, the deformations of $X$ are unobstructed even though $H^2(T_X)$ does not vanish. Consequently, $X$ belongs to a unique irreducible component of the Gieseker moduli space. These irreducible components are uniruled. As a result of all this, we show the existence of infinitely many moduli spaces, satisfying the strict Beauville inequality $p_g > 2q-4$, with an irreducible component that has a proper "quadruple" sublocus where the degree of the canonical morphism jumps up. 
These components are above the Castelnuovo line, but nonetheless parametrize surfaces with non birational canonical morphisms. The existence of jumping subloci is a contrast with the moduli of surfaces with $K^2 = 2p_g- 4$, studied by Horikawa. Irreducible moduli components with a jumping sublocus also present a similarity and a difference to the moduli of curves of genus $g\geq 3$, for, like in the case of curves, the degree of the canonical morphism goes down outside a closed sublocus but, unlike in the case of curves, it is never birational. \textcolor{black}{Finally, our study shows that there are infinitely many moduli spaces \textcolor{black}{with} an irreducible component \textcolor{black}{whose general elements have} non birational canonical morphism and another irreducible component whose general elements have birational canonical map. }
\end{abstract}

\section{introduction}\label{1}

Canonical covers \textcolor{black}{(i.e., canonical maps which are finite morphisms onto their image)} of varieties of minimal degree have a ubiquitous presence in the geometry of algebraic surfaces and higher dimensional varieties. They appear as extremal cases in a variety of geometric situations. The first, paradigmatic example of a canonical double cover is the canonical morphism of a hyperelliptic curve. In the case of surfaces the works of Enriques and Horikawa (see \cite{En49}, \cite{Hor}) show that minimal surfaces of general type on the Noether line $K_X^2 = 2p_g - 4$ are all canonical double covers of surfaces of minimal degree.
These results show that the deformations of canonical double covers of surfaces of minimal degree are again canonical double covers of surfaces of minimal degree, unlike what happens for canonical double covers in the case of curves of $g\geq 3$. Therefore, the natural, next question to ask 
is: 
\begin{question}\label{question0}
\textcolor{black}{Are there cases in which}  the degree $n$ of  a canonical \textcolor{black}{morphism} 
\begin{equation*}
    \varphi: X \longrightarrow Y \hookrightarrow \mathbb P^N 
\end{equation*} of a surface $Y$ of minimal degree changes, when we deform \textcolor{black}{$\varphi$}? 
\end{question}

If $Y$ is not of minimal degree, there are  some interesting families of examples constructed by Catanese, Beauville, Schreyer, Ciliberto--Pardini--Tovena, Ashikaga--Konno, Gallego--Gonz\'alez--Purnaprajna (\cite{Cat81}, \cite{Cat87}, \cite{CS02},  \cite{CPT00}, \cite{AK90}, \cite{GGP2}, \cite{GGP13c}) providing \textcolor{black}{positive} answers to Question 1.1. However, if $Y$ is a surface of minimal degree, apart from the trivial case of $X$ being a product of curves and the $n=2$ and $n=3$ cases, the answer to Question 1.1 in general is unknown. 
Indeed, Enriques and Horikawa settled the matter for canonical double covers.
Canonical triple covers $X$ of surfaces of minimal degree (these covers satisfy $K_X^2=3p_g(X)-6$) are very few (their geometric genus $p_g(X)$ is bounded by $5$ and their images are singular surfaces; \textcolor{black}{see \cite{Hor76Inv},}  \cite{Hor77}, \cite{Hor82} and \cite{Kon}),  and, when $p_g(X) \leq 4$, their deformations 
are again canonical triple covers of surfaces of minimal degree. Thus, degree $n=4$ is the next case of study and for this we settle Question 1.1 for all cases except one, when $X$ is irregular.
The geometry of canonical quadruple covers of minimal degree (these covers 
satisfy $K_X^2=4p_g(X)-8$) display a wide range of behaviors. Indeed, 
they act like general surfaces of general type from a number of geometric perspectives. 
\textcolor{black}{Quadruple canonical covers} are the first case among low degree  covers where \textcolor{black}{irregular} families appear. Moreover, 
as Remark~\ref{remark.unboundedness} below indicates, \textcolor{black}{they} are the only ones, among canonical covers of smooth surfaces of minimal degree, having both unbounded geometric genus and irregularity (see Theorem~\ref{gpmain}), with the possible exception of degree $6$ covers. All this makes quadruple covers stand out as the most  {interesting} case among canonical covers of surfaces of minimal degree and are natural \textcolor{black}{(and non trivial, except, obviously, when $X$ is a product of curves, see Type $(3)_m$ of Theorem~\ref{gpmain})}  candidates  for testing Question \ref{question0}.

\begin{remark}\label{remark.unboundedness}
{{It follows from {a more general result (see \cite{GPtriple}, Theorem 3.2)}, that there are no odd degree canonical covers of smooth surfaces of minimal degree {other than $\mathbb P^2$}.}}  
{This together with \cite{Beau}  implies} that, if $\chi(X) \geq 31$, then the degree of  a canonical cover of a smooth surface of minimal degree could only be $2,4,6$ or $8$ (if $\chi(X) \leq 30$, then $q(X)$ is bounded). 
 Since the irregularity of degree $8$ canonical covers {is} bounded above by $3$ when {$p_g\geq 115$} (see \cite{Xia}),
 degree $4$ canonical covers are the only ones, {among covers of smooth surfaces of minimal degree,} having unbounded irregularity (and thus {unbounded geometric genus}, since $p_g(X)\geq 2q(X)-4$ by \cite{BeauD}) except possibly the degree six canonical covers. One can show that there are no smooth regular degree $6$ abelian covers of smooth surfaces of minimal degree (see \cite{BGMR20}).    Evidence seems to indicate that there are no such irregular covers as well.
\end{remark}

\textcolor{black}{In moduli problems, it is usual to choose that special member whose deformations describe the general element of the moduli component.} In this article we focus on the study of the deformations of irregular quadruple {Galois} canonical covers of smooth surfaces of minimal degree. \textcolor{black}{Tables~\ref{table1} and \ref{table2} show that} the deformations of quadruple Galois covers capture the full complexity of the situation and provide a very interesting answer to Question \ref{question0}.   \textcolor{black}{This} shows that there \textcolor{black}{is} no need to deal with quadruple covers at large, even though \textcolor{black}{we show that} Galois covers do deform to the so called non Galois natural covers. We completely figure out \textcolor{black}{their} behavior under deformations for all but one family \textcolor{black}{Theorem~\ref{gpmain} (}for \textcolor{black}{that} one family, we have some partial results\textcolor{black}{)}. 
\textcolor{black}{Thus our results make} quadruple canonical covers of surfaces of minimal degree the lowest degree covers, with the possible exception of triple covers of geometric genus $5$, providing a positive answer to Question  \ref{question0}. Regular quadruple Galois canonical covers of surfaces of minimal degree also provide {positive} answers to Question \ref{question0}, as the authors will show in a forthcoming  article. From all this we derive interesting consequences for the moduli of surfaces of general type.

\subsection{Classification of irregular quadruple Galois canonical covers of surface scrolls.}
The classification  of irregular quadruple Galois canonical covers of surfaces of minimal degree was done by the first two authors in \cite{GP}. We need the technical details of their classification results, for the purpose of this article, so we will summarize them here. The image of these covers are smooth rational normal scrolls $Y$. Recall that a smooth rational normal scroll is a Hirzebruch surface $\mathbb{F}_e$ ($e\geq 0$), which is, by definition, $\mathbb{P}(\mathcal{E})$, where $\mathcal{E}=\mathcal{O}_{\mathbb{P}^1}\oplus\mathcal{O}_{\mathbb{P}^1}(-e)$. Let $p:\mathbb{P}(\mathcal{E})\to\mathbb{P}^1$ be the natural projection. The line bundles on $\mathbb{F}_e$ are of the form $\mathcal{O}_Y(aC_0+bf)$ where $\mathcal{O}_Y(C_0)=\mathcal{O}_{\mathbb{P}(\mathcal{E})}(1)$ and $\mathcal{O}_Y(f)=p^*\mathcal{O}_{\mathbb{P}^1}(1)$. The line bundle $\mathcal{O}_Y(aC_0+bf)$ is very ample if $b\geq ae+1$. 
\begin{theorem}\label{gpmain}
(\cite{GP}, Theorem $0.1$) Let $X$ be an irregular surface
with at worst canonical singularities and let $Y$ be a smooth surface of
minimal degree. If the canonical bundle of $X$ is ample and base-point-free, and $\varphi: X\to Y$ 
is a quadruple Galois canonical cover, then $Y$ is the Hirzebruch surface $\mathbb{F}_0$, embedded by $|C_0 + mf|$, $(m \geq 1)$. Let $G$
be the Galois group of $\varphi$.
\begin{itemize}
    \item[(a)] If $G = \mathbb{Z}_4$, then $\varphi$ is the composition of two double covers $p_1:X_1\to Y$
branched along a divisor $D_2$ and $ p_2:X\to X_1$, branched along the ramification
of $p_1$ and $p_1^*D_1$, where $D_1$ is a divisor on $Y$ and with trace zero module
$p_1^*\mathcal{O}_Y(-\frac{1}{2}D_1-\frac{1}{4}D_2)$.
    \item[(b)] If $G = \mathbb{Z}_{2}^{\oplus 2}$, then $X$ is the fiber product over $Y$ of two double covers of $Y$
branched along divisors $D_1$ and $D_2$, and $\varphi$ is the natural morphism from
the fiber product to $Y$.
\end{itemize}
More precisely, $\varphi$ has one of the sets of invariants shown in the following table.
Conversely, if $\varphi:X\to Y$ is either
\begin{itemize}
    \item[(1)] the composition of two double covers $p_1:X_1\to Y$, branched along a divisor
$D_2$, and $p_2:X\to X_1$, branched along the ramification of $p_1$ and $p_1^*D_1$, and with trace zero module $p_1^*\mathcal{O}_Y(-\frac{1}{2}D_1-\frac{1}{4}D_2)$, with $D_1$ and $D_2$ as described
in rows $1$ of the table below; or 
\item[(2)] the fiber product over $Y$ of two double covers $p_1:X_1\to Y$ and $p_2:X_2\to Y$,
branched respectively along divisors $D_2$ and $D_1$, as described in rows $2$, $3$,
and $4$ of the table below,
\end{itemize}
 then $\varphi:X \to Y$ is a Galois canonical cover whose Galois group is $\mathbb{Z}_4$ in case $1$ and
$\mathbb{Z}_2^{\oplus 2}$ in case $2$.
\vspace{5pt}

\begin{center}
 \begin{tabular}{c|c|c|c|c|c|c} 
 \hline
 Type & $p_g(X)$ & $Y$ & $G$ & $D_1\sim$ & $D_2\sim$ & $q(X)$\\ [0.5ex] 
 \hline\hline
 $(1)_m$ & $2m+2$ & $\mathbb{F}_0$ & $\mathbb{Z}_4$ & $(2m+4)f$ & $4C_0$ & $1$\\
 \hline 
 $(1')_m$ & $2m+2$ & $\mathbb{F}_0$ & $\mathbb{Z}_2^{\oplus 2}$ & $2C_0+(2m+4)f$ & $4C_0$ & $1$\\
 \hline 
 $(2)_m$ ($m\geq 2$)& $2m+2$ & $\mathbb{F}_0$ & $\mathbb{Z}_2^{\oplus 2}$ & $(2m+2)f$ & $6C_0+2f$ & $m$\\
 \hline 
 $(3)_m$ & $2m+2$ & $\mathbb{F}_0$ & $\mathbb{Z}_2^{\oplus 2}$ & $(2m+4)f$ & $6C_0$ & $m+3$\\
 \hline
 \end{tabular}
\end{center}
\end{theorem}

\textcolor{black}{Surfaces of type $(3)_m$ are clearly products of a genus $2$ curve and a genus $m+1$ hyperelliptic curve, so their study is easier.  Surfaces 
of types $(1)_m$, $(1')_m$ and $(2_m)$ are not product of curves (for one thing, \textcolor{black}{they satisfy the equality $K_X^2=4p_g(X)-8$ but} not the equality $p_g(X) = 2q(X)-4$; see \cite{BeauD}) so studying them is much more complex.} Surfaces of type $(1')_m$, $(2)_m$ and $(3)_m$
are smooth (in fact all surfaces of type $(3)_m$ are smooth) while surfaces of type $(1)_m$ are necessarily singular, the general ones having only $A_1$ singularities. 
The \textcolor{black}{above} table shows the existence of families of quadruple Galois canonical covers with \textcolor{black}{unbounded geometric genus and} unbounded irregularity. 
In addition, some of the families of Theorem~\ref{gpmain} are extremal cases for several inequalities concerning irregular surfaces of general type, such as  $K_X^2\geq 2\chi(X)$, 
the slope inequality (see \cite{LP}), $K_X^2\geq 2p_g(X)$ (see \cite{Deb}) and $p_g(X)\geq 2q(X)-4$ (see \cite{BeauD}). Because of all of this, quadruple Galois covers are interesting from the perspective of the geography of irregular surfaces as well. 

\begin{remark}
{Although the covers of Theorem~\ref{gpmain} are 
\emph{simple iterated double covers} in the sense of \cite{Ma1}, they 
are not \emph{good sequences} \textcolor{black}{(see again \cite{Ma1}, Definition C)}.  
Moreover, the point of view of our article is to study the deformation of canonical morphisms  to projective spaces, rather than the deformation of finite morphisms between two surfaces.
Therefore our study distinctly differs from the study 
of simple iterated double covers carried out by Catanese and Manetti.}
\end{remark}

\subsection{Statements of the main results.} 
Let $X$ be a surface as in Theorem ~\ref{gpmain}, and let $\varphi$ be the canonical morphism of $X$. First we present a description of our results about the algebraic formally semiuniversal deformation space of $\varphi$ (which exists by Remark ~\ref{can}) in the following table (see Theorem ~\ref{A}).
\begin{small}
\begin{table}[H] 
\caption{Deformations of $\varphi$}
\begin{center}\phantomsection\label{table1}
 \begin{tabular}{c|c|c|c} 
 \hline
 \makecell{$X$ is of type} & \makecell{Degree of {\bf any}\\ deformation\\ of ${\varphi}$} & \makecell{Description of $\varphi_t$\\ for  a {\bf general} $t$ in the \\ deformation space of $\varphi$} & \makecell{Normalization of  {the} image \\ of $\varphi_t$ for a general $t$ in the \\ deformation space of $\varphi$} \\  
 \hline\hline 
  $(1)_m$ & $\geq 2$ & \makecell{Double cover onto\\  a non-normal surface} & \makecell{ Elliptic ruled surface \\ with invariant $0$}  \\
 \hline
$(1')_m$ & $\geq 2$ & \makecell{Double cover onto\\ a non-normal surface} & \makecell{Elliptic ruled surface \\ with invariant $0$}  \\
 \hline 
 $(2)_2$ & $\geq 2$ & \makecell{Double cover onto \\ a smooth surface} & \makecell{Ruled surface over \\ a curve of genus $2$ \\with invariant $-2$}  \\
 \hline 
  $(2)_m$ $(m\geq 3)$ & $\geq 2$ & \makecell{\textcolor{black}{factors through a double cover} \\ \textcolor{black}{onto a smooth surface}}  & \makecell{\textcolor{black}{Ruled surface over} \\ \textcolor{black}{a curve of genus $m$}  \textcolor{black}{or $\mathbb{P}^1\times\mathbb{P}^1$}}  \\
 \hline 
 $(3)_1$ & $4$ & \makecell{Quadruple cover onto\\ a smooth surface} & $\mathbb{P}^1 \times \mathbb{P}^1$   \\
 \hline
$(3)_m$ $(m\geq 2)$ & $\geq 2$ & \makecell{Double cover onto \\ a smooth surface} & \makecell{Product of \textcolor{black}{$\mathbb P^1$} with a curve of genus $m+1$}  \\
 \hline
\end{tabular}
\end{center}
\end{table}
\end{small}

\normalsize

\vspace{5pt}

The structure of the covers described in Theorem ~\ref{gpmain} easily
implies the existence of pencils in $X$:

\begin{remark}\label{fibration}
(See also \cite{GP01}, Remark 3.4) Let $X$ be as in Theorem ~\ref{gpmain}. 
\begin{itemize}
     \item[(1)] If $X$ is of type $(1)_m$ or $(1')_m$, then $X$ contains a rational pencil of genus $3$ curves and an irrational pencil (over an elliptic curve) of genus $m+1$ hyperelliptic curves. 
     \item[(2)] If $X$ is of type $(2)_m$, then $X$ contains a rational pencil of genus $2m+1$ curves and an irrational pencil (over a genus $m$ curve) of genus $2$ curves.
     \item[(3)] If $X$ is of type $(3)_m$, \textcolor{black}{since} $X$ is a product of genus $2$ curve and a genus $m+1$ hyperelliptic curve, \textcolor{black}{then} it  contains, obviously, two irrational pencils of genus $2$ 
     and genus $m+1$ hyperelliptic curves.
     \end{itemize}
The deformation of families $(1)_m$, $(1')_m$, $(2)_m$ show the persistence of irrational pencils under deformation. However, unlike in the case of surfaces with $K_X^2=2p_g-4$, 
for surfaces of type $(1)_m$, $(1')_m$ and $(2)_1$, the rational pencils vanish under deformation.
\end{remark}

It is illustrative to compare the results summarized in the above table with the deformations of lower degree canonical covers of surfaces of minimal degree.  As already mentioned, any deformation of a degree 2 canonical morphism is again of degree 2. More generally, deformations of double or triple canonical covers of embedded projective bundles over $\mathbb{P}^1$ of arbitrary dimension are, respectively, of degree 2 and 3 (see \cite{GGP13}, \cite{GGP3} and \cite{GGP4}). Thus Theorem ~\ref{A} is in sharp contrast with these results 
{and, as pointed out before, quadruple  covers are the lowest degree examples} of canonical covers  for which the degree of the canonical map of a general deformation drops down, 
with the possible exception of degree $3$ covers with $p_g=5$.

Now we summarize our results on the moduli (see Theorems ~\ref{B} and ~\ref{C} and Remark~\ref{remark.Ops}). For a given 
$m$, there is a unique component of moduli of surfaces of general type that contains {\it all} surfaces of type $(1)_m$ and $(1')_m$ (see Theorem ~\ref{B}). There is a unique moduli component that contains all surfaces of type $(2)_2$ and, for given $m$, there is also a unique moduli component that contains all surfaces of type $(3)_m$.
In the following table we outline the description of these 
moduli components:

\vspace{-5pt}

\begin{small}

\begin{table}[H]
\caption{Moduli of $X$}
\begin{center}
\phantomsection\label{table2}
\begin{tabular}{c|c|c|c|c} 
 \hline
 \makecell{$X$ is of type} & \makecell{Obstructions to\\ deformations\\ of $X$} & \makecell{Geometry of \\ unique moduli \\ component $\mathcal{M}_{[X]}$\\ containing $X$} &\makecell{Canonical morphism \\of a general surface \\ in $\mathcal{M}_{[X]}$} & \makecell{Normalization of \\ the image of the\\ general canonical \\ morphism} \\  
 \hline\hline 
$(1)_m$ or  $(1')_m$ & Unobstructed & \makecell{Uniruled of \\ dimension \\$8m+20$} & \makecell{Double cover onto\\  a non-normal surface} & \makecell{ Elliptic ruled surface \\ with invariant $0$}  \\
  \hline 
 $(2)_2$ & Unobstructed & \makecell{Uniruled of \\ dimension \\$28$} & \makecell{Double cover onto \\ a smooth surface} & \makecell{Ruled surface over \\ a curve of genus $2$ \\with invariant $-2$}\\
 \hline
\makecell{$(3)_1$} & Unobstructed & \makecell{Uniruled of \\ dimension \\$6$} & \makecell{Quadruple cover onto \\ $\mathbb{P}^1\times\mathbb{P}^1$} & $\mathbb{P}^1\times\mathbb{P}^1$\\
 \hline
\makecell{$(3)_m$} & Unobstructed & \makecell{Uniruled of \\ dimension \\$3m+3$} & \makecell{Double cover onto \\ $\mathbb{P}^1\times C_2$ \\ $g(C_2) = m+1$} & $\mathbb{P}^1\times C_2$\\
 \hline
\end{tabular}
\end{center}
\end{table}
\end{small}

\vspace{5pt}

\smallskip

The fact that the surfaces $X$ as described in the above table lie on a unique component of the moduli of surfaces of general type is a consequence of the unobstructedness of $X$. We prove the unobstructedness of $X$ in these cases despite $H^2(T_X)$ being non zero. 
\smallskip

The results summarized in Table \ref{table2}  have further consequences for the moduli  of surfaces of general type. 

\begin{corollary}
For each one of the moduli spaces $\mathcal M_{(p_g,q,K^2)}=\mathcal M_{(2m+2,1,8m)}$, $m \geq 1$ and  $\mathcal M_{(6,2,16)}$, 
there exists a reduced, irreducible uniruled component for which the degree of the canonical morphism jumps up, from $2$ to $4$, at a proper locally closed "quadruple" sublocus. 
The points of this "quadruple" sublocus correspond to surfaces  whose canonical morphism is a quadruple cover of a surface of minimal degree.
The surfaces parametrized by these moduli spaces are not product of curves. 
\end{corollary}

In this sense, quadruple canonical covers (unlike, for instance, Horikawa surfaces), resemble hyperelliptic curves of genus bigger than $2$  and the existence of a proper hyperelliptic sublocus in the moduli of these curves. In Remark~\ref{remark.dimensions.loci} we estimate the 
dimension of the quadruple loci and give the dimension of the smaller loci parametrizing surface of general type whose canonical morphism is a quadruple Galois cover of a surface of minimal degree.  
\smallskip 

\color{black}
Recall that outside the region defined by 
Castelnuovo's inequality $K_X^2 \geq 3p_g + q-7$ the canonical map of a minimal surface of general type is not birational. Our results (see Subsection 4.4) 
show the existence of infinitely many irreducible components \emph{inside} the region defined by Castelnuovo's inequality $K_X^2 \geq 3p_g + q-7$:

\begin{corollary}\label{corollary.above.Castelnuovo}
For any $m \geq 1$, the moduli spaces $\mathcal M_{(p_g,q,K^2)}=\mathcal M_{(2m+2,1,8m)}$, $\mathcal M_{(2m+2,m,8m)}$ and \linebreak 
$\mathcal M_{(2m+2,m+3,8m)}$, have  an irreducible 
component that pa\-ra\-me\-trizes surfaces whose canonical map is a non birational morphism. 
\end{corollary}
\begin{remark}
Except for 
$\mathcal M_{(4,4,8)}$, the invariants of the moduli spaces of 
Corollary~\ref{corollary.above.Castelnuovo} satisfy 
$K_X^2 \geq 3p_g + q-7$. Among these moduli spaces, only  $\mathcal M_{(2m+2,m+3,8m)}$, $m \geq 1$, lie on the line $p_g=2q-4$, hence parametrize products of curves (see \cite{BeauD}).
\end{remark}
As a consequence of this and the results of Ashikaga (see \cite[Theorem 3.2]{Ashikaga}) we get the following corollary 
(see Corollary~\ref{corollary.two.components}):
\begin{corollary}
There are infinitely many 
moduli spaces \textcolor{black}{with} an irreducible component  \textcolor{black}{whose general elements have} non birational canonical morphism and another irreducible component whose general elements have birational canonical map. 
\end{corollary}

The classification of quadruple canonical covers shows the existence of fibrations in all genus, as is illustrated in Remark~\ref{fibration}. Thus, the results in this article apply to fibrations of all genus. Deformation and moduli of genus two fibrations have been studied in \cite{Sei95} and \cite{GGP13a}. Therefore, for the special case of genus two fibrations, namely, the surfaces of type $(1)_1$, $(1')_1$ and $(2)_2$, their results on unobstructedness and persistence of genus two fibrations upon deformation do apply. But in this article we show uniruledness not only of these moduli components but also of the moduli components of fibrations of all genera. \textcolor{black}{In the final section of the article, we show that infinitesimal Torelli holds for smooth families of type $(1')_m$. It is known that for $(3)_m$ it holds for $m=1$ and does not hold for $m>2$ (see Remark \ref{producttorelli}). We end that section with an interesting question for families of type $(2)_m$.}

\smallskip

\textcolor{black}{In \cite{GGP2} a \textcolor{black}{new} framework, connecting deformations of morphisms and the smoothing of multiple structures on algebraic varieties, \textcolor{black}{was developed. There,} a general criterion on when a finite morphism deforms to a one to one map \textcolor{black}{was found}.  
\textcolor{black}{Even though this started a different way of looking at several interesting, natural situations concerning {the moduli of} surfaces of general type and higher dimensional varieties,  the criterion in \cite{GGP2} is only the first step in proving results in that realm.}
One such natural situation is the moduli of \textcolor{black}{surfaces of general type with} $K^2=4p_g-8$, where we need to bring in new ideas and methods, to study the deformation of the canonical map. \textcolor{black}{This} in turn leads to \textcolor{black}{the} description of moduli components of these surfaces. We describe these \textcolor{black}{new ideas and methods} briefly below.} 
\color{black}
\subsection{Brief sketch of the arguments} We outline the main ideas behind our proofs. We study the deformations of the canonical morphism $\varphi: X \xrightarrow{4:1} Y \hookrightarrow \mathbb{P}^N$. Such a deformation need not be a composition of the deformations of the two morphisms. This brings into play the concepts of existence of multiple structures on a scheme and their smoothings (making our study markedly different from deformations of finite covers). The proof in general is comprised of the following five steps:

\noindent\underline{Step 1.} We show in almost all cases (excepting family $(2)_m$, $m \geq 2$) that the morphism $\varphi_1$ in one of the factorizations $X \xrightarrow[2:1]{\pi_1} X_1 \xrightarrow[2:1]{\varphi_1} Y \hookrightarrow \mathbb{P}^N$ of $\varphi$ (see \cite{GP}) can be deformed to a finite birational morphism  
\begin{equation*}
\begin{tikzcd}
\mathscr{X}_{1t} \arrow[d, "1:1"] \arrow[r] & \mathscr{X}_1 \arrow[d]  & X_1 \arrow[d, "2:1"] \arrow[l] & X\arrow[l, "\pi_1"]\\
\mathbb{P}^N \arrow[r] & \mathbb{P}_{T}^N & \mathbb{P}^N \arrow[l] &
\end{tikzcd}
\end{equation*}
For this we relate the existence and the deformations of double structures to the deformations of $\varphi_1$.

%\smallskip

\noindent\underline{Step 2.} We show that one can complete the above diagram into

\begin{equation*}
\begin{tikzcd}
\mathscr{X}_t \arrow[r] \arrow[d, "2:1"] & \mathscr{X} \arrow[d] & X \arrow[d, "2:1"] \arrow[l]   \\
\mathscr{X}_{1t} \arrow[d, "1:1"] \arrow[r] & \mathscr{X}_1 \arrow[d]  & X_1 \arrow[d, "2:1"] \arrow[l] \\
\mathbb{P}^N \arrow[r] & \mathbb{P}_{T}^N & \mathbb{P}^N \arrow[l]
\end{tikzcd}
\end{equation*}
This requires showing smoothness of the forgetful map between functors ${\bf Def}_{\pi_1} \to {\bf Def}_{X_1}$. This shows that the quadruple cover deforms to a double cover.

\smallskip

\noindent\underline{Step 3.} We then show that any deformation $\mathscr{X} \to \mathbb{P}_T^N$ factors as $\mathscr{X} \to \mathscr{X}_1 \to \mathbb{P}_T^N$ . This requires the smoothness of the forgetful map of functors $\bf Def_{\pi_1/\mathbb{P}^N} \to \bf Def_{\varphi}$ where $\bf Def_{\pi_1/\mathbb{P}^N}$ are deformations of $\pi_1$ relative to $\mathbb{P}^N$ as introduced by Flenner. This shows that the morphism does not further deform to a birational morphism.

\smallskip

\noindent\underline{Step 4.} To prove the unobstructedness of $\varphi$ and $X$, 
we first we show that ${\bf Def}_{\varphi_1}$ is smooth. Then we show that, the following chain of implications hold:
\begin{small}
\begin{center}
${\bf Def}_{\varphi_1}$ is smooth $\qquad\aimplies\qquad$ ${\bf Def}_{\pi_1/\mathbb{P}^N}$ is smooth $\qquad\bimplies\qquad$ ${\bf Def}_{\varphi}$ is smooth $\qquad\cimplies\qquad$ ${\bf Def}_X$ is smooth.
\end{center}
\end{small}
It is to be noted that we show that the varieties are unobstructed inspite of $H^2(T_X) \neq 0$. Finally we show the uniruledness of the algebraic formally semiuniversal deformation space of ${\bf Def}_{\pi_1/\mathbb{P}^N}$, which we construct in Theorem~\ref{Def_1}. From this, we prove the uniruledness of the moduli component of $X$. 

\smallskip

\noindent\underline{Step 5.} There is a subtle point, though, that needs to be addressed: since our surfaces are irregular, a general deformation of the canonical morphism need not be the canonical morphism. We overcome this difficulty and show nevertheless the existence of a deformation of $\varphi$ which is two-to-one and, indeed, canonical.

\color{black}
\subsection{Organization of this article.} In Section ~\ref{2}, we prove the main results we need to carry out our study. In \S ~\ref{prelim 1} we recall some basic facts about the normal sheaf, and in \S ~\ref{prelim 2} we recall the basics of abelian covers. In \S ~\ref{prelim 3}, we describe the techniques of deformations of finite morphisms developed by the first two authors. In \S ~\ref{prelim 4}, we describe how to apply this technique in our situation, namely to study the covers we are interested in. Finally \S ~\ref{prelim 5}  is devoted to study the geometry of the moduli components. We prove the results stated in Table \ref{table1} (see Theorem ~\ref{A}) in Section ~\ref{3}. The proof of the results stated in Table \ref{table2} (see Theorems ~\ref{B} and ~\ref{C}) appear in Section ~\ref{4}. Finally, we prove that the infinitesimal Torelli theorem holds for smooth surfaces of type $(1')_m$ in Section ~\ref{5}. {The results show that there are families for which infinitesimal Torelli hold, and families for which they do not, the latter being well known for the case of product of curves, which is one of the families.}

\medskip
\noindent\textbf{Acknowledgements.} \textcolor{black}{We are grateful to Edoardo  Sernesi for generously giving \textcolor{black}{us} his time for helpful discussions and suggestions which improved the exposition \textcolor{black}{of this paper}.} We thank  Fabrizio Catanese for pointing out to us papers containing examples in the context of Question ~\ref{question0} when $Y$ is not of minimal degree.

\smallskip
The first author thanks
by the General Research Fund of the University of Kansas for partially supporting his research.
The second author was partially supported
by Spanish Government grant PID2021-124440NB-I00 and by Santander-UCM grant PR44/21. The third author was supported by the National Science Foundation, Grant No. DMS-1929284 while in residence at ICERM in Providence, RI, as part of the ICERM Bridge program. The research of the fourth author was partially supported by a Simons postdoctoral fellowship from the Fields Institute for Research in Mathematical Sciences.

\newpage 
%\bigskip
\begin{center}
    \textit{Notations and conventions.}\phantomsection\label{notation}
\end{center}
\begin{itemize}[leftmargin=0.7cm]
  \item[1.] We will always work over the field of complex numbers $\mathbb{C}$ and a variety is an integral separated scheme of finite type over $\mathbb{C}$.
  \item[2.] The symbol `$\sim$' denotes linear equivalence and `$\equiv$' denotes numerical equivalence.
  \item[3.] We will use the multiplicative and the additive notation of line bundles interchangeably. Thus, for line bundles $L_1,L_2$, $L_1\otimes L_2$ and $L_1+L_2$ are the same. $L^{-r}$, $L^{\otimes -r}$ (or $-rL$) denotes $(L^{\vee})^{\otimes r}$.
  \item[4.] If $L_i$ is a line bundle on the variety $X_i$ for $i=1,2$, $L_1\boxtimes L_2$ is by definition, the line bundle $p_1^*L_1\otimes p_2^*L_2$ on $X_1\times X_2$ where $p_i:X_1\times X_2\to X_i$ is the $i$-th projection for $i=1,2$. When $X_i=\mathbb{P}^1$ for $i=1,2$, then $\mathcal{O}_{\mathbb{P}^1\times\mathbb{P}^1}(a,b):=\mathcal{O}_{\mathbb{P}^1}(a)\boxtimes\mathcal{O}_{\mathbb{P}^1}(b)$.
  \item[5.] For a morphism $X\to Y$ between algebraic schemes, $\Omega_{X/Y}$ (or $\Omega_{X/Y}^1$) is the sheaf of relative differentials and $T_{X/Y}=\mathcal{H}\textrm{om}(\Omega_{X/Y},\mathcal{O}_X)$ is the relative tangent sheaf. By convention, $\Omega_X$ (or $\Omega_X^1$) and $T_X$ is obtained by taking $Y=\textrm{Spec}(\mathbb{C})$.
  \item[6.] For an algebraic scheme $X$, let $T^i(X)$ and $\mathcal{T}^i(X)$ be the local and global cohomology of the cotangent complex respectively. It is well known that when $X$ is a projective variety, the first cotangent sheaf $\mathcal{T}_{X}^1=\mathcal{E}\textrm{xt}^1(\Omega_{X},\mathcal{O}_X)$ and $T^1(X)=\textrm{Ext}^1(\Omega_X,\mathcal{O}_X)$. 
  \item[7.] For an algebraic scheme, ${\bf Def}_X$ (resp. \textbf{Def}$'_X$) is the functor of deformations (resp. locally trivial deformations) of $X$.
  \item[8.] For an algebraic scheme $X$ and a line bundle $L$ on it, ${\bf Def}_{(X,L)}$ (resp. ${\bf Def}'_{(X,L)}$) is the functor of deformations (resp. locally trivial deformations) of the scheme and the line bundle i.e. the pair $(X,L)$.
  \item[9.] For a morphism $\varphi:X\to Z$ between algebraic schemes, ${\bf Def}_{\varphi}$ (resp. ${\bf Def}'_{\varphi}$) is the functor of deformations (resp. locally trivial deformations) of $\varphi$ with fixed target.
  \item[10.] For morphism $X\xrightarrow{\pi} Y\to Z$ of algebraic schemes, ${\bf Def}_{\pi/Z}$ (resp. ${\bf Def}'_{\pi/Z}$) is the functor of $Z$-deformati-\\ons (resp. locally trivial $Z$-deformations) of $\pi$ with varying target.
\end{itemize}
For a scheme $X$, if ${\bf Def}_X={\bf Def}_X'$, then ${\bf Def}_{(X,L)}={\bf Def}_{(X,L)}'$, ${\bf Def}_{\varphi}={\bf Def}_{\varphi}'$ and ${\bf Def}_{\pi/Z}={\bf Def}'_{\pi/Z}$ respectively in the situation of 8, 9 and 10. Indeed, by definition (see e.g. \cite{AC} \S 3, and \cite{Ser} \S 3.4.2),  deformations of $(X,L)$, $\varphi$ and $\pi/Z$ are locally trivial if and only if the induced deformation of $X$ is locally trivial.

\section{\textcolor{black}{Results on deformations of finite morphisms}}\label{2}
In this section, we provide the main results regarding the deformations of morphisms that are essential to carry out our study. 
\subsection{Preliminaries on normal sheaves}\label{prelim 1} Locally trivial deformation theory of a morphism is governed by the \textit{normal sheaf} that we define below.
\begin{definition}
  (\cite{Ser}, 3.4.5) To a morphism $f:X\to Z$ of algebraic schemes there exists an exact sequence of coherent sheaves which defines the sheaf $\mathcal{N}_f$ called the {\it normal sheaf} of $f$;
 \begin{equation*}
     0\to T_{X/Z}\to T_X\to f^*T_Z\to \mathcal{N}_f\to 0.
 \end{equation*}
  The morphism $f$ is called {\it non-degenerate} if $T_{X/Z}=0$.
\end{definition}
A morphism being non-degenerate is equivalent to being unramified in an open dense set. Thus, a finite flat morphism between normal Cohen-Macaulay varieties is non-degenerate and so is the composition of non-degenerate morphisms between normal Cohen-Macaulay varieties. The following is the general version of \cite{Go}, Lemma 3.3 whose proof we omit.
\begin{lemma}\label{exactgonzalez}
Let $X,Y,Z$ be normal Cohen-Macaulay varieties. Let $\pi:X\to Y$ be a non-degenerate morphism for which $\pi^*$ is an exact functor (this happens if $\pi$ is finite and flat) and let $\psi:Y\to Z$ be a non-degenerate morphism. Suppose $\varphi:=\psi\circ\pi$. Then there is an exact sequence;
\begin{equation*}
    0\to\mathcal{N}_{\pi}\to\mathcal{N}_{\varphi}\to\pi^*\mathcal{N}_{\psi}\to 0.
\end{equation*}
\end{lemma}

\subsection{Normal abelian covers of smooth varieties}\label{prelim 2} Our objects of study are canonical morphisms that factor through abelian covers. We recall some basic facts about these covers, see \cite{Pa} for further details.

\begin{definition}\label{defgalois}
 Let $Y$ be a variety and let $G$ be a finite abelian group. A {\it Galois cover of $Y$ with Galois group $G$} is a finite flat morphism $\pi:X\to Y$ together with a faithful action of $G$ on $X$ that exhibits $Y$ as a quotient of $X$ via $G$. 
\end{definition}

Let $\pi:X\to Y$ be a Galois $G$ cover of a smooth variety $Y$ with $X$ normal. Then $\pi_*\mathcal{O}_X$ splits as a direct sum indexed by the characters. More precisely, $$\pi_*\mathcal{O}_X=\bigoplus\limits_{\chi\in G^*}L_{\chi}^{-1}.$$ Let $D$ be the branch divisor of $\pi$. Let $\mathcal{C}$ be the set of cyclic subgroups of $G$ and for all $H\in\mathcal{C}$, denote by $S_H$ the set of generators of the group of characters $H^*$. Then, we may write $$D=\sum\limits_{H\in\mathcal{C}}\sum\limits_{\psi\in S_H}D_{H,\psi}$$ where $D_{H,\psi}$ is the sum of all the components of $D$ that have inertia group $H$ and character $\psi$. The sheaves $L_{\chi}$ and the divisors $D_{H,\psi}$ are called the \textit{building data} of the cover. For every pair $\chi,\chi'\in G^*$, for every $H\in \mathcal{C}$ and for every $\psi\in S_H$, one may write $\chi|_H=\psi^{i_{\chi}}$ and $\chi'|_H=\psi^{i_{\chi'}}$, $i_{\chi},i_{\chi'}\in\{0,\dots,m_H-1\}$ where $m_H$ is the order of $H$. Let 
$S_{\chi} = \{(H,\psi): \chi|_H \neq \psi^{m_H-1}\}$.

\subsection{A generalization of Atiyah exact sequence} 
For any scheme $X$, there is a natural map $\mathcal{O}_X^*\to\Omega_X$ defined by $u\mapsto du/u$. For every line bundle $L$ on $X$, the natural map induced between the cohomology groups $H^1(\mathcal{O}_X^{*})\to H^1(\Omega_X)\cong \textrm{Ext}^1(\mathcal{O}_X,\Omega_X)$ gives an extension 
$$e_L:\quad 0\to \Omega_X\to \mathcal{Q}_L\to \mathcal{O}_X\to 0.$$
We set $\mathcal{E}_L:=\mathcal{H}\textrm{om}(\mathcal{Q}_L,\mathcal{O}_X)$, and we obtain the following exact sequence which is known as the Atiyah exact sequence when $X$ is smooth:
\begin{equation}\label{atiyahg}
    0\to\mathcal{O}_X\to\mathcal{E}_L\to T_X\to 0.
\end{equation}
Altmann and Christophersen has generalized \cite{Ser}, Theorem 3.3.11 and showed that $H^1(\mathcal{E}_L)$ parametr-\\izes first order locally trivial deformations of $(X,L)$ when $X$ is reduced.

\begin{theorem}\label{ac} (\cite{AC}, Theorem 3.1 (ii)) If $X$ is a reduced projective scheme, then
${\bf Def}'_{(X,L)}(\mathbb{C}[\epsilon])=H^1(\mathcal{E}_L)$.
\end{theorem}

In fact, one can follow the treatment of \cite{Ser}, Section 3.3.4 to define a map 
$M:\mathcal{E}_L\to H^0(L)^{\vee}\otimes L$ that fits into the following commutative diagram where the left vertical map is the one obtained in ~\eqref{atiyahg}.
\begin{equation}
\begin{tikzcd}
\mathcal{O}_X \arrow[r,"m"]\arrow[d]  & H^0(L)^{\vee}\otimes L\arrow[d, equal] \\
\mathcal{E}_L \arrow[r,"M"] & H^0(L)^{\vee}\otimes L
\end{tikzcd}
\end{equation}          
The proof of the following follows by repeating the argument of \cite{Ser}, Proposition 3.3.14 word by word.

\begin{proposition}\label{sectionlifting}
Let $X$ be a reduced projective scheme and let $L$ be a line bundle on $X$. Assume $(\mathcal{X},\mathcal{L})$ is a first order locally trivial deformation of $(X,L)$ defined by a cohomology class $\eta_1\in H^1(\mathcal{E}_L)$ according to Theorem ~\ref{ac}. Let $s\in H^0(L)$ be a section of $L$. Then $s$ lifts to a section $\mathcal{L}$ if and only if $s\in \textrm{ker}(M_1(\eta_1))$, where $M_1:H^1(\mathcal{E}_L)\to \textrm{Hom}(H^0(L),H^1(L))$ is induced by $M$.
\end{proposition}

When $L$ is base point free, in order to verify the above section-lifting-criterion, we will make use of \cite{Ser}, (3.39) diagram with exact rows and columns, which is given for smooth case, and whose existence is a routine computation when $X$ is reduced and projective. In particular, we have the exact sequence
\begin{equation}\label{atiyah'g}
    0\to \mathcal{E}_L\to H^0(L)^{\vee}\otimes L\to \mathcal{N}_{\varphi_L}\to 0,
\end{equation}
where $\varphi_L:X\to \mathbb{P}(H^0(L)^{\vee})$ is the morphism induced by $|L|$.

\subsection{Techniques to reduce the degree along a deformation}\label{prelim 3} One of the central techniques for deforming a finite morphism to a morphism of smaller degree is to construct a suitable multiple structure on the image of the morphism which are called ropes.

\begin{definition}\phantomsection\label{defropes}
Let $Y$ be a reduced connected scheme and let $\mathcal{E}$ be a vector bundle of rank $m-1$ on $Y$. A {\it rope of multiplicity $m$ on $Y$ with conormal bundle $\mathcal{E}$} is a scheme $\widetilde{Y}$ with $\widetilde{Y}_{\textrm{red}}=Y$ such that $\mathcal{I}_{Y/\widetilde{Y}}^2=0$, and $\mathcal{I}_{Y/\widetilde{Y}}=\mathcal{E}$ as $\mathcal{O}_Y$ modules. If $\mathcal{E}$ is a line bundle then $\widetilde{Y}$ is called a {\it ribbon} on $Y$.
\end{definition}

\begin{remark}\label{extendedribbons}
Recall that, for a morphism $\varphi:X\to\mathbb{P}^N$ from a smooth, projective variety which is finite onto a smooth variety $Y\hookrightarrow\mathbb{P}^N$, the space $H^0(\mathcal{N}_{\varphi})$ parametrizes the space of infinitesimal deformations of $\varphi$. Suppose $\mathcal{E}$ is the trace zero module of the induced morphism $\pi:X\to Y$. It is shown in \cite{Go}, Proposition 2.1, that the space  $H^0(\mathcal{N}_{Y/\mathbb{P}^N}\otimes\mathcal{E})$ parametrizes the pairs  $(\widetilde{Y},\widetilde{i})$ where $\widetilde{Y}$ is a rope on $Y$ with conormal bundle $\mathcal{E}$ and $\widetilde{i}:\widetilde{Y}\to\mathbb{P}^N$ is a morphism that extends $i$. 
\end{remark}

The relation between these two cohomology groups is given by the following proposition.

\begin{proposition}\label{propngonzalez}
 (\cite{Go}, Proposition 3.7) Let $X$ be a normal Cohen-Macaulay projective variety and let $\varphi: X \to \mathbb{P}^N$ be a morphism that factors as $\varphi = i \circ \pi$, where $\pi$ is a finite cover of a smooth variety $Y$
and $i:Y\hookrightarrow\mathbb{P}^N$ is an embedding. Let $\mathcal{E}$ be the trace zero module of $\pi$ and let $\mathcal{I}$ be the ideal sheaf of $i(Y)$. There exists a homomorphism
\begin{equation*}
    H^0(\mathcal{N}_{\varphi})\xrightarrow{\psi}\textrm{Hom}(\pi^*(\mathcal{I}/\mathcal{I}^2),\mathcal{O}_X)
\end{equation*}
that appears when taking cohomology on the commutative diagram \cite{Go} $(3.3.2)$. Since
\begin{equation*}
    \textrm{Hom}(\pi^*(\mathcal{I}/\mathcal{I}^2),\mathcal{O}_X)=H^0(\mathcal{N}_{Y/\mathbb{P}^N})\oplus H^0(\mathcal{N}_{Y/\mathbb{P}^N}\otimes\mathcal{E}),
\end{equation*}
the homomorphism $\psi$ has two components;
\begin{equation*}
    H^0(\mathcal{N}_{\varphi})\xrightarrow{\psi_1}H^0(\mathcal{N}_{Y/\mathbb{P}^N})\textrm{ and }H^0(\mathcal{N}_{\varphi})\xrightarrow{\psi_2}H^0(\mathcal{N}_{Y/\mathbb{P}^N}\otimes\mathcal{E}).
\end{equation*}
\end{proposition}

We will make use of the following fundamental theorem of the deformation theory of finite morphisms to reduce the degree of a general deformation of the canonical morphism.

\begin{theorem}\label{degreductionggp} (\cite{GGP2}, Theorem 1.4)
Let $X$ be a smooth irreducible projective variety and let $\varphi:X\to \mathbb{P}^N$ be a morphism that factors through an embedding $Y\hookrightarrow\mathbb{P}^N$ with $Y$ smooth and let $\pi:X\to Y$ be the induced morphism which we assume to be finite of degree $n\geq 2$. Let $\widetilde{\varphi}:\widetilde{X}\to\mathbb{P}^N_{\Delta}$ $(\Delta=\textrm{Spec}\left(\frac{\mathbb{C}[\epsilon]}{\epsilon^2}\right))$ be a first order infinitesimal deformation of $\varphi$ and let $\nu\in H^0(\mathcal{N}_{\varphi})$ be the class of             
 $\widetilde{\varphi}$.  If
\begin{itemize}
    \item[(a)] the homomorphism $\psi_2(\nu)$ has rank $k> \frac{n}{2}-1$, and 
    \item[(b)] there exists an algebraic formally semiuniversal deformation of $\,\,{\bf Def}_{\varphi}$ and ${\bf Def}_{\varphi}$ is smooth,
\end{itemize}
then there exists a flat family of morphisms, $\Phi: \mathcal{X}\to\mathbb{P}^N_T$ over $T$, where
$T$ is a smooth irreducible algebraic curve with a distinguished point $0$, such that
\begin{enumerate}
    \item[(1)] $\mathcal{X}_t$ is a smooth irreducible projective variety,
    \item[(2)] the restriction of $\Phi$ to the first order infinitesimal neighbourhood of $0$ is $\widetilde{\varphi}$, and 
    \item[(3)] for $t\neq 0$, $\Phi_t$ is finite and one-to-one onto its image in $\mathbb{P}^N$.
\end{enumerate}
\end{theorem}
We will use the theorem above to study the deformations of the canonical morphisms of the varieties we are interested in. However, we are also interested in the degree of the canonical morphisms of the moduli components of these varieties. We remark that a general deformation of the canonical morphism of a regular variety remains canonical by \cite{GGP2}, Lemma 2.4 (the statement requires smoothness, but it holds for varieties with canonical singularities as well, see \cite{BCG21}, proof of Proposition 2.9). 
\par

\subsection{Deformations of iterated double covers of embedded varieties}\label{prelim 4} Throughout this subsection, we will work with the following diagram where $X$, $Y$ and $Z$ are normal local complete intersection (abbreviated as lci) projective varieties, $i:Z\hookrightarrow\mathbb{P}^N$ is an embedding, and $\varphi:X\to Z$ is a twice iterated double cover:
\[\begin{tikzcd}
X\arrow[r, "\pi"]\ar[rr,out=-30,in=200,swap,"\varphi"] & Y\ar[r,"p"] & Z\hookrightarrow\mathbb{P}^N
\end{tikzcd}
\]
We set $\psi:=i\circ p$, and $\varphi:=\psi\circ\pi$. Our objective is to determine the degree of a general deformation of $\varphi$. We will show that under suitable hypotheses, $\varphi$ can be deformed to a two-to-one morphism onto its image. We first need the following technical fact that we will put as a remark for future reference. 

\begin{remark}\label{jayan} Let $Y$ be a normal lci projective variety. Let $L$ be a line bundle on $Y$ and $B$ be a divisor in $|L|$. Assume $H^1(L) = 0$. Let $f:\mathcal{Y} \to T$ be a deformation of $Y$ over a smooth affine pointed variety $(T,0)$ ($f$ is assumed to be proper and flat). Assume that $L$ lifts to a line bundle $\mathcal{L}$ on $\mathcal{Y}$. Then (possibly after shrinking $T$) $f_*(\mathcal{L})$ is locally free of rank $h^0(L)$ on $T$. We have a Cartesian diagram as shown below.
\[ 
\begin{tikzcd}
\mathcal{Y} \times_T \mathbb{P}(f_*(\mathcal{L})) \arrow[r,"p"] \arrow[d,"q"] & \mathbb{P}(f_*(\mathcal{L})) \arrow[d,"g"] \\
\mathcal{Y} \arrow[r,"f"] & T 
\end{tikzcd}
\]
Consequently,  $p: \mathcal{Y} \times_T \mathbb{P}(f_*(\mathcal{L})) \rightarrow \mathbb{P}(f_*(\mathcal{L}))$ is a deformation of $Y$ with $q^*(\mathcal{L})$ and the incidence divisor $\mathcal{B} \in |q^*(\mathcal{L})|$ giving natural lifts of $L$ and $B$ respectively on $\mathcal{Y} \times_T \mathbb{P}(f_*(\mathcal{L}))$. Now since $f_*(\mathcal{L})$ is locally free, by shrinking $T$, we can always construct a section $s: T \to \mathbb{P}(f_*(\mathcal{L}))$ and $\mathcal{B} \times_{\mathbb{P}(f_*(\mathcal{L}))} T$ is a lift of the divisor $B$ to $\mathcal{Y}$. Conversely any lift $\mathcal{B}$ of $B$ on $\mathcal{Y}$ is obtained by a pullback induced by a section $s: T \to \mathbb{P}(f_*(\mathcal{L}))$.
\end{remark}

\begin{remark}\label{trace zero module for canonical morphism}
Let $X$ and $Y$ be normal lci projective varieties with $Y$ smooth and consider morphisms $X \xrightarrow{\pi} Y \xrightarrow{\psi} \mathbb{P}^N$. Let $\pi$ be a finite $2:1$ morphism with trace zero module $L^*$ on $Y$. Let $\varphi = \psi \circ \pi$. Let $L^* = \omega_Y \otimes \psi^*(\mathcal{O}_{\mathbb{P}^N}(-1))$. Then $\varphi^*(\mathcal{O}_{\mathbb{P}^N}(1)) = \omega_X$.
\end{remark}

\color{black}
\begin{theorem}\label{def to 2:1}

Let $X$ be a normal lci projective variety and let $Y$ and $Z$ be smooth projective varieties. Let $\pi:X\to Y$ be a finite, flat morphism of degree two onto $Y$ with trace zero module $\mathscr{E}_{\pi} = L^*$ and branched along a divisor $B \in |L^{\otimes 2}|$, and let $p:Y\to Z$ be a finite (hence flat) morphism of degree two onto $Z$ with trace zero module $\mathcal{E}_p$. Let $i:Z\hookrightarrow\mathbb{P}^N$ be an embedding, $\varphi=i\circ p\circ\pi$ and $\psi=i\circ p$. Suppose 
\begin{itemize}
   \item[(a)] $H^2(\mathcal{O}_{Y})=0$,
   \item[(b)] ${\bf Def}_{\psi}$ is smooth,
   \item[(c)] $\psi_2:H^0(\mathcal{N}_{\psi})\to H^0(\mathcal{N}_{Z/\mathbb{P}^N}\otimes \mathcal{E}_{p})$ is non-zero,
   \item[(d)] $H^1(L^{\otimes 2}) = 0$,
\end{itemize}
then there exist a flat family $\mathcal{X}\to T$ of deformations of $X$ over a smooth pointed affine algebraic curve $(T,0)$ and a $T$-morphism $\Phi:\mathcal{X} \to \mathbb{P}^N_T $ satisfying:
\begin{itemize}
    \item[(1)] $\Phi=\Psi\circ\Pi$, where $\mathcal{Y}\to T$ is a flat family, $\Psi:\mathcal{Y}\to\mathbb{P}^N_T$, and $\Pi:\mathcal{X}\to \mathcal{Y}$ are $T$-morphisms with $\mathcal{Y}_0=Y$, $\Pi_0=\pi$, and $\Psi_0=\psi$, 
    \item[(2)] $\Pi_t$ is a finite morphism of degree 2 for all $t$, and $\Psi_t$ is birational onto its image for all $t\in T-\{0\}$,
    \item[(3)] Suppose that $\varphi$ is the canonical morphism of $X$ and that $\mathcal{E}_{\pi} = L^* = \omega_Y \otimes \psi^*(\mathcal{O}_{\mathbb{P}^N}(-1))$, then $\Phi_t$ can be taken to be the canonical morphism of $\mathcal{X}_t$.
\end{itemize}      
\end{theorem}
\noindent\textit{Proof.} We will prove the assertions (1) and (2) in two steps.\par 
\noindent\underline{\textit{Step $1$.}} In this step, we deform $\psi$ into a birational morphism.  
Notice that ${\bf Def}_{\psi}$ is unobstructed, and has an algebraic formally semiuniversal deformation by \cite{BGG20}, Proposition 1.5. Moreover, $\psi_2$ is non-zero, hence we apply Theorem ~\ref{degreductionggp} and we get that there exists a family $\mathcal{Y}$ of smooth projective varieties, proper and flat over a smooth pointed affine algebraic curve $(T,0)$ and a $T$-morphism $\Psi :\mathcal{Y} \to \mathbb{P}^N_T $ with:
\begin{itemize}
    \item[(1)] $\Psi_0=\psi$, 
    \item[(2)] $\Psi_t$ is birational onto its image for all $t\in T-\{0\}$, 
\end{itemize}      
\noindent\underline{\textit{Step $2$.}} We construct a  deformation $\mathcal{X} \to \mathcal{Y} \to \mathbb{P}_T^N \to T$ of $\varphi$. For this we need to construct a deformation $\Pi: \mathcal{X} \to \mathcal{Y} \to T$ of the finite morphism $\pi: X \to Y$. Let $q: \mathcal{Y} \to T$ be the deformation obtained by applying the forgetful map to $\mathcal{Y} \to \mathbb{P}^N_T \to T$. We need to construct lifts $\mathcal{L}^{\otimes 2}$ and $\mathcal{B}$ of the line bundle $L^{\otimes 2}$ and the divisor $B$ respectively on $\mathcal{Y}$. Note that since $H^2(\mathcal{O}_Y) = 0$, we have that the map ${\bf Def}_{(Y,L)} \to {\bf Def}_{Y}$ is smooth by (\cite{Ser}, Proposition $2.3.6$). Hence by \cite{Ser}, Proposition $2.2.5$, (iv), we have a lift $\mathcal{L}$ of $L$ on $\mathcal{Y}$. The conclusion follows from Remark ~\ref{jayan}. This proves statements $(1)$ and $(2)$.\par 
For part $(3)$ we note that $\omega_{\mathcal{Y}/T} \otimes \Psi^*(\mathcal{O}_{\mathbb{P}^N_T}(-1))$ is a deformation of $L^*$, since $L^* = \omega_Y \otimes \psi^*(\mathcal{O}_{\mathbb{P}^N}(-1))$. Hence $\omega_{\mathcal{Y}/T}^{- 2} \otimes \Psi^*(\mathcal{O}_{\mathbb{P}^N_T}(2))$ is a lift of $L^{\otimes 2}$. Now we apply Remark ~\ref{jayan} to construct a relative cover using a lift $\mathcal{B} \in |\omega_{\mathcal{Y}/T}^{- 2} \otimes \Psi^*(\mathcal{O}_{\mathbb{P}^N_T}(2))|$ of the divisor $B \in |L^{\otimes 2}|$. Since for each $t$, the trace zero module $L_t^* = \omega_{Y_t} \otimes \Psi_t^*(\mathcal{O}_{\mathbb{P}_t^N}(-1))$, we have that for each $t$, $\Phi_t$ is the canonical morphism of $X_t$ by Remark ~\ref{trace zero module for canonical morphism}.     \QEDB\par 

\vspace{5pt}

\begin{remark} \label{generic degree of canonical morphism}
Let $X$ be a projective variety with at worst canonical singularities, and with ample and base point free canonical bundle $\omega_X$. Let $T$ be a smooth affine curve and $\mathcal{Y} \xrightarrow{\Phi} \mathbb{P}^N_T \to T$ be  {such that} for all $t \in T$, $\Phi_t$ is given by the complete linear series $\omega_{\mathcal{Y}_t}$.
Suppose that the degree of the finite morphism $\Phi_t$ is $d$ for a general $t \in T$. Then there exists an irreducible component $U_X$ of the universal deformation space of $X$ such that for a general closed point $u \in U_X$, the canonical morphism of the fibre $\mathcal{X}_u$ of the universal family over $U_X$ has degree less than or equal to $d$.  
\end{remark}
\noindent\textit{Proof.} Choose an irreducible component $U_X$ containing $T$. Since $T$ is smooth this embedding factors through the reduced induced structure $U_X^0$ of $U_X$. Consider the pullback $\mathcal{X} \xrightarrow{p} U_X^0$ of the universal family to $U_X^0$. Since $p_g(\mathcal{X}_s)$ is constant for $s \in U_X^0$ (see 
 \cite{Kol}, Theorem 1),
we have that $p_*(\omega_{\mathcal{X}/U_X^0})$ is locally free of rank $h^0(\omega_X)$. Then $\mathcal{X} \xrightarrow{\Psi} \mathbb{P}(p_*(\omega_{\mathcal{X}/U_X^0})) \to U_X^0$ is a deformation of the canonical morphism of $X$ such that for each $s \in U_X^0$, $\Psi_s$ is the canonical morphism of $\mathcal{X}_s$. Since $U_X^0$ is integral, we have that degree of $\Psi_s$ is upper semicontinuous. Now $\mathcal{Y} \xrightarrow{\Phi} \mathbb{P}^N_T \to T$ is obtained by pulling back $\mathcal{X} \xrightarrow{\Psi} \mathbb{P}(p_*(\omega_{\mathcal{X}/U_X^0})) \to U_X^0$ to $T$ by the embedding $T \hookrightarrow U_X^0$. This shows that degree of $\Psi_s$ is less than or equal to $d$ for a general $s \in U_X^0$. Since closed points of $U_X$ are the same as closed points of $U_X^0$, we are done. \QEDA

\vspace{5pt}

Now we will find the condition following the proof of \cite{Weh}, Proposition 1.10, under which any deformation of $\varphi$ factors through a deformation of $\pi$ (with varying target).

\begin{proposition}\label{corweh}
Let $\pi:X\to Y$ be a finite, flat  morphism with trace zero module $\mathcal{E}$ between projective varieties with $X$ normal lci, and $Y$ smooth. Let $\psi:Y\to Z$ be a non-degenerate morphism to a smooth projective variety $Z$. Let $\varphi=\psi\circ\pi$ be the composed morphism. Assume $H^0(\mathcal{N}_{\psi}\otimes\mathcal{E})=0$. Then the natural map between the functors ${\bf Def}_{\pi/Z}\to{\bf Def}_{\varphi}$ is smooth.
\end{proposition}
\noindent\textit{Proof.}  We apply \cite{Weh}, Proposition 1.10 to the following commutative diagram.
\begin{equation*}\label{eq.comm.square.triangle}
		\xymatrix@C-20pt@R-12pt{
	X  \ar[rr]^{\pi} 
	\ar[dr]_{\varphi} 
	&& Y \ar[dl]^{\psi} \\
	& Z}
\end{equation*}
The maps $\beta_1$ and $\beta_2$ of \cite{Weh}, Proposition 1.10
become the following: $$\beta_1: H^0(\mathcal{N}_{\psi}) \to H^0(\pi^*\mathcal{N}_{\psi})\cong H^0(\mathcal{N}_{\psi})\oplus H^0(\mathcal{N}_{\psi}\otimes\mathcal{E}),$$
$$\beta_2: H^1(\mathcal{N}_{\psi}) \to H^1(\pi^*\mathcal{N}_{\psi})\cong H^1(\mathcal{N}_{\psi})\oplus H^1(\mathcal{N}_{\psi}\otimes\mathcal{E}).$$ The assertion follows since the map $\beta_2$ is always injective and $\beta_1$  is surjective if $H^0(\mathcal{N}_{\psi}\otimes\mathcal{E})=0$.\QEDB\par 

\vspace{5pt}

The following is the main result that we will use to prove Theorem ~\ref{A}. The proof of this result is an immediate consequence of Theorem ~\ref{def to 2:1} and Proposition ~\ref{corweh}.

\begin{corollary}\label{maincor}
Assume the hypotheses (a), (b), (c) and $(d)$ of Theorem ~\ref{def to 2:1}. Furthermore, assume {\bf Def}$_{\varphi}$ has an algebraic formally semiuniversal deformation space and $H^0(\mathcal{N}_{\psi}\otimes\mathcal{E})=0$. Then a general deformation of $\varphi$ is a composition of a double cover over a  deformation $Y'$ of $Y$ followed by a morphism of $Y'\to \mathbb{P}^N$ that is birational onto its image, consequently, it is a two-to-one morphism onto its image.
\end{corollary}

\begin{remark}\label{can}
Let $X$ be a surface with ample and globally generated canonical bundle $\omega_X$ with at worst canonical singularities. Let $\varphi$ be the canonical morphism of $X$. Then $H^0(T_X)=0$. Furthermore, {\bf Def}$_X$, {\bf Def}$_{\varphi}$, and {\bf Def}$_{(X,\omega_X)}$ have algebraic formally universal deformation spaces.
\end{remark}

\begin{proposition}\label{lt-any}
Let $X$ be a normal lci projective variety. Assume $\pi:X\to Y$ is a double cover of a smooth projective variety $Y$ with ramification divisor $R$ and branch divisor $B\in |L^{\otimes 2}|$. Then we have the following exact sequence
\begin{equation}\label{eq.lt-any}
0\to \mathcal{N}_{\pi}\to  \pi^*\mathcal{O}_Y(B)|_R \to \mathcal{T}_X^1 \to 0.
\end{equation}
\color{black}
\end{proposition}

\noindent\textit{Proof.} Let $Y':=\mathbb{V}(L^{-1}):=\textrm{Spec}(\textrm{Sym}(L^{-1}))$ denote the total space of the line bundle $L$. Let $p': Y' \to Y$ be the projection. We have an embedding of $i: X \hookrightarrow Y'$ as a divisor in $Y'$ such that $\pi = p' \circ i $. The conormal sheaf $\mathcal{N}^*_{X/Y'}:=\mathcal{I}/\mathcal{I}^2$ of $X$ in $Y'$ is given by $\pi^*(\mathcal{O}_Y(-B))$ (since $X$ is defined as the scheme of zeroes of $t^{2}-p'^*\varsigma$ where $t \in p'^*(L)$ and $B$ is the zero locus of $\varsigma$). Since $X$ is a local complete intersection, we have an exact sequence 
\begin{equation}\label{ses1}
    0 \to \pi^*(\mathcal{O}_Y(-B)) \to \Omega_{Y'} \otimes \mathcal{O}_X \to \Omega_X \to 0.
\end{equation}  
Now since $p': Y' \to Y$ is a smooth morphism, we have that $\Omega_{Y'/Y}$ is an invertible sheaf isomorphic to $p'^*(L^{-1})$ and we have an exact sequence
\begin{equation}\label{ses2}
    0 \to \pi^*(\Omega_Y) \to \Omega_{Y'} \otimes \mathcal{O}_X \to \pi^*(L^{-1}) \to 0.
\end{equation} 
We also have another exact sequence as follows
\begin{equation}\label{ses3}
    0 \to \pi^*(\mathcal{O}_Y(-B)) \to \pi^*(L^{-1}) \to \pi^*(L^{-1}) \otimes \mathcal{O}_R \to 0.
\end{equation} 
Now apply snake lemma to the following diagram where the first row is ~\eqref{ses1}
\[ 
\begin{tikzcd}
     0 \arrow{r} & \pi^*(\mathcal{O}_Y(-B)) \arrow{r} \arrow{d} & \Omega_{Y'} \otimes \mathcal{O}_X \arrow{d} \arrow{r} & \Omega_X \arrow{r} \arrow{d} & 0 \\
     0 \arrow{r} & \pi^*(L^{-1}) \arrow{r}  & \pi^*(L^{-1}) \arrow{r} & 0 \arrow{r} & 0 
\end{tikzcd} 
\]
and use the previous two short exact sequences ~\eqref{ses2} and ~\eqref{ses3} to get the following exact sequence.
$$ 0 \to \pi^*(\Omega_Y) \to \Omega_X  \to \pi^*(L^{-1})\otimes \mathcal{O}_R \to 0 $$
Since $\pi$ is non-degenerate, we get that $T_{X/Y} = 0$ and hence by dualizing the above sequence we have the following exact sequence 
\begin{equation}\label{fcotex}
    0 \to T_X \to \pi^*T_Y \to \mathcal{E}\textrm{xt}^1(\pi^*(L^{-1})\otimes \mathcal{O}_R, \mathcal{O}_X)\to \mathcal{T}_X^1\to 0.
\end{equation}
Notice that $\mathcal{E}\textrm{xt}^1(\pi^*(L^{-1})\otimes \mathcal{O}_R, \mathcal{O}_X) =  \pi^*(L) \otimes \mathcal{O}_R(R) = \pi^*(\mathcal{O}_Y(B)) \otimes \mathcal{O}_R $, and consequently ~\eqref{fcotex} becomes
\begin{equation}\label{ses4}
    0 \to T_X \to \pi^*T_Y \to \pi^*(\mathcal{O}_Y(B)) \otimes \mathcal{O}_R\to \mathcal{T}_X^1\to 0.
\end{equation}

The exact sequence \eqref{eq.lt-any} follows from the fact that $\mathcal{N}_{\pi}=\textrm{ker}\left(\pi^*(\mathcal{O}_Y(B)) \otimes \mathcal{O}_R\to \mathcal{T}_X^{1}\right)$.\QEDB\par

\subsection{Geometry of the deformation spaces}\label{prelim 5} One of our objectives is to describe the moduli components of surfaces of type $(1)_m$ and $(1')_m$. We will see that for a fixed $m$, there is a unique component of the moduli space that contains all surfaces of both types, and that this component is uniruled. The proof of this fact is based on the following result. 

\begin{theorem}\label{Def_1}
Assume the hypothesis (a), (b) and (d) of Theorem ~\ref{def to 2:1}. Then ${\bf Def}_{\pi/\mathbb{P}^N}$ has a smooth uniruled algebraic formally semiuniversal deformation space $V_{\pi/\mathbb{P}^N}$.
\end{theorem}

\noindent\textit{Proof.} We construct an algebraic formally semiuniversal family of deformations of the functor ${\bf Def}_{\pi/\mathbb{P}^N}$
$$\mathcal{X}_{\pi/\mathbb{P}^N}\to \mathbb{P}^N_{V_{\pi/\mathbb{P}^N}}\to V_{\pi/\mathbb{P}^N}$$  
over a smooth pointed irreducible base $(V_{\pi/\mathbb{P}^N},0)$. \par 

Let $ \mathcal{Y}_{\psi} \to \mathbb{P}_{U_{\psi}}^N \to U_{\psi}$ be the algebraic formally semiuniversal family of deformations of the functor ${\bf Def}_{\psi}$ (this space exists, see for example the proof of Theorem ~\ref{def to 2:1}).\par

Let $(\mathscr{Y}_L,\mathscr{L}) \to U_L$ be the algebraic formally semiuniversal deformation space of the functor ${\bf Def}_{(Y,L)}$. Let $\mathscr{Y} \to U$ be the algebraic formally semiuniversal deformation space of the functor ${\bf Def}_{Y}$. Forgetful maps between functors induce a cartesian diagram, which in turn induces a Cartesian diagram of algebraic formally semiuniversal deformation spaces as shown below. 

\begin{minipage}{.35\textwidth}
\[ 
\begin{tikzcd}
{\bf Def}_{\psi} \times_{{\bf Def}_{Y}} {\bf Def}_{(Y,L)} \arrow{r} \arrow{d} & {\bf Def}_{(Y,L)} \arrow{d} \\
{\bf Def}_{\psi} \arrow{r} & {\bf Def}_{Y}
\end{tikzcd}
\]
\end{minipage}
\begin{minipage}{.3\textwidth}
$\qquad\qquad\qquad\xymatrix@C=4em{{}\ar@{~>}[r]&{}}$
\end{minipage}
\begin{minipage}{.3\textwidth}
\[ 
\begin{tikzcd}
U_{\psi} \times_{U} U_L \arrow{r} \arrow{d} & U_L \arrow{d} \\
U_{\psi} \arrow{r} & U 
\end{tikzcd}
\]
\end{minipage}

Since $H^2(\mathscr{O}_Y) = 0$, we have that the forgetful map ${\bf Def}_{(Y,L)} \to {\bf Def}_{Y}$ is smooth (see \cite{Ser}, Proposition $2.3.6$) and hence the map ${\bf Def}_{\psi} \times_{{\bf Def}_{Y}} {\bf Def}_{(Y,L)} \to {\bf Def}_{\psi} $ is smooth. Now using the fact that ${\bf Def}_{\psi}$ is smooth, we have that ${\bf Def}_{\psi} \times_{{\bf Def}_{Y}} {\bf Def}_{(Y,L)}$ is smooth and hence $U_{\psi} \times_{U} U_L$ is smooth. We set $U_{(\psi,L)} := U_{\psi} \times_{U} U_L$. The semiuniversal families  form the following cartesian diagram
\[ \begin{tikzcd}
 (\mathcal{Y}_{\psi} \times_{\mathcal{Y}} \mathcal{Y}_L \to \mathbb{P}_{U_{(\psi,L)}}^N, \mathcal{L}_{\psi}) \arrow{r} \arrow{d} & (\mathcal{Y}_L,\mathscr{L}) \arrow{d} \\
(\mathcal{Y}_{\psi} \to \mathbb{P}_{U_{\psi}}^N) \arrow{r} & \mathcal{Y} 
\end{tikzcd}
\]
Hence $(\mathcal{Y}_{\psi} \times_{\mathcal{Y}} \mathcal{Y}_L \to \mathbb{P}_{U_{(\psi,L)}}^N \to U_{(\psi,L)}, \mathcal{L}_{\psi})$ is a smooth algebraic formally semiuniversal deformation of ${\bf Def}_{\psi} \times_{{\bf Def}_{Y}} {\bf Def}_{(Y,L)}$ where $\mathcal{L}_{\psi}$ is the pullback of $\mathcal{L}$ under the morphism $\mathcal{Y}_{\psi} \times_{\mathcal{Y}} \mathcal{Y}_L \to \mathcal{Y}_L$. Let the map $\mathcal{Y}_{\psi,L}:= \mathcal{Y}_{\psi} \times_{\mathcal{Y}} \mathcal{Y}_L \to U_{(\psi,L)}$ be denoted by $p_{(\psi,L)}$.  Since $H^1(L^{\otimes 2}) = 0$, we have that $p_{(\psi,L)*}(\mathcal{L_{\psi}}^{\otimes 2})$ is free after possibly shrinking $U_{(\psi,L)}$. Let $V_{\pi/\mathbb{P}^N} := \mathbb{P}(p_{(\psi,L)*}(\mathcal{L_{\psi}}^{\otimes 2}))$ and consider the Cartesian diagram 
\[ \begin{tikzcd}
 \mathcal{Y}_{\psi,L} \times_{U_{\psi,L}} V_{\pi/\mathbb{P}^N} \arrow{r} \arrow{d} & V_{\pi/\mathbb{P}^N} \arrow{d} \\
\mathcal{Y}_{\psi,L} \arrow{r} & U_{\psi,L} 
\end{tikzcd}
\]
Choose a basis $\bigoplus\limits_{i=0}^M \mathcal{O}_{U_{(\psi,L)}}s_i$ of $p_{(\psi,L)*}(\mathcal{L_{\psi}}^{\otimes 2})$. Let $X_i \in H^0(\mathcal{O}_{V_{\pi/\mathbb{P}^N}}(1)) =   H^0(p_{(\psi,L)*}(\mathcal{L_{\psi}}^{\otimes 2})^*)$, with $0 \leq i \leq M$ be the dual basis. Now on $\mathcal{Y}_{\psi,L} \times_{U_{(\psi,L)}} V_{\pi/\mathbb{P}^N}$, consider the divisor $\mathcal{B} = \sum\limits_{i=0}^M X_is_i$. One can construct a relative Galois double cover $\mathcal{X}_{\pi/\mathbb{P}^N} \to \mathcal{Y}_{\psi,L} \times_{U_{(\psi,L)}} V_{\pi/\mathbb{P}^N}$ given by the equation $t^2-\mathcal{B}$ in the total space of $q^*(\mathcal{L}_{\psi})$ where $t$ is the tautological section of $q^*(\mathcal{L}_{\psi})$ and $q: \mathcal{Y}_{\psi,L} \times_{U_{(\psi,L)}} V_{\pi/\mathbb{P}^N} \to \mathcal{Y}_{\psi,L} $. The fibre of this relative double cover over a point $(u,[r]) \in V_{\pi/\mathbb{P}^N}$ with $u \in U_{(\psi,L)}$ and $r \in H^0(\mathcal{L}_{{\psi},u}^{\otimes 2})$ is the double cover $\mathcal{X}_{\pi/\mathbb{P}^N,u} \to \mathcal{Y}_{\psi,L,u}$ given by the line bundle $L_{\psi,u}$ and  $r \in H^0(\mathcal{L}_{{\psi},u}^{\otimes 2})$. This is therefore a smooth algebraic deformation of the functor ${\bf Def}_{\pi/\mathbb{P}^N}$.\par 

Now note that given a flat family of polarized schemes $f: (\mathcal{C}, \mathcal{M}) \to S$ over an affine scheme $S$ with $f_*(\mathcal{M})$ free, giving a divisor $\mathcal{D}$ defined as the zero locus of a section in $H^0(\mathcal{M})$ flat over $S$ is equivalent to giving a unique $S$-valued point in $\mathbb{P}(f_*(\mathcal{M}))$ and hence a section $S \to \mathbb{P}(f_*(\mathcal{M}))$. This along with the fact that $U_{(\psi,L)}$ is formally semiuniversal implies that $V_{\pi/\mathbb{P}^N}$ is a smooth algebraic formally semiuniversal deformation of the functor ${\bf Def}_{\pi/\mathbb{P}^N}$. Also since it is a projective bundle over a smooth affine scheme, it is uniruled. \QEDB 

\begin{remark}
Under the assumptions of Theorem ~\ref{Def_1}, if $X$ is smooth, it is easy to prove the smoothness of ${\bf Def}_{\pi/\mathbb{P}^N}$ only using the existence of a formally semiuniversal deformation, without the explicit construction. Indeed, since $H^1(\mathcal{N}_{\pi})=0$ by \cite{Pa}, Corollary 4.1 or \cite{GGP2}, (2.11), and the assumptions $H^1(L^{\otimes 2})=H^2(\mathcal{O}_Y)=0$, it follows from the following exact sequence (see the top row of the  commutative diagram in \cite[Proof of Prop. 1.10]{Weh}):
$$\cdots \to H^0(\mathcal{N}_{\pi}) \to T^1(\pi/\mathbb{P}^N) \to H^0(\mathcal{N}_{\psi})\to H^1(\mathcal{N}_{\pi})\to T^2(\pi/\mathbb{P}^N) \to H^1(\mathcal{N}_{\psi}) \to \cdots$$
Since $H^1(\mathcal{N}_{\pi})=0$, and \cite{Ser}, Proposition 2.3.6 that the forgetful map ${\bf Def}_{\pi/\mathbb{P}^N}\to {\bf Def}_{\psi}$ is smooth. Consequently ${\bf Def}_{\pi/\mathbb{P}^N}$ is smooth as ${\bf Def}_{\psi}$ is smooth by hypothesis.
\end{remark}

The following corollary shows that if ${\bf Def}_{\varphi}$ has an algebraic formally universal deformation space, then that space is also smooth and uniruled under suitable assumptions. In fact, one can also expect to determine the degree of a general deformation of $\varphi$.

\begin{corollary}\label{Defphi}
Assume the hypotheses (a), (b), (c), and (d) of Theorem ~\ref{def to 2:1}. Assume 
\begin{itemize}
    \item[(a)] {\bf Def}$_{\varphi}$ has an algebraic formally universal deformation space,
    \item[(b)] $H^0(N_{\psi}\otimes \mathscr{E}_{\pi}) = 0$.
\end{itemize}
 Then the the following happens:
 \begin{itemize}
     \item[(1)] the natural forgetful map  ${\bf Def}_{\pi/\mathbb{P}^N} \to {\bf Def}_{\varphi}$ is smooth,
     \item[(2)] the algebraic formally universal deformation space of ${\bf Def}_{\varphi}$ is smooth and  uniruled, and 
     \item[(3)] a general deformation of ${\varphi}$ is a composition of a double cover over a  deformation $Y'$ of $Y$ followed by a morphism of $Y'\to \mathbb{P}^N$ that is birational onto its image.
 \end{itemize} 
\end{corollary}

\noindent\textit{Proof.} The smoothness of ${\bf Def}_{\pi/\mathbb{P}^N} \to {\bf Def}_{\varphi}$ is a consequence of Proposition ~\ref{corweh}, thanks to assumption (b) and Theorem \ref{Def_1}. Moreover, ${\bf Def}_{\varphi}$ is smooth since ${\bf Def}_{\pi/\mathbb{P}^N}$ is smooth (see \cite{Ser}, Proposition 2.2.5 (iii)), thanks to Theorem ~\ref{Def_1}.\par 

Now we show that the algebraic formally semiuniversal deformations space of ${\bf Def}_{\varphi}$, which we denote by $U_{\varphi}$, is uniruled. In the notation of Theorem ~\ref{Def_1}, after possibly shrinking $V_{\pi/\mathbb{P}^N}$ we can assume that $V_{\pi/\mathbb{P}^N} = U_{(\psi,L)} \times \mathbb{P}^m$ where over a point $u \in U_{(\psi,L)}$, the fibre which is a projective space that parametrizes the divisors in the linear system of $\mathcal{L}_{\psi,u}$ which are branch divisors of the finite morphism $X_u \to Y_u$. The conclusion follows since a branch divisor is uniquely determined by the finite morphism.\par
Finally, part (3) follows from Corollary ~\ref{maincor}. \QEDB\par

\vspace{5pt}

Now we provide the consequences of the above results on the deformations of $X$.

\begin{corollary}\label{defcanonicalmorphism }
Assume all the hypotheses of Corollary ~\ref{Defphi}. Furthermore assume ${\bf Def}_X$ has an algebraic formally semiuniversal deformation space. If the natural forgetful map ${\bf Def}_{\varphi} \to {\bf Def}_{X}$ has surjective differential map then ${\bf Def}_{X}$ is smooth and the algebraic formally semiuniversal deformation space of $X$ is uniruled.
\end{corollary}
\noindent\textit{Proof.} Since ${\bf Def}_{\varphi}$ is smooth by Corollary ~\ref{Defphi}, the smoothness of ${\bf Def}_X$ follows from \cite{Ser}, Proposition $2.3.7$. Composing by the smooth surjection $V_{\pi/\mathbb{P}^N} \to U_{\varphi}$, we have a smooth surjection $V_{\pi/\mathbb{P}^N} \to U$ where $U$ is the algebraic formally universal deformation space of ${\bf Def}_X$. Lastly, $U$ is uniruled since $X$ is normal and for a normal abelian cover, the branch divisors are uniquely determined by $X$.\QEDB\par 

\section{Deformations of irregular covers of surface scrolls} \label{3}

The objective of this section is to study the deformations of the canonical morphisms of  surfaces of each of the four 
types $(1)_m$ , $(1')_m$, $(2)_m$ and $(3)_m$ described in Theorem ~\ref{gpmain}. In particular, we aim to prove the following 

\begin{theorem}\label{A}
Let $X$ be an irregular surface with at worst canonical singularities. Assume the canonical bundle $\omega_X$ is ample and globally generated, and the canonical morphism $\varphi$ is a quadruple Galois cover onto a smooth surface of minimal degree, i.e, $X$ belongs to one of the four families described in Theorem ~\ref{gpmain}. Then  we have the following description of the algebraic formally semiuniversal deformation space \textcolor{black}{$U_\varphi$} of $\varphi$ (that exists by Remark ~\ref{can}).
\begin{itemize}
    \item[(1)] \textcolor{black}{Any deformation of $\varphi$ factors through a double cover of a ruled surface over a smooth curve of genus
    \begin{itemize}
        \item[(I)] $g = 1$ if $X$ is of type $(1)_m$ or $(1')_m$;
        \item[(II)] $g = m$ if $X$ is of type $(2)_m$;
        \item[(III)] $g = m+1$ if $X$ is of type $(3)_m$.
    \end{itemize}
    In particular, there do not exist any irreducible component in the algebraic formally semiuniversal deformation space of $\varphi$, such that its general element is birational onto its image.}

    \item[(2)] If $X$ belongs to the family of type $(1)_m$ $(m\geq 1)$, $(1')_m$ $(m\geq 1)$, $(2)_2$ or $(3)_m$ $(m\geq 2)$, then
    \textcolor{black}{$U_\varphi$ is irreducible} \textcolor{black}{and a} general element $\varphi'$ of $U_\varphi$ is a two-to-one morphism onto its image,  
    \begin{itemize}
        \item[(I)] which is a non-normal variety whose normalization is an elliptic ruled surface with invariant $e = 0$, if $X$ is a surface of type $(1)_m$ or $(1')_m$;
        \item[(II)] which is a smooth surface ruled over a smooth curve of genus $2$ with invariant $e = -2$, if $X$ is a surface of type $(2)_2$ ;
        \item[(III)] which is a product of a smooth curve of genus $2$ with a smooth non-hyperelliptic curve of genus $m+1$ if $X$ is a surface of type $(3)_m$ $(m\geq 2)$. 
    \end{itemize} 
 Moreover, $\varphi'$ is induced by  the complete linear series of a line bundle numerically equivalent to \textcolor{black}{the canonical} \textcolor{black}{(in case (III), any element of $U_\varphi$ is a canonical morphism)}.
    \item[(3)] If $X$ is of type $(3)_1$, any deformation of $\varphi$ is a \textcolor{black}{canonical} morphism of degree four onto its image which is $\mathbb{F}_0$. 
\end{itemize}
\end{theorem}

\textcolor{black}{We will give the proof of Theorem~\ref{A} at the end of this section, as a consequence of the results we will be proving in it.} First we fix the notations that we are going to use throughout this section. It follows from Theorem ~\ref{gpmain} that if $\pi:X\to Y$ is an irregular quadruple Galois canonical cover of a smooth surface of minimal degree $Y$ with trace zero module $\mathcal{E}$, then  $i:Y=\mathbb{P}^1\times\mathbb{P}^1\hookrightarrow\mathbb{P}^N$ and the embedding is given by the complete linear series $|\mathcal{O}_Y(C_0+mf)|$. We have $N=2m+1$ and we identify $\mathcal{O}_Y(C_0)$ with $\mathcal{O}_Y(1,0)$ and $\mathcal{O}_Y(f)$ with $\mathcal{O}_Y(0,1)$. Notice that $T_Y=\mathcal{O}_Y(2,0)\oplus\mathcal{O}_Y(0,2)$, whence $h^0(T_Y)=6$. One has the following two exact sequences;
\begin{equation}\label{ttnany}
    0\to \mathcal{O}_Y(2,0)\oplus\mathcal{O}_Y(0,2)\to T_{\mathbb{P}^N|Y}\to\mathcal{N}_{Y/\mathbb{P}^N}\to 0,
\end{equation}
\begin{equation}\label{eulerany}
    0\to\mathcal{O}_Y\to\mathcal{O}_Y(1,m)^{\oplus N+1}\to T_{\mathbb{P}^N|Y}\to 0.
\end{equation}

\begin{lemma}\label{cohomology Y}
Let $Y=\mathbb{P}^1\times\mathbb{P}^1\hookrightarrow\mathbb{P}^N$ be the embedding given by the complete linear series $|\mathcal{O}_Y(C_0+mf)|$. 
\begin{itemize}
    \item[(1)] $H^0(T_{\mathbb{P}^N|Y})=(N+1)^2-1$, $H^0(\mathcal{N}_{Y/\mathbb{P}^N})=(N+1)^2-7$.
    \item[(2)] $H^1(T_{\mathbb{P}^N|Y})=0$, $H^1(\mathcal{N}_{Y/\mathbb{P}^N})=0$.
\end{itemize}
\end{lemma}
\noindent\textit{Proof.} Since $Y$ is regular with $H^2(\mathcal{O}_Y)=0$, the assertions about $H^j(T_{\mathbb{P}^N|Y})$ for $j=0,1$ follows from ~\eqref{eulerany}. Consequently, it is easy to compute $H^j(\mathcal{N}_{Y/\mathbb{P}^N})$ for $j=0,1$ using ~\eqref{ttnany}.\QEDB\par  

\vspace{5pt}

Now we fix our notations for nonrational ruled surfaces. A nonrational ruled surface over a nonrational smooth curve $C$ of genus $g\neq 0$ is by definition a projective bundle $W=\mathbb{P}(\mathcal{E'})$ where $\mathcal{E'}$ is a rank 2 vector bundle on $C$. We will always assume that $\mathcal{E'}$ is normalized, i.e., $\mathcal{E'}$ has sections, but any twist of $\mathcal{E'}$ by any line bundle of negative degree has no section. By definition $e:=-\textrm{deg}(\textrm{det}(\mathcal{E'}))$ is the invariant of $W$. A section of $p':W\to C$ determines a sectional curve $C_0'$ with self intersection $-e$, and let $f'$ be the numerical class of a fiber of $p'$. It is known that $\textrm{Pic}(W)=\mathbb{Z}C_0'\oplus p'^*\textrm{Pic}(C)$. In particular, the N\'eron-Severi group $\textrm{NS}(W)=\mathbb{Z}C_0'\oplus\mathbb{Z}f'$ satisfying
$$C_0'^2=-e,\quad C_0'f'=1,\quad f'^2=0.$$
If ${\bf a}$ is a divisor on $C$, ${\bf a}f'$ denotes the pull-back $p'^*{\bf a}$. The canonical bundle $\omega_W=\mathcal O_W(-2C_0'+({\bf e}+\omega_C)f')$, where ${\bf e}:=\textrm{det}(\mathcal{E'})$, consequently $$\omega_W\equiv -2C_0'-(e+2-2g)f'.$$

\subsection{Deformations of canonical morphisms for types \texorpdfstring{$(1)_m$}{TEXT} and \texorpdfstring{$(1')_m$}{TEXT}}\label{study11'} For these surfaces, we have the following diagram. 
\[\begin{tikzcd}
X\arrow[r, "\pi_1"]\ar[rr,out=-30,in=200,swap,"\pi"] & X_1\ar[r,"p_1"] & Y\overset{{i}}\hookrightarrow\mathbb{P}^N
\end{tikzcd}
\]
We also know that ${p_1}_*\mathcal{O}_{X_1}=\mathcal{O}_Y\oplus\mathcal{O}_Y(-2C_0)$. Since we have identified $\mathcal{O}_Y(C_0)$ with $\mathcal{O}_Y(1,0)$, we can write ${p_1}_*\mathcal{O}_{X_1}=\mathcal{O}_Y\oplus\mathcal{O}_Y(-2,0)$. It is easy to see that  $X_1=E\times\mathbb{P}^1$ where $\psi:E\to\mathbb{P}^1$ is a smooth double cover, with $\psi_*\mathcal{O}_E=\mathcal{O}_{\mathbb{P}^1}\oplus\mathcal{O}_{\mathbb{P}^1}(-2)$, i.e. $E$ is a smooth elliptic curve. We set $\varphi_1:=i\circ p_1$ and call $B$ the branch divisor of $\pi_1$.

\begin{proposition}\label{11'}
Let $X$ be a surface of type $(1)_m$ or $(1')_m$. Then the following happens:
\begin{itemize}
    \item[(1)] $h^1(\mathcal{O}_{X_1})=1$ and $h^2(\mathcal{O}_{X_1})=0$,
    \item[(2)] $h^0(\mathcal{N}_{p_1})=4$ and $h^1(\mathcal{N}_{p_1})=0$,
    \item[(3)] $h^0(T_{\mathbb{P}^N|Y}\otimes\mathcal{O}_Y(-2,0))=1$ and $h^1(T_{\mathbb{P}^N|Y}\otimes\mathcal{O}_Y(-2,0))=0$,
    \item[(4)]  $h^0(\mathcal{N}_{Y/\mathbb{P}^N}\otimes\mathcal{O}_Y(-2C_0))=3$ and $h^1(\mathcal{N}_{Y/\mathbb{P}^N}\otimes\mathcal{O}_Y(-2C_0))=0$,
    \item[(5)] $h^0(\mathcal{N}_{\varphi_1})=(N+1)^2$ and $h^1(\mathcal{N}_{\varphi_1})=0$; consequently $\varphi_1$ is unobstructed. 
\end{itemize}
\end{proposition}
\noindent\textit{Proof.} (1) Follows from $h^j(\mathcal{O}_{X_1})=h^j(\mathcal{O}_Y)\oplus h^j(\mathcal{O}_Y(-2,0)$ and K\"unneth formula.\par 
(2) We apply \cite{Pa}, Corollary 4.1 or \cite{GGP2}, (2.11). Since $Y$ is regular, it follows that $h^0(\mathcal{N}_{p_1})=h^0(\mathcal{O}_Y(4,0))-1=4$. Furthermore, since $H^2(\mathcal{O}_Y)=0$, we obtain $h^1(\mathcal{N}_{p_1})=h^1(\mathcal{O}_Y(4,0))=0$.\par 
(3) Tensor ~\eqref{eulerany} by $\mathcal{O}_Y(-2,0)$ to obtain the following exact sequence $$0\to \mathcal{O}_Y(-2,0)\to \mathcal{O}_Y(-1,m)^{\oplus N+1}\to T_{\mathbb{P}^N|Y}\otimes\mathcal{O}_Y(-2,0)\to 0.$$ It follows that $h^0(T_{\mathbb{P}^N|Y}\otimes\mathcal{O}_Y(-2,0))=1$, $h^1(T_{\mathbb{P}^N|Y}\otimes\mathcal{O}_Y(-2,0))=0$.\par 
(4) This is a consequence of the long exact sequence associated to the 
{exact sequence}~\eqref{ttnany} tensored by $\mathcal{O}_Y(-2,0)$ and part (3).\par 
(5) We have the following short exact sequence by Lemma ~\ref{exactgonzalez}
\begin{equation}\label{exactgonzalezp1}
    0\to\mathcal{N}_{p_1}\to \mathcal{N}_{\varphi_1}\to p_1^*\mathcal{N}_{Y/\mathbb{P}^N}\to 0.
\end{equation}
It follows that $h^0(\mathcal{N}_{\varphi_1})=h^0(\mathcal{N}_{p_1})+h^0(p_1^*\mathcal{N}_{Y/\mathbb{P}^N})=h^0(\mathcal{N}_{p_1})+h^0(\mathcal{N}_{Y/\mathbb{P}^N})+h^0(\mathcal{N}_{Y/\mathbb{P}^N}\otimes\mathcal{O}_Y(-2C_0))$, since $h^1(\mathcal{N}_{p_1})=0$ by part (2). We obtain $h^0(\mathcal{N}_{\varphi_1})=(N+1)^2$ thanks to part (2) and Lemma ~\ref{cohomology Y}.\par 
The fact $h^1(\mathcal{N}_{\varphi_1})=0$ follows from the vanishings of $h^1(\mathcal{N}_{p_1})$ (proven in part (2)), $h^1(\mathcal{N}_{Y/\mathbb{P}^N})$ (proven in Lemma ~\ref{cohomology Y}), and $h^1(\mathcal{N}_{Y/\mathbb{P}^N}\otimes\mathcal{O}_Y(-2C_0))$ (proven in part (4)).\QEDB\par 

\begin{corollary}\label{11't}
Let $X$ be a  surface of type $(1)_m$ or $(1')_m$. Then 
there exists a smooth, affine irreducible algebraic curve $T$ and a flat family of morphisms $\Phi: \mathcal{X}\to\mathbb{P}^N_T$ over $T$ for which the following happens;
\begin{itemize}
    \item[(a)] $\Phi_t:\mathcal{X}_t \to \mathbb{P}^{2m+1}$ is a morphism of degree two from a normal projective surface with at worst canonical singularities for all $t \in T-\{0\}$. Further for any $t \in T-\{0\}$, the normalization of Im($\varphi_t$) is an elliptic ruled surface which is the projectivization of a rank two split vector bundle on an elliptic curve and has invariant $e = 0$. Further one can take $\Phi_t$ to be the canonical morphism of $\mathscr{X}_t$.
    \item[(b)] $\Phi_0:\mathcal{X}_0 \to \mathbb{P}^{2m+1}$ is the canonical morphism $\varphi:X \to \mathbb{P}^N$.
\end{itemize}
Moreover the forgetful map from ${\bf Def_{\pi_1/\mathbb{P}^N}} \to {\bf Def_{\varphi}}$ is smooth and hence any deformation of $\varphi$ \textcolor{black}{factors through a double cover of an elliptic ruled surface and is hence} a morphism of degree $\geq 2$. Hence in particular $\varphi$ cannot be deformed to a birational morphism.
\end{corollary}
\noindent\textit{Proof.} We check the hypotheses of Theorem ~\ref{def to 2:1}. Hypothesis (a) has been checked in Proposition ~\ref{11'} (1). To check hypothesis (b), we need to prove that ${\bf Def}_{\varphi_1}$ is smooth, which we have showed in Proposition ~\ref{11'} (5). To check hypothesis (c), we need to check $H^0(\mathcal{N}_{\varphi_1})\to H^0(\mathcal{N}_{Y/\mathbb{P}^N}\otimes\mathcal{O}_Y(-2,0))$ is non-zero. This is a consequence of the long exact sequence associated to ~\eqref{exactgonzalezp1}, and the facts that $h^1(\mathcal{N}_{p_1})=0$ (proven in Proposition ~\ref{11'} (2)), and $h^0(\mathcal{N}_{Y/\mathbb{P}^N}\otimes\mathcal{O}_Y(-2,0))\neq 0$ (proven in Proposition ~\ref{11'} (4)). The fact that $\mathcal{X}_t$ is a normal projective surface with at worst canonical singularities follows thanks to \cite{Kaw}. 
Now note that since $h^1(T_{\mathbb{P}^N|Y}\otimes\mathcal{O}_Y(-2,0))= h^1(T_{\mathbb{P}^N|Y})=0$ (by Proposition ~\ref{11'} (3) and Lemma ~\ref{cohomology Y} (2)), we have that $h^1(\varphi_1^*(T_{\mathbb{P}^N}))=0$ and ${\bf Def_{\varphi_1}} \to {\bf Def}_{X_1}$ is smooth. Hence the map $H^0(\mathcal{N}_{\varphi_1}) \to H^1(T_{X_1})$ is surjective. By \cite{Sei92}, Lemma $12$, there exist an open set in $H^1(T_{X_1})$ such that for a smooth curve along a first order deformation belonging to the open set a general deformation of $X_1$ along the curve is an elliptic ruled surface which is the projectivization of a split rank two vector bundle with invariant $e = 0$. Also $H^0(\mathcal{N}_{\varphi_1})\to H^0(\mathcal{N}_{Y/\mathbb{P}^N}\otimes\mathcal{O}_Y(-2,0))$ is surjective and there exist an open set of non-zero elements in $H^0(\mathcal{N}_{Y/\mathbb{P}^N}\otimes\mathcal{O}_Y(-2,0))$. Hence one can choose an element (in fact an open set of elements) from $H^0(\mathcal{N}_{\varphi_1})$ such that it maps to a non-zero element in      $H^0(\mathcal{N}_{Y/\mathbb{P}^N}\otimes\mathcal{O}_Y(-2,0))$ and the general induced deformation $\mathcal{X}_{1t}$ of $X_1$ is an elliptic ruled surface which is the projectivization of a split rank two vector bundle with invariant $e = 0$.  Finally to check hypothesis ($d$), note that $H^1(L^{\otimes 2}) = H^1(B) = 0$ by Proposition ~\ref{moduli11'}, ($1$).\par
Note that $\mathcal{E}_{\pi_1} = p_1^*(\mathcal{O}_Y(-C_0-(m+2)f))$. Let $\mathcal{E}_{p_1} = \mathcal{O}_Y(-2C_0)$. Then $$\omega_{X_1} = p_1^*(\omega_{Y} \otimes \mathcal{O}_Y(2C_0)) = p_1^*(\mathcal{O}_Y(-2C_0-2f) \otimes \mathcal{O}_Y(2C_0)) = p_1^*(-2f).$$ Thus, we obtain 
\begin{equation}\label{equation.11't}
\omega_{X_1} \otimes  \varphi_1^*(\mathcal{O}_{\mathbb{P}^N}(-1)) = p_1^*(-2f) \otimes p_1^*(-C_0-mf) = p_1^*(\mathcal{O}_Y(-C_0-(m+2)f)) = \mathcal{E}_{\pi_1}.    
\end{equation}
 Hence by Theorem ~\ref{def to 2:1}, we can take $\Phi_t$ to be the canonical morphism of $X_t$.   \par
The second assertion follows from Corollary ~\ref{maincor}. The existence of an algebraic formally semiuniversal deformation space of $\varphi$ follows from Remark ~\ref{can}.  
To finish the proof,  we need to show that $H^0(\mathcal{N}_{\varphi_1}\otimes\mathcal{E}_{\pi_1})=0$ where $\mathcal{E}_{\pi_1}$ is the trace zero module of $\pi_1$.
We make use of the fact that $X_1=E\times \mathbb{P}^1$ for an elliptic curve $E$. Recall that $\psi:E\to\mathbb{P}^1$ is the morphism induced by the restriction of $p_1$, satisfies $\psi_*\mathcal{O}_E=\mathcal{O}_{\mathbb{P}^1}\oplus \mathcal{O}_{\mathbb{P}^1}(-2)$. We have $${\pi_1}_*\mathcal{O}_X=\mathcal{O}_{X_1}\oplus(\psi^*\mathcal{O}_{\mathbb{P}^1}(-1)\boxtimes\mathcal{O}_{\mathbb{P}^1}(-m-2)).$$
Also recall that $T_{X_1}=(\mathcal{O}_E\boxtimes\mathcal{O}_{\mathbb{P}^1}(2))\oplus\mathcal{O}_{X_1}$ and $\mathcal{E}_{\pi_1}=\psi^*\mathcal{O}_{\mathbb{P}^1}(-1)\boxtimes\mathcal{O}_{\mathbb{P}^1}(-m-2)$. It is easy to check that $H^1(T_{X_1}\otimes\mathcal{E}_{\pi_1})=0$. One has the following pullback of the Euler sequence;
\begin{equation}\label{inteu}
    0\to\mathcal{O}_{X_1}\to p_1^*\mathcal{O}_Y(1,m)^{\oplus N+1}\to\varphi_1^*T_{\mathbb{P}^N}\to 0.
\end{equation}
We tensor ~\eqref{inteu} by $\mathcal{E}_{\pi_1}$ and take the long exact sequence of cohomology. Notice that $H^1(\mathcal{E}_{\pi_1})=0$, and $H^0(p_1^*\mathcal{O}_Y(1,m)\otimes\mathcal{E}_{\pi_1})=0$, consequently $H^0(\varphi_1^*T_{\mathbb{P}^N}\otimes\mathcal{E}_{\pi_1})=0$. Now consider the exact sequence;
\begin{equation}\label{ttn11'}
    0\to T_{X_1}\to\varphi_1^*T_{\mathbb{P}^N}\to\mathcal{N}_{\varphi_1}\to 0.
\end{equation}
It follows from the long exact sequence of cohomology that $H^0(\mathcal{N}_{\varphi_1}\otimes\mathcal{E}_{\pi_1})=0$. \QEDB\par 

\vspace{5pt}

Before moving on to the next case, we make a remark that will help us to see that for these surfaces $H^1(\mathcal{N}_{\varphi})\neq 0$.

\begin{remark}\label{nonzero!!}
It follows from the vanishing of $H^1(T_{X_1}\otimes\mathcal{E}_{\pi_1})$ and the long exact sequence associated to ~\eqref{ttn11'} that $h^1(\mathcal{N}_{\varphi_1}\otimes\mathcal{E}_{\pi_1})\geq h^1(\varphi_1^*T_{\mathbb{P}^N}\otimes\mathcal{E}_{\pi_1})\geq  h^1(T_{\mathbb{P}^N|Y}\otimes\mathcal{O}_Y(-1,-m-2))=N+1$, where the last equality follows from ~\eqref{eulerany}.
\end{remark}

\color{black}
\begin{corollary}\label{existence of morphism deformation component $(1)$, $(1')$}
Let $X$ be a surface of type $(1)_m$ or $(1')_m$. 
Then there exists an irreducible component $U_{\varphi}$ of the algebraic formally semiuniversal deformation space of ${\varphi}$ (that exists by Remark ~\ref{can}) whose general elements are  two-to-one morphisms onto their image whose normalization is an elliptic ruled surface with invariant $e = 0$. Further, there does not exist any component of the algebraic formally semiuniversal deformation space of ${\varphi}$ whose general elements are morphisms that are birational onto their image.
\end{corollary}

\noindent\textit{Proof.} Since the curve constructed in Corollary ~\ref{11't} is irreducible, it is contained in an irreducible component. Now the assertion follows by applying semicontinuity to the reduced induced structure of the irreducible component (note that a general closed point of an irreducible scheme is the same as a general closed point of its reduced induced structure).\QEDB

\vspace{5pt}

The corollary of the following proposition shows that the image of a general morphism in the irreducible component $U_{\varphi}$ constructed above is necessarily non-normal. We remark that what we prove in the following proposition is a slightly stronger statement than what we need in order to prove Corollary ~\ref{non-normality for $(1)$, $(1')$}; to prove Corollary ~\ref{non-normality for $(1)$, $(1')$}, we only need the conclusion of the following proposition for $e=0$. 

\begin{proposition}\label{non-existence of canonical double covers on elliptic ruled}
There does not exist a surface of general type $X'$ with at worst canonical singularities and $K_{X'}^2=4p_g(X')-8$ that satisfies both of the following properties.
\begin{itemize}
    \item[(1)] There exist an ample and base point free line bundle $K\equiv\omega_{X'}$  with $h^0(K) = p_g(X')$.
    \item[(2)] The morphism $\varphi'$ induced by the complete linear series $|K|$ is two-to-one onto its image which is a smooth elliptic ruled surface with invariant $e \geq 0$.
\end{itemize}
\end{proposition}

\noindent\textit{Proof.} Suppose there exists such a surface $X'$ with a numerically canonical bundle $K$ satisfying the properties in the proposition. Let the image of the morphism $\varphi$ given by $|K|$ be $W$ so that the morphism $\varphi$ factors as
\begin{equation*}
    X' \xrightarrow{\pi'} W \hookrightarrow \mathbb{P}^{N'}
\end{equation*}
where $N'+1 = p_g(X')$. Let the very ample line bundle on $W$ be denoted by $aC_0'+bf'$. Note that we have $\varphi'^*(\mathcal{O}_{\mathbb{P}^{N'}}(1)) = K$. The morphism is induced by the complete linear series and hence $$h^0(K) = h^0(aC_0'+bf') + h^0((aC_0'+bf') \otimes \mathcal{E}_{\pi'}) = N'+1,$$ where $\mathcal{E}_{\pi'}$ is the trace zero module of $\pi'$. But now $\mathcal{E}_{\pi'} \equiv K_W \otimes (-aC_0'-bf')$. Hence $h^0((aC_0'+bf') \otimes \mathcal{E}_{\pi'}) = 0$ which gives $h^0(aC_0'+bf') = p_g(X')$. Since $h^i(aC_0'+bf')=0$ for $i=1,2$ (see for example \cite{GP96}, Proposition 3.1), we obtain by Riemann-Roch  
\begin{equation*}
    \frac{1}{2}(-a^2e-ae+2ab+2b) = p_g(X').
\end{equation*}
Now note that $K = \varphi^*(\mathcal{O}_{\mathbb{P}^N}(1)) = \pi^*(aC_0'+bf')$. Hence $K^2 = 2(aC_0'+bf')^2 = 2(-a^2e+2ab)$. Using the relation $K^2 = 4p_g(X')-8$, we obtain 
\begin{equation*}
    2(-a^2e+2ab) = 2(-a^2e-ae+2ab+2b) - 8.
\end{equation*}
This gives $2b-4 = ae$. But very ampleness of $aC_0'+bf'$ implies $a\geq 1$, $b \geq ae+3$ which implies $-ae \geq 2$ which is a contradiction if $e \geq 0$.\QEDB

\begin{corollary}\label{non-normality for $(1)$, $(1')$}
Consider the irreducible component $U_{\varphi}$ obtained in Corollary ~\ref{existence of morphism deformation component $(1)$, $(1')$}. There exist an open set $U_{\varphi}^0 \subseteq U_{\varphi}$ such that for a closed point $t \in U_{\varphi}^0$, Im($\varphi_t$) is non-normal, whose normalization is an elliptic ruled surface with $e=0$ which is the projectivization of a rank two split vector bundle on the elliptic curve.
\end{corollary}

\noindent\textit{Proof.} Since we are concerned with closed points $t \in U_{\varphi}$, we can take the reduced induced structure of $U_{\varphi}$ and consider the pullback of the formally semiuniversal family. Thus, without loss of generality, one can assume that $U_{\varphi}$ is integral. Let $\mathcal{X} \xrightarrow{\Phi} \mathbb{P}_{U_{\varphi}}^N \to U_{\varphi}$ be the algebraic formally semiuniversal family of $\varphi$. Since the forgetful map ${\bf Def_{\pi_1/\mathbb{P}^N}} \to {\bf Def_{\varphi}}$ is smooth we have that the above deformation factors as $\mathcal{X} \xrightarrow{\Pi_1} \mathcal{X}_1 \xrightarrow{\Phi_1} \mathbb{P}_{U_{\varphi}}^N \to U_{\varphi}$ ($\Phi = \Pi_1 \circ \Phi_1$). Let $\mathcal{Y} = \textrm{Im}(\Phi) = \textrm{Im}(\Phi_1)$. Since $\mathcal{X}$ is integral, $\mathcal{Y}$ is integral. Since $U_{\varphi}$ is integral, we have by generic flatness that $\mathcal{Y} \to U_{\varphi}$ is flat (after possibly shrinking $U_{\varphi}$). Consider the induced deformation $\mathcal{X}_1 \to U_{\varphi}$. By our choice of $U_{\varphi}$, we have that there exist $t \in U_{\varphi}$ such that $\mathcal{X}_{1t}$ is an elliptic ruled surface which is the projectivization of a rank two split vector bundle on the elliptic curve with invariant $e = 0$. Then by \cite{Sei92}, Lemma $12$, we have that for a general $t \in U_{\varphi}$, $\mathcal{X}_{1t}$ has the same property. Now for a general $t \in U_{\varphi}$, Im$(\Phi_t)$ = Im$(\Phi_{1t}) = \mathcal{Y}_t$. Also for a general $t \in U_{\varphi}$, we have that $\mathcal{X}_{1t} \to \mathcal{Y}_t$ is the normalization map. Assume for a general $t$, that $\mathcal{Y}_t$ is smooth. Then $\mathcal{Y}_t \cong \mathcal{X}_{1t} = \textrm{Im}(\varphi_t)$. But this is a contradiction to Proposition ~\ref{non-existence of canonical double covers on elliptic ruled}.
\QEDB

\color{black}
\subsection{Deformation of canonical morphism for type \texorpdfstring{$(2)_m$}{TEXT} (\texorpdfstring{$m\geq 2$}{TEXT})}\label{study2} Now we do analogous calculations for surfaces type $(2)_m$. In this case, We have the following diagram where ${p_1}_*\mathcal{O}_{X_1}=\mathcal{O}_Y\oplus\mathcal{O}_Y(-3,-1)$ and ${p_2}_*\mathcal{O}_{X_2}=\mathcal{O}_Y\oplus\mathcal{O}_Y(0,-m-1)$.
\begin{equation}\label{pic0}
\begin{tikzcd}
    X \arrow{r}{\pi_1} \arrow{d}[swap]{\pi_2} & X_1 \arrow{d}{p_1} \\
    X_2 \arrow{r}[swap]{p_2} & Y \arrow[hookrightarrow]{r}[swap]{i} & \mathbb{P}^N 
\end{tikzcd} 
\end{equation}
Notice $X_2=\mathbb{P}^1\times C$ for a smooth curve $C$ that is a double cover $\psi:C\to\mathbb{P}^1$ with $\psi_*\mathcal{O}_C=\mathcal{O}_{\mathbb{P}^1}\oplus \mathcal{O}_{\mathbb{P}^1}(-m-1)$. Set $\varphi_j=p_j\circ i$ for $j=1,2$ and call $B$ the 
branch divisor of $\pi_2$. 

\begin{proposition}\label{2m}
Let $X$ be a surface of type $(2)_m$ ($m\geq 2$). Then the following happens;
\begin{itemize}
    \item[(1)] $h^1(\mathcal{O}_{X_2})=m$ and $h^2(\mathcal{O}_{X_2})=0$,
    \item[(2)] $h^0(\mathcal{N}_{p_2})=2m+2$ and $h^1(\mathcal{N}_{p_2})=0$,
    \item[(3)] $h^0(T_{\mathbb{P}^N|Y}\otimes\mathcal{O}_Y(0,-m-1))=h^1(\mathcal{O}_Y(0,-m-1))=m$ and $h^1(T_{\mathbb{P}^N|Y}\otimes\mathcal{O}_Y(0,-m-1))=0$,
    \item[(4)] $h^0(\mathcal{N}_{Y/\mathbb{P}^N}\otimes\mathcal{O}_Y(0,-m-1))=5m-2$ and $h^1(\mathcal{N}_{Y/\mathbb{P}^N}\otimes\mathcal{O}_Y(0,-m-1))=0$,
    \item[(5)] $h^0(\mathcal{N}_{\varphi_2})=(N+1)^2+7m-7$ and $h^1(\mathcal{N}_{\varphi_2})=0$ where $\varphi_2=p_2\circ i$; consequently, $\varphi_2$ is unobstructed. 
\end{itemize}
\end{proposition}
\noindent\textit{Proof.} 
($1$) Follows from $h^j(\mathcal{O}_{X_2})=h^j(\mathcal{O}_Y)\oplus h^j(\mathcal{O}_Y(0,-m-1)$ and K\"unneth formula.\par 
($2$) We apply \cite{Pa}, Corollary 4.1 or \cite{GGP2}, (2.11). Since $Y$ is regular, it follows that $h^0(\mathcal{N}_{p_0})=h^0(\mathcal{O}_Y(0,2m+2))-1=2m+2$. Furthermore, since $H^2(\mathcal{O}_Y)=0$, we obtain $h^1(\mathcal{N}_{p_2})=h^1(\mathcal{O}_Y(0,2m+2))=0$.\par 
($3$) The assertion follows by tensoring the exact sequence ~\eqref{eulerany} by $\mathcal{O}_Y(0,-m-1)$ and taking the long exact sequence of cohomology.\par 
($4$) This is obtained by tensoring ~\eqref{ttnany} by $\mathcal{O}_Y(0,-m-1)$ and taking cohomology.\par 
($5$) As before, the assertion follows from the following exact sequence (see Lemma ~\ref{exactgonzalez})
\begin{equation*}
    0\to\mathcal{N}_{p_2}\to \mathcal{N}_{\varphi_2}\to p_2^*\mathcal{N}_{Y/\mathbb{P}^N}\to 0.
\end{equation*}
Since $h^1(\mathcal{N}_{p_2})=0$, we obtain $h^0(\mathcal{N}_{\varphi_2})=h^0(\mathcal{N}_{p_2})+h^0(\mathcal{N}_{Y/\mathbb{P}^N})+h^0(\mathcal{N}_{Y/\mathbb{P}^N}\otimes\mathcal{O}_Y(0,-m-1))$. We get the value of $h^0(\mathcal{N}_{\varphi_2})$ from part (2), (4) and Lemma ~\ref{cohomology Y}. Finally, $h^1(\mathcal{N}_{\varphi_2})=0$ by part (2), (4) and Lemma ~\ref{cohomology Y}. \QEDB\par

\begin{corollary}\label{2t}
Let $X$ be a  surface of type $(2)_m$ ($m\geq 2$).
\begin{itemize}
    \item[(1)] Suppose $m = 2$, then there exists a smooth, affine irreducible algebraic curve $T$ and a flat family of morphisms $\Phi: \mathcal{X}\to\mathbb{P}^N_T$ over $T$ for which the following happens;
      \begin{itemize}
        \item[(a)] $\Phi_t:\mathcal{X}_t \to \mathbb{P}^{5}$ is a morphism of degree two from a normal projective surface with at worst canonical singularities for all $t \in T-\{0\}$. Further for any $t \in T-\{0\}$, the normalization of Im($\varphi_t$) is a ruled surface over a smooth curve of genus $2$ and has invariant $e = -2$. Further one can take $\Phi_t$ to be the canonical morphism of $\mathscr{X}_t$. 
    
        \item[(b)] $\Phi_0:\mathcal{X}_0 \to \mathbb{P}^{5}$ is the canonical morphism $\varphi:X \to \mathbb{P}^N$.
    \end{itemize}
   \item[(2)] 
The forgetful map from ${\bf Def_{\pi_2/\mathbb{P}^N}} \to {\bf Def_{\varphi}}$ is smooth and hence any deformation of $\varphi$ \textcolor{black}{factors through a double cover of a ruled surface over a curve of genus $m$ and is hence} a morphism of degree bigger than or equal to $2$ onto its image. In particular \textcolor{black}{$\varphi$ cannot be deformed to a birational morphism.}
\end{itemize}
\end{corollary}
\noindent\textit{Proof.} (1) The existence of the curve $T$ so that for a general $t \in T$, $\Phi_t$ has degree two follows from Theorem ~\ref{def to 2:1}, Proposition ~\ref{2m} and \cite{Kaw}. 
Note that since $h^1(T_{\mathbb{P}^5}\vert_Y\otimes\mathcal{O}_Y(0,-3))= h^1(T_{\mathbb{P}^5}\vert_Y)=0$ (by Proposition ~\ref{2m} (3), Lemma ~\ref{cohomology Y} (2)), we have that $h^1(\varphi_2^*(T_{\mathbb{P}^5}))=0$ and ${\bf Def_{\varphi_2}} \to {\bf Def}_{X_2}$ is smooth. Hence the map $H^0(\mathcal{N}_{\varphi_2}) \to H^1(T_{X_2})$ is surjective. By \cite{Sei92}, Lemma $12$, there exist an open set in $H^1(T_{X_2})$ such that for a smooth curve along a first order deformation belonging to the open set a general deformation of $X_2$ along the curve is a ruled surface over a smooth curve of genus $2$ with invariant $e = -2$. Also $H^0(\mathcal{N}_{\varphi_2})\to H^0(\mathcal{N}_{Y/\mathbb{P}^5}\otimes\mathcal{O}_Y(0,-3))$ is surjective and there exist an open set of non-zero elements in $H^0(\mathcal{N}_{Y/\mathbb{P}^5}\otimes\mathcal{O}_Y(0,-3))$. Hence one can choose an element (in fact an open set of elements) from $H^0(\mathcal{N}_{\varphi_2})$ such that it maps to a non-zero element in      $H^0(\mathcal{N}_{Y/\mathbb{P}^5}\otimes\mathcal{O}_Y(0,-3))$ and the general induced deformation $\mathcal{X}_{2t}$ of $X_2$ is a ruled surface over a smooth curve of genus $2$ with invariant $e = -2$.
The trace zero module of the finite map $\pi_2$ is $\mathcal{E}_{\pi_2} = \omega_{X_2} \otimes  \varphi_2^*(\mathcal{O}_{\mathbb{P}^N}(-1))$ (this follows from a computation identical to the one in
the proof of Corollary \ref{11't}  when we showed \eqref{equation.11't}). The fact that $\Phi_t$ can be taken to be the canonical morphism of $X_t$ follows from Theorem ~\ref{def to 2:1}. 

(2) Recall that $X_2=\mathbb{P}^1\times C$ for a smooth curve $C$ that is a double cover $\psi:C\to\mathbb{P}^1$ with $\psi_*\mathcal{O}_C=\mathcal{O}_{\mathbb{P}^1}\oplus \mathcal{O}_{\mathbb{P}^1}(-m-1)$. We have the following splitting of ${\pi_2}_*\mathcal{O}_X$: $${\pi_2}_*\mathcal{O}_X=\mathcal{O}_{X_2}\oplus(\mathcal{O}_{\mathbb{P}^1}(-3)\boxtimes\psi^*\mathcal{O}_{\mathbb{P}^1}(-1)).$$ Also, $T_{X_2}=(\mathcal{O}_{\mathbb{P}^1}\boxtimes\psi^*\mathcal{O}_{\mathbb{P}^1}(1-m))\oplus\mathcal{O}_{\mathbb{P}^1}(2)\boxtimes\mathcal{O}_{C}$, and $\mathcal{E}_{\pi_2}=\mathcal{O}_{\mathbb{P}^1}(-3)\boxtimes\psi^*\mathcal{O}_{\mathbb{P}^1}(-1)$ is the trace zero module of $\pi_2$. It is easy to check that $H^1(T_{X_2}\otimes\mathcal{E}_{\pi_2})=0$. One has the following pullback of the Euler sequence;
\begin{equation*}
    0\to\mathcal{O}_{X_2}\to p_2^*\mathcal{O}_Y(1,m)^{\oplus N+1}\to\varphi_2^*T_{\mathbb{P}^N}\to 0.
\end{equation*}
By tensoring the above exact sequence by $\mathcal{E}_{\pi_2}$ and taking cohomology, we obtain $H^0(\varphi_2^*T_{\mathbb{P}^N}\otimes\mathcal{E}_{\pi_2})=0$. Consequently, the exact sequence;
\begin{equation*}
    0\to T_{X_2}\to\varphi_2^*T_{\mathbb{P}^N}\to\mathcal{N}_{\varphi_2}\to 0.
\end{equation*}
shows that $H^0(\mathcal{N}_{\varphi_2}\otimes\mathcal{E}_{\pi_2})=0$, and the assertion follows from Proposition ~\ref{corweh} (and Corollary ~\ref{maincor}). \QEDB\par 

\vspace{5pt}

The following corollary is analogous to Corollary ~\ref{existence of morphism deformation component $(1)$, $(1')$} and follows immediately from Corollary ~\ref{2t}.
\color{black}
\begin{corollary}\label{existence of morphism deformation component for $(2)$}
Let $X$ be a  surface of type $(2)_m$ $(m\geq 2)$.
\begin{itemize}
\item[(1)] If $m = 2$, there exists an irreducible component $U_{\varphi}$ of the algebraic formally semiuniversal deformation space of ${\varphi}$ (that exists by Remark ~\ref{can}) whose general elements are a two-to-one morphisms onto their image, whose normalization is a ruled surface over a smooth curve of genus $2$ and has invariant $e = -2$.
\item[(2)] There does not exist any component of the algebraic formally semiuniversal deformation space of ${\varphi}$ whose general elements are morphisms that are birational onto their image.
\end{itemize}
\end{corollary}

\color{black}

The following propositions and corollary show that the image of a general morphism in the irreducible component constructed above is smooth. It also shows that this open set intersects the locally closed subloci where the deformed morphism is again the canonical morphism. \par
Recall that for a ruled surface $X \to C$ of invariant $e$ over a smooth curve of genus $g\geq 1$, we denote $C_0$ and $f$ denote the numerical classes of a section and a fibre respectively satisfying $C_0'^2 = -e, f'^2 = 0$ and $C_0' \cdot f' = 1$.

\begin{proposition}\label{deformation of polarized ruled surfaces}
Suppose that there exist a flat family $(\mathcal{X}_2 \to T, \mathcal{L})$ of polarized surfaces ruled over a curve $C$ of genus $m$  over a smooth one dimensional base $T$ with 
\begin{itemize}
\item[(1)] $\mathcal{X}_{20}$ has invariant $0$ and $\mathcal{L}_0 \equiv C_0'+2mf'$
\item[(2)] $\mathcal{X}_{2t}$ has invariant $-m$ or $-(m-1)$ accordingly as $m$ is even or odd.
\end{itemize}
Then after possibly shrinking $T$, $\mathcal{L}_t \equiv C_0'+ \frac{3m}{2}f'$ if $m$ is even and $\mathcal{L}_t \equiv C_0'+ \frac{3(m-1)}{2}f'$ if $m$ is odd.
\end{proposition}
\noindent\textit{Proof.} We prove the statement for $m$ even. The proof is identical for $m$ odd. Let $\mathcal{X}_{20} = X_2$. Let $X_2 \xrightarrow{q} C$ be the smooth morphism. Consider the following commutative diagram.  
\[
\adjustbox{scale=0.8,center}{
\begin{tikzcd}
\dots \arrow[r] & T^1(X_2/C) = H^1(T_{X_2/C}) \arrow[r]\arrow[equal]{d} & T^1(q)\arrow[r]\arrow[d, "\alpha_1"] & H^1(T_C)\arrow[r]\arrow[d, "\beta_1"] & T^2(X_2/C) = H^2(T_{X_2/C})\arrow[r]\arrow[equal]{d} & T^2(q) \arrow[r]\arrow[d,"\alpha_2"] & H^2(T_C) \arrow[r]\arrow[d, "\beta_2"] & \dots \\
\dots \arrow[r] & T^1(X_2/C) = H^1(T_{X_2/C})\arrow[r] & H^1(T_{X_2})\arrow[r] & H^1(q^*(T_C)) \arrow[r] & T^2(X_2/C) = H^2(T_{X_2/C})\arrow[r] & H^2(T_{X_2}) \arrow[r] & H^2(q^*(T_C)) \arrow[r] & \dots
\end{tikzcd}}
\]
Note that since $q_*(\mathcal{O}_{X_2}) = \mathcal{O}_C$, $H^i(q^*(T_C)) = H^i(T_C)$ for $i = 1,2$. Hence the map $\beta_1$ is surjective and the map $\beta_2$ is injective and therefore $\alpha_1$ is surjective and $\alpha_2$ is injective. This implies that the forgetful map ${\bf Def}_{q} \to {\bf Def}_{X_2}$ is smooth. Hence there exist a deformation of $\mathcal{C} \to T$ of $C$ so that $\mathcal{X}_2 \to T$ factors as $\mathcal{X}_2 \xrightarrow{Q} \mathcal{C} \to T$. Fix a line bundle $\mathcal{O}_C(1)$ of degree one on $C$. Since $H^2(\mathcal{O}_C) = 0$, we have that the line bundle $\mathcal{O}_C(1)$ lifts to a line bundle $\mathcal{O}_{\mathcal{C}}(1)$. Now the numerical class of $Q^*(\mathcal{O}_{\mathcal{C}}(1))$ restricts to the numerical class of $f'$ on the central fibre and it is the pullback of a degree one line bundle on $\mathcal{C}_t$ on a general fibre. Hence the numerical class of $f'$ on $X_2$ deforms to the numerical class of $f'$ on $\mathcal{X}_{2t}$. Now suppose that a line bundle of numerical class $C_0'$ in $X_2$ deforms to a line bundle of numerical class $aC_0'+bf'$ on $\mathcal{X}_{2t}$ for $t \neq 0$. Using the fact that their self intersections are the same and noting that $C_0'^2 = 0$ on $X_2$ while $C_0'^2 = m$ on $\mathcal{X}_{2t}$, we have that $a^2m + 2ab = 0$. Suppose that $a = 0$. Then for sufficiently large $k$, $C_0'+kf'$ which is very ample on $X_{20}$ deforms to $(k+b)f'$ which is not ample. Hence $a \neq 0$ and $b =  \frac{-am}{2}$. Then $C_0'$ on $X_{20}$ deforms to $a(C_0'-\frac{m}{2}f')$ on $\mathcal{X}_{2t}$. Considering that on $\mathcal{X}_{20}$, $C_0' \cdot f' = 1$, we have that $a = 1$. Hence our statement is proven.  \QEDB

\vspace{5pt}

We prove a slightly stronger version of a result we need to prove Corollary
~\ref{smoothness of image for (2)_2}. More precisely we will use the result proven below for $m = 2$.
\begin{proposition}\label{very ampleness on ruled surfaces}
Suppose $X$ is a ruled surface over a curve $C$ of genus $m$ with invariant $-m$ or $-(m-1)$ accordingly as $m$ is even or odd. Then a line bundle $L\equiv C_0+ \frac{3m}{2}f$ is very ample if $m$ is even and a line bundle $L\equiv C_0+ \frac{3(m-1)}{2}f$ is very ample if $m$ is odd.
\end{proposition}

\noindent\textit{Proof}. We use the following criterion for very ampleness (see \cite{LM05} Corollary $2.13$): let $|H|$ be the complete linear series of a line bundle $H \equiv C_0+bf$. Then $|H|$ is very ample if and only for any two points $P$ and $Q$ on $C$, $h^0(H-(P+Q)f) = h^0(H)-4$. \par 
Let $m$ be even and consider $L\equiv C_0+\frac{3m}{2}f$ on a ruled surface with invariant $-m$ over a curve of genus $m$ with $m \geq 2$. We need to show, for any two points $P$ and $Q$ on $C$,
\begin{equation*}
    h^0(L-(P+Q)f)=h^0(L)-4.
\end{equation*}
Let $X = \mathbb{P}(\mathcal{E})$ where $\mathcal{E}$ is normalized and since $e < 0$, we have that $\mathcal{E}$ is stable. Since higher pushforward of any bundle of numerical equivalence class $C_0+bf$ is zero we can compute the above cohomology by pushing forward to the base curve of genus $m$. Notice that $L-(P+Q)f\equiv C_0+(\frac{3m}{2}-2)f$. 
Hence it is enough to show that for any two line bundles $L_1$ and $L_2$ on $C$ with degree $3m/2-2$ and $3m/2$ respectively, we have
\begin{equation*}
h^0(L_1\otimes\mathcal{E})=h^0(L_2\otimes\mathcal{E})-4.
\end{equation*}
For simplicity, let us denote $L_1$ by $\mathcal{O}_C(\frac{3m}{2}-2)$, and $L_2$ by $\mathcal{O}_C(\frac{3m}{2})$. Applying Riemann-Roch to the vector bundles $\mathcal{O}_C(\frac{3m}{2}-2) \otimes \mathcal{E}$ and $\mathcal{O}_C(\frac{3m}{2}) \otimes \mathcal{E}$ and subtracting we get 
\begin{equation*}
   h^0\left(\mathcal{E}\left(\frac{3m}{2}\right)\right)- h^0\left(\mathcal{E}\left(\frac{3m}{2}-2\right)\right) = h^1\left(\mathcal{E}\left(\frac{3m}{2}\right)\right)- h^1\left(\mathcal{E}\left(\frac{3m}{2}-2\right)\right) + c_1\left(\mathcal{E}\left(\frac{3m}{2}\right)\right) - c_1\left(\mathcal{E}\left(\frac{3m}{2}-2\right)\right).
\end{equation*}
Since $\mathcal{E}$ is of rank two we have that $c_1(\mathcal{E}(\frac{3m}{2}) - c_1(\mathcal{E}(\frac{3m}{2}-2)) = 4$. Hence we are done if we show $$h^1\left(\mathcal{E}\left(\frac{3m}{2}\right)\right)= h^1\left(\mathcal{E}\left(\frac{3m}{2}-2\right)\right) = 0.$$ \par 
Note that $h^1(\mathcal{E}(\frac{3m}{2})) = h^0(\mathcal{E}^*(\frac{m}{2}-2)$. Now the slope of the vector bundle $\mu(\mathcal{E}^*(\frac{m}{2}-2)) = -2$. Also since $\mathcal{E}$ is stable we have that $\mathcal{E}^*(\frac{m}{2}-2)$ is stable. Since its slope is negative we have that $h^0(\mathcal{E}^*(\frac{m}{2}-2)) = 0$. \par
Now note that $h^1(\mathcal{E}(\frac{3m}{2}-2)) = h^0(\mathcal{E}^*(\frac{m}{2}))$. Note that $\textrm{deg}(\mathcal{E}^*(\frac{m}{2})) = 0$ and $\mathcal{E}^*(\frac{m}{2})$
is stable since $\mathcal{E}$ is stable. Then $h^0(\mathcal{E}^*(\frac{m}{2})) = 0$ since for degree $0$ vector bundles the existence of a section contradicts stability. 
The proof for the case $m$ odd follows exactly along the same lines. \QEDB

\begin{corollary}\label{smoothness of image for (2)_2}
In Corollary ~\ref{2t} (1), we can choose the curve $T$ so that after possibly shrinking $T$, for $t \in T$, $t \neq 0$, $\textrm{Im}(\Phi_t)$ is smooth and $\Phi_t$ can be taken to be the canonical morphism of $\mathcal{X}_t$.
\end{corollary}

\noindent\textit{Proof.} We resume notations of Corollary ~\ref{2t}. Consider the factorization $X \xrightarrow{\pi_2} X_2 \xrightarrow{p_2} Y \xrightarrow{i} \mathbb{P}^N$. Let $\varphi_2 = i \circ p_2$. Note that since $X_2 = C \times \mathbb{P}^1$, where $C$ is a smooth curve of genus $2$, we have that it is a ruled surface over $C$ and has invariant $e = 0$. Consider $L = \varphi_2^*(\mathcal{O}_{\mathbb{P}^N}(1)) = p_2^*(\mathcal{O}_Y(C_0+2f))$ (recall $C_0$ and $f$ are the classes of a section and fibre of $Y$). Since $C_0^2=f^2=0$ and $C_0 \cdot f=1$, we have that $p_2^*(\mathcal{O}_Y(C_0)) \equiv aC_0'$, {and} $p_2^*(\mathcal{O}_Y(f)) \equiv bf'$, with $ab = 2$. Then since $h^0(p_2^*\mathcal{O}_Y(C_0)) = 2$ we have that $a = 1$ and hence $b = 2$. Hence $p_2^*(\mathcal{O}_Y(C_0+2f)) \equiv C_0'+4f'$. 
\par 
Note that the pair $(X_2,L)$ is unobstructed since $h^2(\mathcal{E}_L) = 0$ (since $h^2(\mathcal{O}_{X_2}) = h^2(T_{X_2}) = 0$). Also since $h^2(\mathcal{O}_{X_2}) = 0$, ${\bf Def}_{(X_2,L)} \to {\bf Def}_{X_2}$ is smooth. Choose a smooth curve $T$ from the smooth versal deformation space of $(X_2,L)$. Let $(\mathcal{X}_2 \xrightarrow{\sigma} T, \mathcal{L})$ be the family obtained. Then for a general such curve, for $t \neq 0$, $\mathcal{X}_{2t}$ has invariant $-2$. By Proposition ~\ref{deformation of polarized ruled surfaces}, $\mathcal{L}_t \equiv C_0'+3f'$ which is very ample by Proposition ~\ref{very ampleness on ruled surfaces}. Since $H^1(L)=0$, (easy to check by projection formula) we have that (after shrinking $T$), $\sigma_*(\mathcal{L})$ is locally free of rank $h^0(L)$ and  we get a morphism $\mathcal{X}_2 \xrightarrow{\Phi_2} \mathbb{P}(\sigma_*(\mathcal{L})) \to T$ which is an embedding for $t \neq 0$ since it is given by the complete linear series of a line bundle numerically equivalent to $C_0'+3f'$, which is very ample. Note that $\mathcal{E}_{\pi_2} = \omega_{X_2} \otimes \varphi_2^*(L^*)$ and let $B \in |(\omega_{X_2} \otimes \varphi_2^*(L^*))^{-2}|$ be the divisor giving $\pi_2$. Note that $(\omega_{\mathcal{X}_2/T}^{-2} \otimes \Phi_2^*(\mathcal{O}_{\mathbb{P}_T^N}(2)))$ is a lift of $(\omega_{X_2} \otimes \varphi_2^*(L^*))^{-2}$. Then by Remark ~\ref{jayan} one can construct a lift $\mathcal{B}$ of $B$ and hence a relative double cover $\Pi_2$ since $H^1(\mathcal{O}_{X_2}(B))=0$ (easy to check, see the proof of Proposition ~\ref{moduli'2} (1)). Consider $\Phi = \Pi_2 \circ \Phi_2$. For $t \neq 0$, $\Phi_t$ is the composition of a double cover $\Pi_{2t}$ followed by an embedding $\Phi_{2t}$ of a smooth surface given by the complete linear series of a very ample line bundle $\mathcal{L}_t$. Moreover $\Pi_{2t}$ is branched along $(\omega_{\mathcal{X}_{2t}} \otimes \mathcal{L}_t^*)^{-2}$. Hence $\Phi_t$ is the canonical morphism of $\mathcal{X}_t$ (by Remark ~\ref{trace zero module for canonical morphism}) and its image is smooth. \QEDB

\color{black}
\subsection{Deformation of canonical morphism for type \texorpdfstring{$(3)_m$}{TEXT}}\label{defproduct} \textcolor{black}{Recall} that surfaces of type $(3)_m$ are of the form $C_1\times C_2$ where $C_1$ is a smooth hyperelliptic curve of genus $2$ and $C_2$ is a smooth hyperelliptic curve of genus $m+1$. 
\color{black} The deformations of $\varphi$ can be studied with 
the same machinery used so far. Thus, they fit in the same theoretical framework as the deformations of the other three types of quadruple Galois covers. However, product of curves can be dealt with in a different, ad-hoc, easier way. We start with a consequence of 
Beauville's numerical characterization of minimal surfaces of general type satisfying $p_g = 2q-4$ (see 
\cite{BeauD}): 

\begin{corollary}\label{cor.Beauville}
Let $C_1'$ be a curve of genus $2$, let $C_2'$ be a curve of genus $g$, $g \geq 2$, let $X'=C_1' \times C_2'$ and let $\mathcal M_{(p_g,q,K^2)}$ be the moduli space of minimal surfaces of general type with invariants 
$(p_g,q,K^2)$ to which $[X']$ belongs. Then 
\begin{enumerate}
    \item  $\mathcal M_{(p_g,q,K^2)}= \mathcal M_{(2g,g+2,8g-8)}=
    \mathcal M_{2} \times \mathcal M_{g}$ if $g > 2$; and
    \item  $\mathcal M_{(p_g,q,K^2)}=\mathcal M_{(4,4,8)}=   \mathrm{Sym}^2(\mathcal M_{2})$ if $g=2$.
\end{enumerate}
In particular, $\mathcal M_{(p_g,q,K^2)}$ is irreducible. Moreover, the canonical morphism of $X'$ is the composition of  $\psi_1' \times \psi_2'$, where $\psi_i'$ is the canonical morphism of $C_i'$, and the Segre embedding of $\mathbb P^1 \times \mathbb P^{g-1}$. 
\end{corollary}

\noindent\textit{Proof.} The corollary follows from the fact that $\omega_{X'}=\omega_{C_1'} \boxtimes \omega_{C_2'}$, Kunneth formula, Beauville's  characterization of minimal  surfaces of general type satisfying $p_g = 2q-4$ (see 
\cite{BeauD}) and \cite{Ops}, Theorems $4.1$, $4.2$.\QEDB

\color{black}
\begin{proposition}\label{prop.product.of.curves}
Let $X$ be a (smooth) surface of type $(3)_m$ and let $U_{\varphi}$ be the algebraic formally semiuniversal deformation space of ${\bf Def}_{\varphi}$. For any $[(X',\varphi')]$ of  $U_{\varphi}$, we have 
$X'=C_1' \times C_2'$, where $C_1'$ is a curve of genus $2$ and $C_2'$ is a curve of genus $g$, $g \geq 2$, \textcolor{black}{and 
 $\varphi'$ is the canonical morphism of $X'$, so $U_\varphi$ is irreducible.} In particular:
\begin{itemize}
    \item[(1)] If $m=1$, then any deformation of $\varphi$ is a morphism of degree $4$ onto its image, which is isomorphic to $\mathbb{P}^1 \times \mathbb{P}^1$.  
    \item[(2)] If $m\geq 2$, then any deformation $\varphi'$ of $\varphi$ has degree $2$ or $4$ (in particular, $\varphi'$ is not birational onto its image). If $\varphi'$ is general, then $\varphi'$ has degree $2$. 
\end{itemize} 
\end{proposition}

\color{black}
\noindent\textit{Proof.} 
The space $U_{\varphi}$ admits a dominant morphism to $\mathcal M_{(2g,g+2,8g-8)}$, so
Corollary~\ref{cor.Beauville} implies $X'=C_1' \times C_2'$, where $C_1'$ is a curve of genus $2$ and $C_2'$ is a curve of genus $g$, $g \geq 2$. 
By abuse of notation we will call $C_2$ any of the fibers of the projection from $X$ to $C_1$. 
{\textcolor{black}{For any $[(X',\varphi')] \in U_\varphi$, a fiber of the projection from $X'$ to $C_1'$, which we will call $C_2'$ by the same abuse of notation, is a deformation of some of the $C_2$, and $\varphi'|_{C_2'}$ is the deformation of 
$\varphi|_{C_2}$, which, by the K\"unneth formula, is the canonical morphism of $C_2$. Then $\varphi'|_{C_2'}$ is the canonical morphism of $C_2'$, so the restriction of $L'=\varphi'(\mathcal O_{\mathbb P^N}(1)$ to each fiber of the projection to $C_1'$ is the canonical bundle of the fiber.}}
Then (see e.g. \cite[\S 11.5]{BL} and \cite[Proposition 3.3.8]{Smith}),  $L' \otimes \omega_{X'}^*$ is the box product of a line bundle of $C_1'$ and a line bundle of $C_2'$; consequently, $L'=(\omega_{C_1'} \otimes \delta_1) \boxtimes 
(\omega_{C_2'} \otimes \delta_2)$, where $\delta_i'$ is a degree $0$ line bundle on $C_i'$. Since  $\varphi'$ is induced by $|L'|$, 
we have, by K\"unneth formula, $2g=h^0((\omega_{C_1'} \otimes \delta_1) \boxtimes 
(\omega_{C_2'} \otimes \delta_2))=h^0(\omega_{C_1'} \otimes \delta_1)\cdot 
h^0(\omega_{C_2'} \otimes \delta_2)$. Then $h^0(\omega_{C_1'} \otimes \delta_1)=2$ and 
$h^0(\omega_{C_2'} \otimes \delta_2)=g$, so $\delta_1=\delta_2=0$. Thus 
$\varphi'$ is the composition of $\psi_1' \times \psi_2'$, where
$\psi_i'$ is the canonical morphism of $C_i'$, and the Segre embedding 
of $\mathbb P^1 \times \mathbb P^{g-1}$. Therefore, $\varphi'$ is the canonical morphism of $X'$. Then $U_\varphi$ is birational to $\mathcal M_{(2g,g+2,8g-8)}$, hence irreducible (see Corollary~\ref{cor.Beauville}). Since $\varphi'$ factors through 
$\psi_1' \times \psi_2'$ and the Segre embedding, (1) and (2) are straight-forward. \QEDB

 \begin{remark}\label{existence of morphism deformation component $(3)$, }
 Using arguments similar to the ones employed in \S 3.1 and \S 3.2, one can prove that, if $X$ is of type $(3)_m$, $m \geq 2$, then any deformation of $\varphi$ factors through a double cover of a ruled surface over a curve of genus $m+1$.
 \end{remark}

\color{black}

\vspace{5pt}

\noindent\textit{Proof of Theorem ~\ref{A}.} It follows immediately from Corollary ~\ref{existence of morphism deformation component $(1)$, $(1')$}, Corollary ~\ref{non-normality for $(1)$, $(1')$}, Corollary ~\ref{existence of morphism deformation component for $(2)$}, Corollary ~\ref{smoothness of image for (2)_2}, \textcolor{black}{Proposition~\ref{prop.product.of.curves}, Remark~\ref{existence of morphism deformation component $(3)$, } and Theorem~\ref{B}}.\QEDB

\color{black}

\section{Moduli components of irregular covers of surface scrolls}\label{4}
{In this section we will study the moduli components of irregular quadruple covers of minimal degree.}
 Furthermore, if the cover is unobstructed, we know that there is a unique component of the moduli of surfaces of general type; in that case we would like to understand the geometry of this moduli component. Regarding  surfaces $X$ of type $(1)_m$ and $(1')_m$, our result is as follows.

\begin{theorem}\label{B} Let $X$ be a surface of type $(1)_m$ or $(1')_m$. Then
    \begin{itemize}
     \item[(1)] the functors ${\bf Def}_{\varphi}$, where $\varphi$ is the canonical morphism of $X$, and ${\bf Def}_X$ are smooth. Hence in particular any $X$ of type $(1)_m$ or $(1')_m$ is contained in a unique irreducible component of the moduli of surfaces of general type; 
     \item[(2)] there exists a unique irreducible component of the moduli of surfaces of general type $\mathcal{M}_{(8m,1,2m+2)}$ containing all surfaces of both types. This component is uniruled of dimension $8m+20$; and
     \item[(3)] the canonical morphism of a general element of this component is a two-to-one morphism onto its image which is a non-normal variety whose normalization is an elliptic ruled surface which is the projectivization of a rank two split vector bundle over an elliptic curve and has invariant $e = 0$. 
     
   \end{itemize}
\end{theorem}

The situation is not as clean as the previous theorem for general surfaces of type $(2)_m$ for $m\geq 3$, even though the result is neat for $m=2$. In particular, we show the following.

\begin{theorem}\label{C}
Let $X$ be a  surface of type $(2)_m$. Then
\begin{itemize}
    \item[(1)] If $m=2$, then ${\bf Def}_{\varphi}$, where $\varphi$ is the canonical morphism of $X$, and ${\bf Def}_X$ are smooth. Hence in particular $X$ is contained in a unique irreducible component of the moduli of surfaces of general type. Furthermore, 
    \begin{itemize} 
     \item[(a)] there exists a unique irreducible component of the moduli space of surfaces of general type $\mathcal{M}_{(16,2,6)}$ containing $X$ (and all other surfaces of type $(2)_2$). This component is uniruled of dimension 28, and
     \item[(b)] the canonical morphism of a general element in that component is a double cover onto its image 
     whose normalization is a ruled surface over a smooth curve of genus $2$ and has invariant $e = -2$.
    \end{itemize}
    \item[(2)] If $m\geq 3$, there do not exist an irreducible component of the moduli of surfaces of general type $\mathcal{M}_{(8m,m,2m+2)}$ containing $X$, such that the canonical morphism of a general element in that component is birational onto its image.
\end{itemize}
\end{theorem}

The proofs of Theorems~\ref{B} and \ref{C} are quite involved and 
we devote the next two subsections to them. In contrast, surfaces of type $(3)_m$, since they are product of two curves, are much easier to handle:

\color{black}
\begin{remark}\label{remark.Ops}
Recall that, if $X$ is a surface of type $(3)_m$, then  $X$ is a product of smooth hyperelliptic curves of genus $2$ and $m+1$. Then, as pointed out in Corollary~\ref{cor.Beauville}, $X$ belongs to 
the moduli space $\mathcal M=\mathcal M_{(2m+2,m+3,8m)}$, which is irreducible.
We  also saw in Corollary~\ref{cor.Beauville} that $\mathcal M_{(4,4,8)}=\textrm{Sym}^2(\mathcal{M}_{2})$ and that, if $m \geq 2$, $\mathcal M=\mathcal M_2 \times \mathcal M_{m+1}$.
 It follows that $\mathcal{M}$ is of dimension $3m+3$ and is uniruled, since $\mathcal{M}_2$ is uniruled (in fact, rational, see \cite{Bog}, \cite{Kat}). By Corollary~\ref{cor.Beauville}, if $m \geq 2$, then the canonical morphism of  a surface,  general in $\mathcal M$, is of degree 2, and its image  is isomorphic to $\mathbb{P}^1\times C'_2$, where $C_2'$ is a curve of genus $m+1$. If $m=1$, the canonical morphism of any surface of $\mathcal M$ is a quadruple cover of  $\mathbb{P}^1\times\mathbb{P}^1$ (see again Corollary~\ref{cor.Beauville}).
\end{remark}

\color{black}

In contrast to Theorems ~\ref{B} $(2)$, ~\ref{C} $(2)$ 
\textcolor{black}{and Remark~\ref{remark.Ops}}, it is shown in \textcolor{black}{\cite{Cat81}, \cite{Cat87}, \cite{CS02},  \cite{CPT00}}, \cite{GGP2} \textcolor{black}{(see also \cite{AK90} and \cite{GGP13c})} 
that there are irreducible components of moduli spaces of surfaces of general type such that the canonical morphism of a general element is birational onto its image. 
\textcolor{black}{At the same time, in these components there is a proper sublocus in which the canonical morphism has degree $2$.}

\subsection{Description of moduli components of surfaces of types $(1)$ and $(1')$} First we aim to prove Theorem ~\ref{B}. Throughout this subsection, we work with the notations of \S ~\ref{study11'}. Let $B\in |\psi^*\mathcal{O}_{\mathbb{P}^1}(2)\boxtimes\mathcal{O}_{\mathbb{P}^1}(2m+4)|$ be the branch divisor of $\pi_1$. In order to do that, we need the following cohomology computations.
\begin{proposition}\label{moduli11'}
Let $X$ be a surface of type $(1)_m$ or $(1')_m$. 

\begin{itemize}
    \item[(1)] $h^0(\mathcal{O}_{X_1}(B))=8m+20$ and    $h^1(\mathcal{O}_{X_1}(B))=h^2(\mathcal{O}_{X_1})=0$.
    \item[(2)] $h^0(\pi_1^*T_{X_1})=4$, $h^1(\pi_1^*T_{X_1})=4$ and $h^2(\pi_1^*T_{X_1})=4m$. 
    \end{itemize}
    Furthermore, if $X$ is smooth (and, hence, of type $(1'_m)$), then:
    \begin{itemize}
    \item[(3)] $h^0(\mathcal{N}_{\pi_1})=8m+20 $ and   $h^1(\mathcal{N}_{\pi_1})=h^2(\mathcal{N}_{\pi_1})=0$.
    \item[(4)] $h^1(T_X)=8m+20 $, $h^2(T_X)=4m$.
    \item[(5)] $h^0(\mathcal{N}_{\varphi})=4m^2+16m+24$, and  $h^1(\mathcal{N}_{\varphi})\geq 2m+2$.
\end{itemize}
\end{proposition}
\noindent\textit{Proof.} (1) It is easy to see from K\"unneth formula, and projection formula that 
$$h^0(\mathcal{O}_{X_1}(B))=8m+20,\textrm{ } h^1(\mathcal{O}_{X_1}(B))=0,\textrm{ and } h^2(\mathcal{O}_{X_1}(B))=0.$$
Then $(3)$ follows from $(1)$ and the following exact sequence 
\begin{equation*}
    0\to\mathcal{O}_{X_1}\to \mathcal{O}_{X_1}(B)\to\mathcal{O}_B(B)\to 0,
\end{equation*}
\cite{Pa}, Corollary 4.1 or \cite{GGP2}, (2.11), and Proposition ~\ref{11'} (1).\par 

(2): One checks this readily by Proposition ~\ref{11'} (1), K\"unneth formula, and projection formula since $h^j(\pi_1^*T_{X_1})$ is nothing but the following sum
$$h^j(\mathcal{O}_{X_1})+h^j(\mathcal{O}_E\boxtimes\mathcal{O}_{\mathbb{P}^1}(2))+h^j(\psi^*\mathcal{O}_{\mathbb{P}^1}(-1)\boxtimes\mathcal{O}_{\mathbb{P}^1}(-m-2))+h^j(\psi^*\mathcal{O}_{\mathbb{P}^1}(-1)\boxtimes\mathcal{O}_{\mathbb{P}^1}(-m)).$$
\indent (4): We use the following exact sequence
\begin{equation}
    0\to T_X\to \pi_1^*T_{X_1}\to\mathcal{N}_{\pi_1}\to 0.
\end{equation}
Since $h^0(T_X)=0$, and $h^1(\mathcal{N}_{\pi_1})=h^2(\mathcal{N}_{\pi_1})=0$ by part (3), we get the following two exact sequences:
$$0\to H^0(\pi_1^*T_{X_1})\to H^0(\mathcal{N}_{\pi_1})\to H^1(T_X)\to H^1(\pi_1^*T_{X_1})\to 0,$$
$$0\to H^2(T_X)\to H^2(\pi_1^*T_{X_1})\to 0.$$
The conclusion now follows from parts (2) and (3).\par 

(5): We get the following exact sequence from Lemma ~\ref{exactgonzalez}:
\begin{equation}\label{exactgonzalezpi1}
    0\to \mathcal{N}_{\pi_1}\to\mathcal{N}_{\varphi}\to\pi_1^*\mathcal{N}_{\varphi_1}\to 0.
\end{equation}

We first compute $h^0(\mathcal{N}_{\varphi})$. From part (1), we get $h^1(\mathcal{N}_{\pi_1})=0$. It follows from projection formula that $$h^0(\mathcal{N}_{\varphi})=h^0(\mathcal{N}_{\pi_1})+h^0(\mathcal{N}_{\varphi_1})+h^0(\mathcal{N}_{\varphi_1}\otimes\mathcal{E}_{\pi_1}).$$
Recall that we have checked the vanishing of $h^0(\mathcal{N}_{\varphi_1}\otimes\mathcal{E}_{\pi_1})$ in the proof of Corollary ~\ref{11't}. Thus, $$h^0(\mathcal{N}_{\varphi})=8m+20+(N+1)^2$$ by part (3) and Proposition ~\ref{11'} (5). The conclusion follows since $N+1=2m+2$.\par 

Notice that $h^2(\mathcal{N}_{\pi_1})=0$ by part (3). From ~\eqref{exactgonzalezpi1}, we obtain that $h^1(\mathcal{N}_{\varphi})=h^1(\mathcal{N}_{\varphi_1})+h^1(\mathcal{N}_{\varphi_1}\otimes\mathcal{E}_{\pi_1})$. The conclusion follows from Remark ~\ref{nonzero!!}.\QEDB\par 

\vspace{5pt}

\noindent\textit{Proof of Theorem ~\ref{B}.} We resume the notations of \S ~\ref{study11'}. Fix a  surface $X$ of type $(1)_m$ or $(1')_m$.  First note that ${\bf Def}_{\varphi}$ and ${\bf Def}_{X}$ have algebraic formally semiuniversal deformation spaces by Remark ~\ref{can}.\par 

We apply Corollary ~\ref{Defphi}. All the hypotheses have been verified in the proofs of Proposition ~\ref{moduli11'} and Corollary ~\ref{11't}. It follows that ${\bf Def}_{\varphi}$ is unobstructed.\par

We first show that ${\bf Def}_X$ is unobstructed for a smooth surface of family $(1')_m$. We will show unobstructedness of ${\bf Def}_X$ for a singular surface of family $(1')_m$ and any surface of family $(1)_m$ after we prove part $(2)$. We use Corollary ~\ref{defcanonicalmorphism }. It remains to verify that the differential of the map ${\bf Def}_{\varphi}\to{\bf Def}_{X}$ is surjective. Since the canonical bundle $\omega_X$ lifts to any first order deformation of $X$, it is enough to show that any section of $\omega_X$ lifts to any first order deformation of $(X,\omega_X)$. We aim to use the section lifting criterion (Proposition ~\ref{sectionlifting}). Consider the Atiyah extension ~\eqref{atiyahg}
\begin{equation}\label{mainatiyah}
    0\to \mathcal{O}_X\to\mathcal{E}_{\omega_X}\to T_X\to 0.
\end{equation}
Since the canonical bundle $\omega_X$ lifts to any first order deformation of $X$, $H^1(\mathcal{E}_{\omega_X})\to H^1(T_X)$ is surjective. This, together with $H^0(T_X)=0$ implies that $$h^0(\mathcal{E}_{\omega_X})=1 \textrm{ and }h^1(\mathcal{E}_{\omega_X})=h^1(\mathcal{O}_X)+h^1(T_X)=8m+21$$ where the last equality follows from Proposition ~\ref{moduli11'} (3) and Theorem \ref{gpmain}. Now consider the following exact sequence ~\eqref{atiyah'g} 
\begin{equation}\label{atiyah}
    0\to\mathcal{E}_{\omega_X}\to H^0(\omega_X)^{\vee}\otimes \omega_X\to \mathcal{N}_{\varphi}\to 0.
\end{equation}
Recall that a section in $H^0(\omega_X)$ lifts to a first order deformation $\eta\in H^1(\mathcal{E}_{\omega_X})$ of the pair $(X,\omega_X)$ if and only if its image under the map $H^1(\mathcal{E}_{\omega_X})\to \textrm{Hom}(H^0(\omega_X),H^1(\omega_X))$ induced from ~\eqref{atiyah} is zero (see Proposition ~\ref{sectionlifting}). Thus, it is enough to show that the map $H^1(\mathcal{E}_{\omega_X})\to H^0(\omega_X)^{\vee}\otimes H^1(\omega_X))$ induced from ~\eqref{atiyah} is zero. Now, ~\eqref{atiyah} gives rise to the following exact sequence 
\begin{equation}\label{sectionlift}
    0\to H^0(\mathcal{E}_{\omega_X})\to H^0(\omega_X)^{\vee}\otimes H^0(\omega_X)\to H^0(\mathcal{N}_{\varphi})\to H^1(\mathcal{E}_{\omega_X}).
\end{equation}
Thus, the dimension of the image of $H^0(\mathcal{N}_{\varphi})\to H^1(\mathcal{E}_{\omega_X})$ is 
$$ h^0(\mathcal{E}_{\omega_X})-(h^0(\omega_X))^2+h^0(\mathcal{N}_{\varphi})=1-(2m+2)^2+4m^2+16m+24=8m+21$$
where the last equality follows from Proposition ~\ref{moduli11'} and Theorem \ref{gpmain}.  But this dimension is the same as $h^1(\mathcal{E}_{\omega_X})$. This shows that ~\eqref{sectionlift} is surjective on the right, and consequently any section of $\omega_X$ lifts to any first order deformation of $(X,\omega_X)$. Thus the differential of ${\bf Def}_{\varphi}\to{\bf Def}_{X}$ is surjective. Thus ${\bf Def}_X$ is unobstructed, and the algebraic formally semiuniversal (in fact universal) deformation space $U_X$ of this functor is smooth, irreducible and uniruled.\par 

(2): We now show that there exist a unique component of moduli of surfaces of general type containing {\it all} surfaces of both types $(1)_m$ and $(1')_m$. We do this by the following few steps.

\smallskip

\noindent\underline{Step 1.} We claim that {\it all} bidouble covers i.e, the surfaces of the family $(1')_m$, are contained in an irreducible component of the moduli. We will show this by showing that the bidouble covers in $(1')_m$ are parametrized by
a subset {and, smooth ones, by  a (non-empty) open set} of $\mathbb{P}(H^0(2C_0+(2m+4)f)) \times \mathbb{P}(H^0(4C_0))$, {which is irreducible}. 
\par
Let $L_1 = -C_0-(m+2)f$ and $L_2 = -2C_0$ on $Y = \mathbb{P}^1 \times \mathbb{P}^1$. Note that $Y$ is rigid and consider the following Cartesian square. 
\[ 
\begin{tikzcd}
Y \times \mathbb{P}(p_*(L_1^{\otimes -2})) \times_k \mathbb{P}(p_*(L_2^{\otimes -2}))  \arrow[r] \arrow[d,"q"] & \mathbb{P}(p_*(L_1^{\otimes -2})) \times_k \mathbb{P}(p_*(L_2^{\otimes -2})) \arrow[d] \\
Y \arrow[r,"p"] & \textrm{Spec}(\mathbb{C}) 
\end{tikzcd}
\]
Furthermore, consider the divisor $\mathscr{B}_1 = \{(x,r,s) \in Y \times \mathbb{P}(p_*(L_1^{\otimes -2})) \times_k \mathbb{P}(p_*(L_2^{\otimes -2})) | r(x) = 0 \} = q^*(L_1^{\otimes -2})$ and $\mathscr{B}_2 = \{(x,r,s) \in Y \times \mathbb{P}(p_*(L_1^{\otimes -2})) \times_k \mathbb{P}(p_*(L_2^{\otimes -2})) | s(y) = 0 \} = q^*(L_2^{\otimes -2})$. Let $$T = Y \times \mathbb{P}(p_*(L_1^{\otimes -2})) \times_k \mathbb{P}(p_*(L_2^{\otimes -2})).$$ Let $\mathscr{L}_1^{\otimes -1} \xrightarrow{f_1} T$ and $\mathscr{L}_2^{\otimes -1} \xrightarrow{f_2} T$ denote the total spaces of the line bundles $q^*(L_1^{\otimes -1})$ and $q^*(L_2^{\otimes -1})$ on $T$. Moreover, let $t_1 \in H^0(f_1^*q^*(L_1^{\otimes -1}))$ and $t_2 \in H^0(f_2^*q^*(L_2^{\otimes -1}))$ be the corresponding tautological sections. Then one can consider relative double covers on $T$ given as the zero locus of $t_i^2-f_i^*\mathscr{B}_i$ inside $\mathscr{L}_i^{\otimes -1}$.
\[ 
\begin{tikzcd}
T_1 = (t_i^2-f_i^*\mathscr{B}_i)_0 \arrow[r, hookrightarrow] \arrow[dr]  & \mathscr{L}_i^{\otimes -1} \arrow[d]  \\
& T 
\end{tikzcd}
\]
Now consider the fibre product of the relative double covers $T_1$ and $T_2$ over $T$. Consider the flat family $T_1 \times_T T_2 \to T \to \mathbb{P}(p_*(L_1^{\otimes -2})) \times_k \mathbb{P}(p_*(L_2^{\otimes -2}))$. Pulling back the composed morphism at a point $(r,s) \in \mathbb{P}(p_*(L_1^{\otimes -2})) \times_k \mathbb{P}(p_*(L_2^{\otimes -2}))$ we have the following Cartesian square
\[ 
\begin{tikzcd}
Y_r \times_Y Y_s \arrow[r, hookrightarrow] \arrow[d]  & T_1 \times_T T_2 \arrow[d]  \\
Y \arrow[r, hookrightarrow] \arrow[d] & T \arrow[d] \\
(r,s) \arrow[r, hookrightarrow] & \mathbb{P}(p_*(L_1^{\otimes -2})) \times_k \mathbb{P}(p_*(L_2^{\otimes -2})) 
\end{tikzcd}
\]
where $Y_r$ and $Y_s$ denote the double covers constructed by $r \in  {|L_1^{\otimes -2}|}$ and $s \in  {|L_2^{\otimes -2}|}$. 
{{By Theorem ~\ref{gpmain},  each surface of family $(1')_m$  is a fiber product $Y_r \times_Y Y_s$,  for uniquely determined divisors $r$ in   $|L_1^{\otimes -2}|$ and $s$ in  $|L_2^{\otimes -2}|$, so surfaces belonging to family $(1')_m$ are parametrized by a subset of $\mathbb{P}(p_*(L_1^{\otimes -2})) \times_k \mathbb{P}(p_*(L_2^{\otimes -2}))$.   A surface of family $(1')_m$ is smooth if and only the divisors $r$  and $s$ are smooth and meet transversally. Such divisors exist by Bertini. Then, generic smoothness implies that
 smooth surfaces  are parametrized by a {(non-empty)} open set} of $\mathbb{P}(p_*(L_1^{\otimes -2})) \times_k \mathbb{P}(p_*(L_2^{\otimes -2}))$.}

\smallskip

\noindent\underline{Step 2.} In this step we note that there is a {\it unique} irreducible component, say $D$, of the moduli containtaing {\it all} surfaces of type $(1')_m$. This comes from Step 1 and the unobstructedness of the smooth surfaces of type $(1')_m$. 

\smallskip

\noindent\underline{Step 3.} We claim that {\it any} cyclic cover in $(1)_m$ can be  deformed to a smooth bidouble cover along an irreducible curve whose general fibre parametrizes smooth surfaces of general type. In particular, given any cyclic cover, there exists an irreducible component $D'$ containing the cyclic cover and a smooth bidouble cover. Indeed, take the cyclic cover $X$. Its intermediate cover $X_1$ is smooth. The cyclic cover is obtained by a (special) choice of branch divisor from $p_1^*(2C_0+2(m+2)f)$. Since $X_1$ is smooth, by a different (special) choice of branch divisor, one can construct a smooth bidouble cover over $X_1$. Since $\mathbb{P}(H^0(p_1^*(2C_0+2(m+2)f)))$ is irreducible and a general member is smooth by Bertini, we have that one can deform the cyclic cover to a smooth bidouble cover along a curve whose general member is a smooth  surface of general type (in this case a smooth double cover over an elliptic ruled surface).

\smallskip

\noindent\underline{Step 4.} We claim that $D' = D$. Indeed, if $D' \neq D$, the smooth bidouble cover lies in both $D$ and $D'$ contradicting its unobstructedness.

\color{black}

\smallskip

The dimension of the moduli component containing surfaces of type $(1)_m$ and $(1')_m$ follows from Proposition ~\ref{moduli11'} (3) and the unobstructedness of ${\bf Def}_X$ for family $(1')_m$. That completes the proof of part (2).\par

Consider now an arbitrary surface of type $(1)_m$ or $(1')_m$ (in particular, $X$ might be singular). By the proof of Proposition ~\ref{lt-any}, we have 
\begin{equation}\label{eq.new.proof}
 0 \to \pi^*(\Omega_{X_1}) \to \Omega_X  \to \pi^*(L^{-1})\otimes \mathcal{O}_R \to 0   
\end{equation}

We apply $\textrm{Hom}(\,\_ \,, \mathscr{O}_X)$ to \eqref{eq.new.proof}. By Remark ~\ref{can} and because $X_1$ is smooth, we have the exact sequence

\begin{center}
\begin{equation}\label{unobcyclic}
0 \to H^0(\pi_1^*(T_{X_1})) \to \textrm{Ext}^1(\pi_1^*(L^{-1})\otimes \mathscr{O}_R, \mathscr{O}_X) \to \textrm{Ext}^1(\Omega_X, \mathscr{O}_X) \to H^1(\pi_1^*(T_{X_1})) \to 0 
\end{equation}
\end{center}

\smallskip
Now note that again from the proof of Proposition ~\ref{lt-any}, we have
$$ 0 \to \pi_1^*(\mathcal{O}_{X_1}(-B)) \to \pi_1^*(L^{-1}) \to \pi_1^*(L^{-1}) \otimes \mathcal{O}_R \to 0. $$
Once again taking $\textrm{Hom}(\, \_, \, \mathscr{O}_X)$ we have
\begin{eqnarray*}
    0 \to \textrm{Hom}(\pi_1^*(L^{-1}), \mathscr{O}_{X}) \to
\textrm{Hom}(\pi_1^*(\mathcal O_{X_1}(-B)), \mathscr{O}_{X}) \to \\ \textrm{Ext}^1(\pi_1^*(L^{-1})\otimes \mathscr{O}_R, \mathscr{O}_X) \to 
\textrm{Ext}^1(\pi_1^*(L^{-1}), \mathscr{O}_X) \to \textrm{Ext}^1(\pi_1^*(\mathcal O_{X_1}(-B)), \mathscr{O}_X) \to 0.
\end{eqnarray*}

Note that $$\textrm{Hom}(\pi_1^*(L^{-1}), \mathscr{O}_{X}) = H^0(\pi_1^*L) = H^0(L) \oplus H^0(\mathscr{O}_{X_1}).$$ Similarly, $$\textrm{Hom}(\pi_1^*\mathscr{O}_{X_1}(-B), \mathscr{O}_{X}) = H^0(\pi_1^*\mathscr{O}_{X_1}(B)) = H^0(\mathscr{O}_{X_1}(B)) \oplus H^0(L).$$ Again $$\textrm{Ext}^1(\pi_1^*(L^{-1}), \mathscr{O}_X) = H^1(\pi_1^*(L)) = H^1(L) \oplus H^1(\mathscr{O}_{X_1}).$$ Similarly $$\textrm{Ext}^1(\pi_1^*\mathscr{O}_{X_1}(B), \mathscr{O}_X) = H^1(\pi_1^*\mathscr{O}_{X_1}(B)) = H^1(\mathscr{O}_{X_1}(B)) \oplus H^1(L).$$ Now using Proposition \ref{moduli11'}, (1) we have $H^1(\mathscr{O}_{X_1}(B)) = 0$. So 
\begin{equation*}
    \textrm{dim}(\textrm{Ext}^1(\pi_1^*(L^{-1})\otimes \mathscr{O}_R, \mathscr{O}_X)) = h^0(\mathscr{O}_{X_1}(B))
\end{equation*}

Putting this back in ~\eqref{unobcyclic} and using Proposition ~\ref{moduli11'}, (1) and (2), we have
\begin{center}
    \begin{equation*}
        \textrm{dim}(\textrm{Ext}^1(\Omega_X, \mathscr{O}_X)) = h^0(\mathscr{O}_{X_1}(B)) = 8m+20
    \end{equation*}
\end{center}
Since this is the dimension of the unique irreducible moduli component of $X$, by \cite{Ser}, Theorem 2.4.1 (iv), this shows that all surfaces of type $(1)_m$ and $(1')_m$  are unobstructed.

\smallskip
\color{black}
(3) Since  there is a unique component of the moduli space containing all surfaces of type $(1)_m$ and $(1')_m$, therefore to describe the canonical morphism of a general surface in this component, it is enough to start with a general surface $X$ of either type. It follows from Corollary ~\ref{11't}, Remark ~\ref{generic degree of canonical morphism}, 
and the fact that $X$ is unobstructed, that for a general surface of the algebraic formally universal deformation space of $X$, the canonical morphism is of degree two onto its image which is non-normal and whose normalization is an elliptic ruled surface which is the projectivization of a split vector bundle over an elliptic curve with invariant $e=0$. Since $X$ is a smooth surface with ample canonical bundle we have that the same holds for its unique irreducible moduli component. 
That completes the proof. \QEDB\par

\begin{remark}\label{nonzerosmooth}
It is interesting to note that $H^1(\mathcal{N}_{\varphi})\neq 0$ by Proposition ~\ref{moduli11'}, but ${\bf Def}_{\varphi}$ is still unobstructed by the above proof. Another example of such an instance is when $\phi: H\to \mathbb{P}^L$ is a morphism that is finite onto its image where $H$ is a hyperk\"ahler variety. It has been proven in \cite{MR}, Lemma 3.1, that in this case ${\bf Def}_{\phi}$ is unobstructed but $H^1(\mathcal{N}_{\phi})\cong H^2(T_H)$ which is non-zero in general.
\end{remark}

\begin{remark}\label{remark.failure.2.4}
From Proposition ~\ref{lt-any} we have the exact sequence
\begin{equation*}
0\to \mathcal{N}_{\pi_1}\to  \pi_1^*\mathcal{O}_Y(B)|_R \to \mathcal{T}_X^1 \to 0.
\end{equation*}
Since $\mathcal{T}_X^1$ is supported on the singular locus, if $X$ is smooth, then $\mathcal{N}_{\pi_1}=  \pi_1^*\mathcal{O}_Y(B)|_R$, while, if $X$ is singular of type $(1_m)$ or $(1'_m)$, then, by \cite{Ser}, \S 3.1.1, $\mathcal{T}_X^1 \neq 0$,  so $\mathcal{N}_{\pi_1}$ and   $\pi_1^*\mathcal{O}_Y(B)|_R$ are not isomorphic in that case. This is why, when proving the smoothness of ${\bf Def}_{X}$,   we proceeded differently for $X$ smooth and for $X$ singular.
\end{remark}
\color{black}

\subsection{Description of moduli components of surfaces of type $(2)_m$} Now we aim to prove Theorem ~\ref{C}. Throughout this subsection, we work with the notations of \S ~\ref{study2}. Recall that $B\in |\mathcal{O}_{\mathbb{P}^1}(6)\otimes\psi^*\mathcal{O}_{\mathbb{P}^1}(2)|$ is the branch divisor of $\pi_2$. In order to do that, we need the following cohomology computations.

\begin{proposition}\label{moduli'2}
Let $X$ be a surface of type $(2)_2$. 
\begin{itemize}
    \item[(1)] $h^0(\mathscr{O}_{X_2}(B)) = 21$ and $h^1(\mathscr{O}_{X_2}(B)) = 0$
    \item[(2)] $h^0(\pi_2^*T_{X_2})=3$ and $h^1(\pi_2^*T_{X_2})=9$. 
\end{itemize}    
Moreover if $X$ is smooth, then
\begin{itemize}
    \item[(3)] $h^0(\mathcal{N}_{\pi_2})=22$ and $h^1(\mathcal{N}_{\pi_2})=0$,
    \item[(4)] $h^1(T_X)=28$,
    \item[(5)] $h^0(\mathcal{N}_{\varphi})=65$.
\end{itemize}
\end{proposition}

\noindent\textit{Proof.} (1) It is easy to see from K\"unneth formula, and projection formula that 
$$h^0(\mathcal{O}_{X_2}(B))=21,\textrm{ and } h^1(\mathcal{O}_{X_2}(B))=0.$$
$(3)$ follows from the following exact sequence 
\begin{equation*}
    0\to\mathcal{O}_{X_2}\to \mathcal{O}_{X_2}(B)\to\mathcal{O}_B(B)\to 0
\end{equation*}
and \cite{Pa}, Corollary 4.1 or \cite{GGP2}, (2.11), and the fact that $h^1(\mathcal{O}_{X_2})=2$ (see Proposition ~\ref{2m} (1)).\par 

(2) This one follows from Proposition ~\ref{2m} (1), K\"unneth formula, and projection formula since $h^j(\pi_2^*T_{X_2})$ is is the following sum
$$h^j(\mathcal{O}_{\mathbb{P}^1}(2)\boxtimes\mathcal{O}_C)+
h^j(\mathcal{O}_{\mathbb{P}^1}\boxtimes\psi^*\mathcal{O}_{\mathbb{P}^1}(-1))+
h^j(\mathcal{O}_{\mathbb{P}^1}(-1)\boxtimes\psi^*\mathcal{O}_{\mathbb{P}^1}(-1))+
h^j(\mathcal{O}_{\mathbb{P}^1}(-3)\boxtimes\psi^*\mathcal{O}_{\mathbb{P}^1}(-2)).$$
\indent (4) We use the following two exact sequence 
\begin{equation}
    0\to T_X\to \pi_2^*T_{X_2}\to\mathcal{N}_{\pi_2}\to 0.
\end{equation}
Since $h^0(T_X)=0$, and $h^1(\mathcal{N}_{\pi_2})=0$ by part (1), we get the following exact sequence:
$$0\to H^0(\pi_2^*T_{X_2})\to H^0(\mathcal{N}_{\pi_2})\to H^1(T_X)\to H^1(\pi_2^*T_{X_2})\to 0.$$
The conclusion now follows from part (2).\par 

(5) We get the following exact sequence from Lemma ~\ref{exactgonzalez}:
\begin{equation}\label{exactgonzalezpi2}
    0\to \mathcal{N}_{\pi_2}\to\mathcal{N}_{\varphi}\to\pi_2^*\mathcal{N}_{\varphi_2}\to 0.
\end{equation}

We first compute $h^0(\mathcal{N}_{\varphi})$. From part (1), we get $h^1(\mathcal{N}_{\pi_2})=0$. It follows from projection formula that $$h^0(\mathcal{N}_{\varphi})=h^0(\mathcal{N}_{\pi_2})+h^0(\mathcal{N}_{\varphi_2})+h^0(\mathcal{N}_{\varphi_2}\otimes\mathcal{E}_{\pi_2}).$$
The vanishing of $h^0(\mathcal{N}_{\varphi_2}\otimes\mathcal{E}_{\pi_2})$ has been shown in the proof of Corollary ~\ref{2t}. Thus, $$h^0(\mathcal{N}_{\varphi})=22+36+7=65$$ by part (1) and Proposition ~\ref{2m} (5).\QEDB\par

\vspace{5pt}

\noindent\textit{Proof of Theorem ~\ref{C}.} We only need to show  the unobstructedness of $X$ for $X$ of type $(2)_2$. All the cohomological criteria have been verified in Proposition ~\ref{2m}, Corollary ~\ref{2t} and in Proposition ~\ref{moduli'2}. The  unobstructedness of $\varphi$ in assertion 1 follows from Corollary~\ref{Defphi} and assertion 1 (b) and (2) are consequences of Corollary ~\ref{2t}, Remark ~\ref{generic degree of canonical morphism} and the unobstructedness of $X$. 

\smallskip
We first show that for a smooth
surface $X$ of type $(2)_2$, $X$ is unobstructed. We need to show any section of $H^0(\omega_X)$ lifts to any first order deformation of $X$. The remaining argument is identical to the proof of Theorem ~\ref{B}. %\par
Using the generalized Atiyah sequence ~\eqref{mainatiyah}, and arguing as in the proof of Theorem ~\ref{B}, we obtain
$$h^0(\mathcal{E}_{\omega_X})=1 \textrm{ and }h^1(\mathcal{E}_{\omega_X})=h^1(\mathcal{O}_X)+h^1(T_X)=30$$ thanks to Proposition ~\ref{moduli'2} (3). Now, arguing as in the proof of Theorem ~\ref{B}, we get that the dimension of the image of $H^0(\mathcal{N}_{\varphi})\to H^1(\mathcal{E}_{\omega_X})$ is 
$$ h^0(\mathcal{E}_{\omega_X})-(h^0(\omega_X))^2+h^0(\mathcal{N}_{\varphi})=1-36+65=30$$ thanks to Proposition ~\ref{moduli'2} (4). Thus, the map $H^1(\mathcal{E}_{\omega_X})\to H^0(\omega_X)^{\vee}\otimes H^1(\omega_X)$ is zero, so any section of $H^0(\omega_X)$ lifts to any first order deformation of $X$. \par

Consider now an arbitrary surface $X$ of type $(2)_2$ (in particular, $X$ might be singular). By the proof of Proposition ~\ref{lt-any}, as in the proof of Theorem ~\ref{B}, we have 

\begin{center}
\begin{equation}\label{unob2}
0 \to H^0(\pi_2^*(T_{X_2})) \to \textrm{Ext}^1(\pi_2^*(L^{-1})\otimes \mathscr{O}_R, \mathscr{O}_X) \to \textrm{Ext}^1(\Omega_X, \mathscr{O}_X) \to H^1(\pi_2^*(T_{X_2})) \to 0 
\end{equation}
\end{center}

and $$ 0 \to \pi_2^*(\mathcal{O}_{X_2}(-B)) \to \pi_2^*(L^{-1}) \to \pi_2^*(L^{-1}) \otimes \mathcal{O}_R \to 0. $$ Calculating
$\textrm{Ext}^1(\Omega_X, \mathscr{O}_X)$ from the above exact sequences as in the proof of Theorem ~\ref{B}, we have
$$ \textrm{Ext}^1(\Omega_X, \mathscr{O}_X) = 28 $$ Then $X$ is unobstructed by \cite{Ser}, Theorem 2.4.1 (iv) and the fact that surfaces of type $(2)_2$ form an irreducible family (the latter follows in the same way as in Step $(1)$ of the proof of Theorem ~\ref{B}, part $(2)$).
\QEDB\par

\vspace{5pt}
\color{black}
We end this section by asking the following natural questions, concerning the deformations of surfaces of type $(2)_m$, as it is evident that our technique of showing unobstructedness does not work for them.

\begin{question}\label{question1}
Let $X$ be a smooth surface of type $(2)_m$.
\begin{itemize}
     \item[(1)] Are ${\bf Def}_{\varphi}$ and ${\bf Def}_X$ unobstructed if $m\geq 3$? The problem that we face for these surfaces is that $H^1(\mathcal{N}_{\pi_2})\neq 0$. Thus, for these surfaces, we know that the forgetful map ${\bf Def}_{\pi_2/\mathbb{P}^N}\to{\bf Def}_{\varphi}$ is smooth; however the smoothness of ${\bf Def}_{\pi_2/\mathbb{P}^N}$ is unknown despite knowing the smoothness of ${\bf Def}_{\varphi_2}$, since the forgetful map ${\bf Def}_{\pi_2/\mathbb{P}^N}\to {\bf Def}_{\varphi_2}$ is not smooth.
    \item[(2)] Can surfaces of type $(2)_m$, for $m \geq 3$ be deformed to canonical double covers over surfaces ruled over smooth curves of genus $m$ ?
\end{itemize} 
\end{question}

\subsection{Dimensions of the quadruple loci}

For $X$ of type $(1)_m$, $(1')_m$, $(2)_2$ or $(3)_m$, we define $\mathcal M_{[X]}^{\textrm{quad}}$
as the locally closed sublocus, in the (unique) irreducible moduli component $\mathcal M_{[X]}$ of $X$, that parametrizes surfaces of general type whose canonical morphism is a quadruple cover of a surface of minimal degree. As pointed out in the introduction, since $\mathcal M_{[X]}^{\textrm{quad}}$ is, except in case $(3)_1$, a proper sublocus, quadruple canonical covers resemble hyperelliptic curves and the existence of a proper hyperelliptic sublocus in the moduli of curves of genus bigger than $2$.
 Thus, it is natural to ask for the dimension of $\mathcal M_{[X]}^{\textrm{quad}}$ and of
other related subloci contained in $\mathcal M_{[X]}^{\textrm{quad}}$, namely, 
the locally closed sublocus $\mathcal{M}_{[X]}^G$ of surfaces whose
canonical morphism is a quadruple Galois cover with the same Galois group $G$ as $X$ and the locally closed sublocus $\mathcal{M}_{[X]}^{\textrm{Nat}}(G)$ that contains $\mathcal{M}_{[X]}^G$ and parametrizes the so-called natural deformations of the cover $\pi: X \longrightarrow Y$ (see \cite{Pa}, Definition $5.1$; the canonical morphism of surfaces of $\mathcal{M}_{[X]}^{\textrm{Nat}}(G)$ is again four--to--one but not Galois in general). Note that for the  moduli of curves and its hyperelliptic locus the 
three analogous subloci coincide as double covers are Galois. In the next remark
we compute the dimensions of $\mathcal{M}_{[X]}^G$ and $\mathcal{M}_{[X]}^{\textrm{Nat}}(G)$ (this yields a lower bound for the dimension 
of $\mathcal M_{[X]}^{\textrm{quad}}$). We omit the proof as it follows from 
the description of $\mathcal{M}_{[X]}^{\textrm{Nat}}(G)$ given in \cite{Pa}, Section 5, and general principles of deformation theory under some vanishing conditions, using methods same as in the proof of Theorem ~\ref{Def_1}.

\color{black}

\begin{remark}\label{remark.dimensions.loci}
The loci $\mathcal{M}_{[X]}^G$ and $\mathcal{M}_{[X]}^{\textrm{Nat}}(G)$ are uniruled and their dimensions are as stated in the following table, where we compare them with the dimension, computed in Theorems~\ref{B} and \ref{C} and Remark~\ref{remark.Ops}, of the irreducible moduli component $\mathcal M_{[X]}$ of $X$, that contains them. We exclude the case $(3)_1$ from the table as in this case ${\mathcal{M}_{[X]}^G=\mathcal M_{[X]}}$. 

\smallskip

\begin{center}
\begin{tabular}{c|c|c|c|c} 
\hline
 \makecell{$X$ is general \\ of type} & \makecell{Dim. of $\mathcal{M}_{[X]}^G$} & \makecell{Dim. of $\mathcal{M}_{[X]}^{\textrm{Nat}}(G)$}   
 & \makecell{Lower bound for \\
 dim. of $\mathcal{M}_{[X]}^{\textrm{quad}}$}
 & \makecell{Dim. of $\mathcal{M}_{[X]}$}  \\  
 \hline\hline 
 $(1)_m$ & \makecell{$2m+1$} &  \makecell{$2m+4$} & \makecell{$6m+18$}   & \makecell{$8m+20$}\\
 \hline 
 $(1')_m$ & \makecell{$6m+13$} & \makecell{$6m+18$} &  \makecell{$6m+18$}  & \makecell{$8m+20$} \\
 \hline
$(2)_2$ & \makecell{$25$} &  \makecell{$25$} & \makecell{$25$} & \makecell{$28$} \\
 \hline
\makecell{$(3)_m$ ($m \geq 2$)} & \makecell{$2m+4$} & \makecell{$2m+4$}  &  \makecell{$2m+4$}  & \makecell{$3m+3$} \\ 
 \hline
\end{tabular}
\end{center}

\end{remark}

\color{black}
\begin{remark}\label{remark.subloci.product.of.curves}
As mentioned before Remark~\ref{remark.Ops}, the moduli space in which surfaces of type $(3)_m$ lie is much easier compared to the other types $(1)_m$, $(1')_m$ and $(2)_m$, due to the fact that 
they are products of curves. Then the existence, when $m \geq 2$, of a quadruple locus in these moduli spaces, follows essentially from 
\cite{BeauD} and \cite[Theorems 4.1 and 4.2]{Ops} (see Corollary~\ref{cor.Beauville}). In fact, products of curves yield another, easy-to-get, examples of irreducible moduli components with proper, locally closed subloci, where the degree of the canonical morphisms jumps up. Indeed, let $g_1, g_2 \geq 3$. Products of curves of genera $g_1$ and $g_2$ are parametrized by an irreducible moduli component $\mathcal M_{g_1,g_2}$ which has three strata: the general one, corresponding to surfaces whose canonical morphism is an embedding; a stratum parametrizing surfaces whose canonical morphism has degree $2$ onto its image; a stratum parametrizing surfaces whose canonical morphism has degree $4$ onto its image. Note however that the image of the canonical morphism of the surfaces of this last stratum is not a surface of minimal degree, as the canonical covers this article is concerned with. Coming back to our surfaces of types $(1)_m$, $(1')_m$ and $(2)_m$ versus surfaces of type $(3)_m$, it is interesting that, while the existence of the quadruple locus for the latter is dictated by being a product, the former, despite not being such, have also a quadruple locus of the same nature.
\end{remark}

\color{black}
\subsection{Remarks on the geography of surfaces of general type}

Recall that the invariants of a minimal surface of general type with
birational canonical map satisfy Castelnuovo's inequality, which, for given irregularity $q$ is 
\begin{equation*}\label{Castelnuovo.inequality}
    K^2 \geq p_g + q -7.
\end{equation*}
It's worth noting that Theorems~\ref{B}, \ref{C} and Remark~\ref{remark.Ops} imply that, for any $m \geq 1$, the moduli spaces $\mathcal M_{(p_g,q,K^2)}=\mathcal M_{(2m+2,1,8m)}$, $\mathcal M_{(2m+2,m,8m)}$ and $\mathcal M_{(2m+2,m+3,8m)}$, have  an irreducible 
component that pa\-ra\-me\-trizes surfaces whose canonical map is a non birational morphism, despite the fact that, except for 
$\mathcal M_{(4,4,8)}$, they are all above Castelnuovo's line. This and 
Ashikaga's construction of surfaces with birational canonical map 
(see \cite[Theorem 3.2]{Ashikaga}) yield the existence of infinitely many moduli spaces having at least two irreducible components of very distinct nature:

\begin{corollary}\label{corollary.two.components}
For any $m \geq 1$, the moduli 
space $\mathcal M_{(2m+2,1,8m)}$ has one irreducuble component parametrizing surfaces with non birational canonical morphisms and another component whose general elements have birational canonical maps.
\end{corollary}

\section{Infinitesimal Torelli theorem}\label{5} 
{The goal of this section is to prove the infinitesimal Torelli theorem for some smooth families of quadruple covers.}
Let $X$ be a smooth algebraic variety of dimension $n$ with ample canonical bundle $\omega_X=\Omega_X^n$. Let $p:\mathcal{X}\to U_X$ be a semiuniversal deformation of $X$, and we assume that $S$ is smooth. The infinitesimal Torelli problem for weight $n$ Hodge structure asks how far the complex structure of $X$ is determined by the decreasing Hodge filtration $$\left(F^p=\bigoplus_{i\geq p}H^{n-i}(X,\Omega_X^i)\right)_{p\in\mathbb{N}}.$$
The Hodge filtrations on the fibres glue together to give a subbundle $\mathcal{F}^p$ of $\mathcal{O}_{U_X}\otimes R^np_*\mathbb{Z}_{\mathcal{X}}$ and define the period map $\Phi_n:U_X\to D_n$ to the space $D_n$ parametrizing Hodge filtrations of weight $n$. The tangent map $T\Phi_n$ at the special point is the composition of an injective map and the sum of the linear maps 
$$\lambda_{i}: H^1(X, T_X)\to \textrm{Hom}(\mathbb{C}; H^{n-i}(X,\Omega_X^i), H^{n-i+1}(X,\Omega_X^{i-1})),\quad n\geq i\geq 1$$
induced by the contraction maps. Hence $T\Phi_n$ is injective if $\lambda_i$ is injective for some $i$, in which case $\Phi_n$ is an immersion. We say that 
\begin{itemize}[leftmargin=0.7cm]
    \item[(1)] the infinitesimal Torelli theorem holds (for weight $n$ Hodge structures) for $X$ if $\Phi_n$ is an immersion;
    \item[(2)] the infinitesimal Torelli theorem for periods of $n$ forms holds for $X$ if $\lambda_n$ is an injection.
\end{itemize}   
The infinitesimal Torelli theorem for periods of $n$ forms holds for $X$ if and only if the following map is a surjection
$$H^{n-1}(X,\Omega_X)\otimes H^0(X,\Omega_X^n)\to H^{n-1}(\Omega_X\otimes \Omega_X^n).$$

Classically, a curve of genus $g\geq 2$ satisfies the infinitesimal Torelli theorem if and only if $g=2$ or it is non-hyperelliptic. When $n=2$ i.e., when $X$ is a surface, the infinitesimal Torelli theorem for periods of 2 forms holds for $X$ if and only if the infinitesimal Torelli theorem for weight 2 Hodge structures holds for $X$. The infinitesimal Torelli problem for abelian covers was studied by Pardini in \cite{Par98} in a very general setting. We refer to the article of Catanese (see \cite{Cat84}) for counterexamples of Torelli problems, the article of Bauer and Catanese (see \cite{BC}) for counterexamples of the infinitesimal Torelli theorems with $\omega_X$ quasi-very ample.\par

The following criterion under which the infinitesimal Torelli theorem for weight $n$ Hodge structures holds for $X$ was developed by Flenner.
\begin{theorem}\label{flenner} (\cite{Fle}, Theorem 1.1)
Let $X$ be a compact $n$-dimensional K\"ahler manifold and assume the existence of a resolution of $\Omega_X$ by vector bundles 
$$0\to\mathcal{G}\to\mathcal{F}\to\Omega_X\to 0.$$
If both conditions are satisfied:
\begin{itemize}
    \item[(a)] $H^{j+1}(S^j\mathcal{G}\otimes \bigwedge^{n-j-1}\mathcal{F}\otimes\omega_X^{-1})=0$ for all $0\leq j\leq n-2$;
    \item[(b)] The pairing $H^0(S^{n-p}\mathcal{G}^{-1}\otimes \omega_X)\otimes H^0(S^{p-1}\mathcal{G}^{-1}\otimes\omega_X)\to H^0(S^{n-1}\mathcal{G}^{-1}\otimes\omega_X^{\otimes 2})$ is surjective for a suitable $p\in\{1,2,\dots,n\}$ 
\end{itemize}
then the canonical map $\lambda_p: H^1(X,T_X)\to \textrm{Hom}(\mathbb{C};H^{n-p}(X,\Omega_X^p),H^{n+1-p}(X,\Omega_X^{p-1}))$ is injective.
\end{theorem}

Notice that if $X$ is a smooth quadruple Galois canonical cover of a smooth surface of minimal degree with irregularity one, then $X$ is necessarily of type $(1')_m$. In order to prove the theorem, we are going to invoke the theorem of Flenner i.e., Theorem ~\ref{flenner}.

\begin{theorem}\label{D}
Let $X$ be a smooth surface with irregularity one. Assume the canonical bundle $\omega_X$ is ample and globally generated, and the canonical morphism $\varphi$ is a quadruple Galois canonical cover onto a smooth surface of minimal degree. Then the infinitesimal Torelli theorem holds for $X$.
\end{theorem}

\noindent\textit{Proof.} We work with the notations of \S ~\ref{study11'}. We have the following commutative diagram
\[
\begin{tikzcd}
X \arrow[hookrightarrow, "j"]{r} \arrow[swap, "\pi_1"]{rd}
& Z \arrow["p"]{d}\\
& X_1=E\times\mathbb{P}^1
\end{tikzcd}
\]
where $Z=\mathbb{V}(\mathcal{E}_{\pi_1})$ is the affine bundle over $X_1$ with the natural projection $p$. We have the following short exact sequence ~\eqref{ses1}
\begin{equation*}
    0\to\pi_1^*\mathcal{E}_{\pi_1}^{\otimes 2}\to\Omega_Z|_X\to\Omega_X\to 0.
\end{equation*}
We will use the criterion of Flenner (Theorem ~\ref{flenner}). \par 
(a) Since dim$(X)=2$, we need to check $H^1(\Omega_Z|_X\otimes\omega_X^{-1})=0$. We also have the following exact sequence ~\eqref{ses2}
\begin{equation}\label{flennerexact}
    0\to \pi_1^*\Omega_{X_1}\to \Omega_Z|_X\to \pi_1^*\mathcal{E}_{\pi_1}\to 0.
\end{equation}
Now, $\Omega_{X_1}=(\mathcal{O}_E\boxtimes\mathcal{O}_{\mathbb{P}^1}(-2))\oplus(\mathcal{O}_E\boxtimes\mathcal{O}_{\mathbb{P}^1})$. Consequently, we obtain:
\begin{equation*} 
H^1(\pi_1^*\Omega_{X_1}\otimes\omega_X^{-1})= H^1(\pi_1^*(\Omega_{X_1}\otimes(\psi^*\mathcal{O}_{\mathbb{P}^1}(-1)\boxtimes\mathcal{O}_{\mathbb{P}^1}(-m))))
\end{equation*}
Now, the last term by projection formula is just
\begin{align*}
    & H^1(\psi^*\mathcal{O}_{\mathbb{P}^1}(-1)\boxtimes\mathcal{O}_{\mathbb{P}^1}(-m-2))\oplus H^1(\psi^*\mathcal{O}_{\mathbb{P}^1}(-1)\boxtimes\mathcal{O}_{\mathbb{P}^1}(-m))\oplus\\ & H^1(\psi^*\mathcal{O}_{\mathbb{P}^1}(-2)\boxtimes\mathcal{O}_{\mathbb{P}^1}(-2m-4)) \oplus H^1(\psi^*\mathcal{O}_{\mathbb{P}^1}(-2)\boxtimes\mathcal{O}_{\mathbb{P}^1}(-2m-2)).
\end{align*}
This term is zero by K\"unneth formula and the projection formula. Notice that we also have
\begin{equation*}
H^1(\pi_1^*\mathcal{E}_{\pi_1}\otimes\omega_X^{-1}) = H^1(\pi_1^*(\psi^*\mathcal{O}_{\mathbb{P}^1}(-2)\boxtimes\mathcal{O}_{\mathbb{P}^1}(-2m-2)))
\end{equation*}
This term is just $H^1(\psi^*\mathcal{O}_{\mathbb{P}^1}(-2)\boxtimes\mathcal{O}_{\mathbb{P}^1}(-2m-2))\oplus H^1(\psi^*\mathcal{O}_{\mathbb{P}^1}(-3)\boxtimes\mathcal{O}_{\mathbb{P}^1}(-3m-4))$ which is zero by K\"unneth formula and the projection formula. Thus, it follows from ~\eqref{flennerexact} that $H^1(\Omega_Z^1|_X\otimes\omega_X^{-1})=0$.\par 

(b) Now we check the surjection corresponding to $p=1$, i.e.,
$$H^0(\pi_1^*(\mathcal{E}_{\pi_1}^*)^{\otimes 2}\otimes\omega_X)\otimes H^0(\omega_X)\to H^0(\pi_1^*(\mathcal{E}_{\pi_1}^*)^{\otimes 2}\otimes\omega_X^{\otimes 2}).$$
Thus we need to check the surjection of 
$$H^0(\pi_1^*(\psi^*\mathcal{O}_{\mathbb{P}^1}(2)\boxtimes\mathcal{O}_{\mathbb{P}^1}(2m+4))\otimes\omega_X)\otimes H^0(\omega_X)\to H^0(\pi_1^*(\psi^*\mathcal{O}_{\mathbb{P}^1}(2)\boxtimes\mathcal{O}_{\mathbb{P}^1}(2m+4))\otimes\omega_X^{\otimes 2}).$$
To check this surjection, we use Castelnuovo-Mumford regularity (see \cite{Mu}). By projection formula, $$H^1(\pi_1^*(\psi^*\mathcal{O}_{\mathbb{P}^1}(2)\boxtimes\mathcal{O}_{\mathbb{P}^1}(2m+4)))=H^1(\psi^*\mathcal{O}_{\mathbb{P}^1}(2)\boxtimes\mathcal{O}_{\mathbb{P}^1}(2m+4))\oplus H^1(\psi^*\mathcal{O}_{\mathbb{P}^1}(1)\boxtimes\mathcal{O}_{\mathbb{P}^1}(m+2))$$
and it is easy to check that both terms are zero by K\"unneth formula and the projection formula.
Now we compute the following cohomology group  $$H^2(\pi_1^*(\psi^*\mathcal{O}_{\mathbb{P}^1}(2)\boxtimes\mathcal{O}_{\mathbb{P}^1}(2m+4))\otimes\omega_X^{-1})=H^2(\psi^*\mathcal{O}_{\mathbb{P}^1}(1)\boxtimes\mathcal{O}_{\mathbb{P}^1}(m+4))\oplus H^2(\psi^*\mathcal{O}_{\mathbb{P}^1}\boxtimes\mathcal{O}_{\mathbb{P}^1}(2))$$
and one checks that both terms are zero. That concludes the proof.\QEDB\par

\begin{remark}\label{producttorelli}
Let $X$ be a surface of type $(3)_m$. It is easy to see in these cases that the infinitesimal Torelli theorem holds only if $X$ is a surface of type $(3)_1$. We give a brief explanation following \cite{BC} for the sake of completeness. We resume the notations of \S ~\ref{defproduct}. Since $\textrm{dim}(X)=2$, infinitesimal Torelli theorem holds $\iff$ infinitesimal Torelli theorem for periods of 2 forms holds $\iff$ $H^1(\Omega_X^1)\otimes H^0(\Omega_X^2)\to H^1(\Omega_X^1\otimes\Omega_X^2)$ is surjective.\par 
Using $\Omega_{X}^1=(\Omega_{C_1}^1\boxtimes\mathcal{O}_{C_2})\oplus(\mathcal{O}_{C_1}^1\boxtimes\Omega_{C_2})$, $\Omega_X^2=\Omega_{C_1}^1\boxtimes\Omega_{C_2}^1$, and K\"unneth formula, we have
\begin{equation*}
    H^1(\Omega_X^1)=\left(H^0(\Omega_{C_1}^1)\otimes H^1(\mathcal{O}_{C_2})\right)\oplus \left(H^1(\mathcal{O}_{C_1})\otimes H^0(\Omega_{C_2}^1)\right),
\end{equation*}
\begin{equation*}
    H^0(\Omega_X^2)=H^0(\Omega_{C_1}^1)\otimes H^0(\Omega_{C_2}^1),
\end{equation*}
\begin{equation*}
    H^1(\Omega_X^1\otimes\Omega_X^2)=H^1((\Omega_{C_1}^1)^{\otimes 2}\boxtimes \Omega_{C_2}^1)\oplus H^1(\Omega_{C_1}^1\boxtimes(\Omega_{C_2}^1)^{\otimes 2})=\left(H^0((\Omega_{C_1}^1)^{\otimes 2})\otimes H^1(\Omega_{C_2}^1)\right)\oplus\left(H^1(\Omega_{C_1}^1)\times H^0((\Omega_{C_2}^1)^{\otimes 2})\right),
\end{equation*}
where the last equality is obtained by using the fact that $H^1((\Omega_{C_j}^1)^{\otimes 2})=0$ for $j=1,2$. Notice also that $H^1(\mathcal{O}_{C_j}^1)\otimes H^0(\Omega_{C_j}^1)\to H^1(\Omega_{C_j}^1)$ is surjective for $j=1,2$ by the non-degeneracy of Serre duality. Thus, $H^1(\Omega_X^1)\otimes H^0(\Omega_X^2)\to H^1(\Omega_X^1\otimes\Omega_X^2)$ is surjective if and only if 
\begin{equation}\label{surjhyper}
    H^0(\Omega_{C_j}^1)\times H^0(\Omega_{C_j}^1)\to H^0(\Omega_{C_j}^{\otimes 2})
\end{equation}
is surjective for $j=1,2$. Notice that ~\eqref{surjhyper} is surjective for $j=1$ since $C_1$ is a hyperelliptic curve of genus 2. Since $C_2$ is also hyperelliptic, ~\eqref{surjhyper} is surjective only when $m=1$. 
\end{remark}

We end this article by asking the natural question regarding the infinitesimal Torelli theorem for smooth surfaces of type $(2)_m$.

\begin{question}
Let $X$ be a smooth surface of type $(2)_m$. Does the infinitesimal Torelli theorem hold for $X$? Let us resume the notations of \S ~\ref{study2}. Let $Z:=\mathbb{V}(\mathcal{E}_{\pi_2})$ be the affine bundle over $X_2$. It is easy to verify that $H^1(\Omega^1_Z|_X\otimes\omega_X^{-1})=0$. However, to apply the criterion of Flenner (i.e., Theorem ~\ref{flenner}), we need the surjectivity of the following multiplication map:
$$H^0(\pi_2^*(\mathcal{E}_{\pi_2}^*)^{\otimes 2}\otimes\omega_X)\otimes H^0(\omega_X)\to H^0(\pi_2^*(\mathcal{E}_{\pi_2}^*)^{\otimes 2}\otimes\omega_X^{\otimes 2}).$$
It follows from \cite{GP01}, Lemma 2.1 that this map is not surjective. 
\end{question}

%\smallskip

\section*{Funding and/or Conflicts of interests/Competing interests.} 

 The first author was partially supported
by the General Research Fund of the University of Kansas.
The second author was partially supported
by Spanish Government grant PID2021-124440NB-I00 and by Santander-UCM grant PR44/21. The third author was supported by the National Science Foundation, Grant No. DMS-1929284 while in residence at ICERM in Providence, RI, as part of the ICERM Bridge program. The research of the fourth author was partially supported by a Simons postdoctoral fellowship from the Fields Institute for Research in Mathematical Sciences.

\smallskip

We have no conflict of interests to disclose.

%\smallskip

\bibliographystyle{plain}

\end{document}